\documentclass[reqno]{amsart}
\usepackage[margin = 1in]{geometry}
\usepackage{amsmath, amssymb, amsthm, fancyhdr, verbatim, graphicx, mathtools}
\usepackage{enumerate}
\usepackage[all]{xy}
\usepackage[usenames,dvipsnames]{xcolor}
\usepackage{mathrsfs}
\usepackage{tikz-cd}
\usepackage[T1]{fontenc} 

\usepackage{framed, hyperref}	
\usepackage[OT2, T1]{fontenc}
\usepackage[titletoc]{appendix}
\usepackage{bbm}
\usepackage{enumitem}

\usepackage{amssymb}

\usepackage{stmaryrd}
\newcommand{\shear}{\mathbin{\mkern-6mu\fatslash}}

\DeclareSymbolFont{cyrletters}{OT2}{wncyr}{m}{n}
\DeclareMathSymbol{\Sha}{\mathalpha}{cyrletters}{"58}

\usepackage{color}

\newcommand{\F}{\mathbf{F}}

\newcommand{\CC}{\mathbf{C}}
\newcommand{\G}{\mathbf{G}}

\newcommand{\wt}[1]{\widetilde{#1}}

\newcommand{\Q}{\mathbf{Q}}
\newcommand{\Z}{\mathbf{Z}}

\newcommand{\mf}[1]{\mathfrak{#1}}

\newcommand{\Gal}{\operatorname{Gal}}

\newcommand{\cl}{\overline}

\newcommand{\ul}[1]{\underline{#1}}
\newcommand{\ol}[1]{\overline{#1}}
\newcommand{\wh}[1]{\widehat{#1}}

\newcommand{\Cal}[1]{\mathcal{#1}}
\newcommand{\A}{\mathbf{A}}

\newcommand{\co}{\colon}
\newcommand{\mrm}[1]{\mathrm{#1}}
\newcommand{\msf}[1]{\mathsf{#1}}

\newcommand{\bs}{\backslash}

\newcommand{\PP}{\mathbf{P}}

\newcommand{\bu}{\bullet} 
\newcommand{\inj}{\hookrightarrow}
\newcommand{\sotimes}{\stackrel{!}\otimes}
\newcommand{\diagotimes}{\stackrel{\Delta}\otimes}

\newcommand{\bL}{\mathbf{L}}
\newcommand{\bT}{\mathbf{T}}

\newcommand{\chX}{\check{X}}
\newcommand{\chY}{\check{Y}}

\newcommand{\chG}{\check{G}}
\newcommand{\chT}{\check{T}}

\newcommand{\chB}{\check{B}}
\newcommand{\chM}{\check{M}}

\newcommand{\chP}{\check{P}}
\newcommand{\chH}{\check{H}}

\newcommand{\chlambda}{\check{\lambda}}
\newcommand{\BunX}{\Bun^X_{G}}
\newcommand{\LocX}{\Loc^{\check{X}}_{\check{G}}}
\newcommand{\LocG}{\Loc_{\check{G}}}
\newcommand{\LocY}{\Loc^{Y}_{\check{G}}}
\newcommand{\LocH}{\Loc_{\check{H}}}
\newcommand{\LocB}{\Loc_{\check{B}}}

\newcommand{\LS}{\Loc_{\check{G}}}
\newcommand{\PX}{\cP_X}
\newcommand{\LX}{\cL_{\chX}}
\newcommand{\pispec}{\pi^{\mrm{spec}}}
\newcommand{\piaut}{\pi^{\mrm{aut}}}
\newcommand{\LL}{\mathbb{L}}
\newcommand{\Ggr}{\mathbb{G}_{gr}}

\newcommand\cA{\mathcal{A}}
\newcommand\cB{\mathcal{B}}

\newcommand\cD{\mathcal{D}}
\newcommand\cE{\mathcal{E}}
\newcommand\cF{\mathcal{F}}
\newcommand\cG{\mathcal{G}}

\newcommand\cI{\mathcal{I}}

\newcommand\cK{\mathcal{K}}
\newcommand\cL{\mathcal{L}}

\newcommand\cO{\mathcal{O}}
\newcommand\cP{\mathcal{P}}

\newcommand\cV{\mathcal{V}}
\newcommand\cW{\mathcal{W}}
\newcommand\cY{\mathcal{Y}}
\newcommand\cZ{\mathcal{Z}}

\RequirePackage{xspace}

\newcommand{\rB}{\ensuremath{\mathrm{B}}\xspace}

\newcommand{\rH}{\ensuremath{\mathrm{H}}\xspace}

\newcommand{\rR}{\ensuremath{\mathrm{R}}\xspace}

\newcommand{\sA}{\ensuremath{\mathscr{A}}\xspace}

\newcommand{\sW}{\ensuremath{\mathscr{W}}\xspace}

\DeclareMathOperator{\GL}{GL}

\DeclareMathOperator{\Frob}{Frob}

\DeclareMathOperator{\coker}{coker}

\DeclareMathOperator{\Hom}{Hom}
\DeclareMathOperator{\Ind}{Ind}

\DeclareMathOperator{\rank}{rank}

\DeclareMathOperator{\Rep}{Rep}

\DeclareMathOperator{\Spec}{Spec}

\DeclareMathOperator{\Lie}{Lie}

\DeclareMathOperator{\red}{red}

\DeclareMathOperator{\Div}{Div}

\DeclareMathOperator{\Bun}{Bun}
\DeclareMathOperator{\Ext}{Ext}

\DeclareMathOperator{\ch}{ch\,}
\DeclareMathOperator{\Id}{Id}

\DeclareMathOperator{\Sym}{Sym}

\DeclareMathOperator{\pt}{pt}

\DeclareMathOperator{\Sat}{Sat}

\DeclareMathOperator{\pr}{pr}

\DeclareMathOperator{\Dmod}{Dmod}
\DeclareMathOperator{\IndCoh}{IndCoh}
\DeclareMathOperator{\Loc}{Loc}
\DeclareMathOperator{\Map}{Map}
\DeclareMathOperator{\dR}{dR}
\DeclareMathOperator{\Ran}{Ran}
\DeclareMathOperator{\QCoh}{QCoh}
\DeclareMathOperator{\ev}{ev}
\DeclareMathOperator{\Hk}{Hk}
\DeclareMathOperator{\Std}{Std}

\DeclareMathOperator{\Eis}{Eis}
\DeclareMathOperator{\spec}{spec}

\DeclareMathOperator{\Mir}{Mir}
\newcommand{\univ}{\mathrm{univ}}
\DeclareMathOperator{\triv}{triv}
\DeclareMathOperator{\Nilp}{Nilp}

\DeclareMathOperator{\colim}{colim}
\DeclareMathOperator{\FinSurj}{FinSurj}
\DeclareMathOperator{\Chev}{Chev}

\DeclareMathOperator{\Perf}{Perf}
\DeclareMathOperator{\aut}{aut}
\DeclareMathOperator{\Coh}{Coh}

\DeclareMathOperator{\op}{op}
\DeclareMathOperator{\disj}{disj}
\DeclareMathOperator{\Fact}{Fact}
\DeclareMathOperator{\Vect}{Vect}
\DeclareMathOperator{\add}{add}
\DeclareMathOperator{\AJ}{AJ}
\DeclareMathOperator{\Mod}{-Mod}
\DeclareMathOperator{\Tot}{Tot}
\DeclareMathOperator{\FT}{FT}

\DeclareMathOperator{\hol}{hol}
\DeclareMathOperator{\supp}{supp}

\DeclareMathOperator{\Shv}{Shv}

\DeclareMathOperator{\CAlg}{ComAlg}

\DeclareMathOperator{\CoCAlg}{ComCoalg}
\DeclareMathOperator{\FS}{FinSurj}
\DeclareMathOperator{\un}{union}
\DeclareMathOperator{\oblv}{oblv}

\DeclareMathOperator{\Sect}{Sect}

\DeclareMathOperator{\fib}{fib}

\newcommand{\AS}{\sA}
\DeclareMathOperator{\norm}{norm}

\newcommand{\sm}[1]{{\scriptstyle  #1}}

\newcommand{\Hdr}{\Gamma_{\dR}}

\newcommand{\tw}[1]{\sm{\langle #1 \rangle}}
\newcommand{\ulk}{\underline{k}}

\newcommand{\cRHom}{\Cal{R}\Cal{H}om}
\DeclareMathOperator{\RHom}{RHom}

\newcommand{\surj}{\twoheadrightarrow}
\newcommand{\hSect}{\Sect}

\newcommand{\otstar}{\otimes^{\star}}
\newcommand{\otch}{\otimes^{\mrm{ch}}}
\newcommand{\Lambdapos}{\Lambda^{+}}

\renewcommand{\cl}{\mrm{cl}}
\renewcommand{\ch}{\mrm{ch}}

\newtheorem{thm}{Theorem}[subsection]
\newtheorem{lemma}[thm]{Lemma}
\newtheorem{lemmaconstr}[thm]{Lemma/Construction}
\newtheorem{prop}[thm]{Proposition}
\newtheorem{cor}[thm]{Corollary}

\newtheorem{conj}[thm]{Conjecture}

\theoremstyle{remark}
\newtheorem{remark}[thm]{Remark} 
\newtheorem{notat}[thm]{Notation} 
\newtheorem{constr}[thm]{Construction} 
\newtheorem{defn}[thm]{Definition}

\newtheorem{example}[thm]{Example}

\makeatletter
\def\th@remark{%
  \thm@headfont{\bfseries}%
  \normalfont 
  \thm@preskip \thm@preskip 
  \thm@postskip\thm@preskip
}
\def\imod#1{\allowbreak\mkern5mu({\operator@font mod}\,\,#1)}
\makeatother

\numberwithin{equation}{subsection}

\title{Geometric Langlands duality for periods}

\author{Tony Feng and Jonathan Wang}

\begin{document}

\begin{abstract}
We study conjectures of Ben-Zvi--Sakellaridis--Venkatesh that categorify the relationship between automorphic periods and $L$-functions in the context of the Geometric Langlands equivalence. We provide evidence for these conjectures in some low-rank examples, by using derived Fourier analysis and the theory of chiral algebras to categorify the Rankin-Selberg unfolding method. 
\end{abstract}

\maketitle

\tableofcontents

\section{Introduction} Many arithmetic questions are encoded in the analytic behavior of $L$-functions, and a powerful tool for understanding $L$-functions is via integral representations as periods of automorphic forms. Some of the first examples of such integral representations include:
\begin{itemize}
\item Riemann's integral representation of the zeta function, as the period of a theta function.
\item Hecke's integral representation of standard $L$-functions for $\GL_2$, as periods of modular forms. 
\item Rankin-Selberg's integral representation of tensor product $L$-functions for $\GL_2$, via Rankin-Selberg convolution of modular forms against mirabolic Eisenstein series. 
\end{itemize}
The purpose of the present paper is to study a \emph{categorification} of automorphic periods and their relation to $L$-functions, in the context of the Geometric Langlands equivalence. In particular, we will focus on the three examples noted above, cast in terms of the \emph{Relative Langlands duality} conjectures of Ben-Zvi--Sakellaridis--Venkatesh \cite{BZSV}, henceforth abbreviated ``BZSV''.

\subsection{Relative Langlands duality} Let us briefly summarize the format of the BZSV conjectures. Langlands defined a duality of split reductive groups, $G \leftrightarrow \chG$. BZSV extend this to a duality  
\[
\{\text{hyperspherical $G$-varieties}\} \longleftrightarrow \{\text{hyperspherical $\chG$-varieties}\}.
\]
We do not define ``hyperspherical'' $G$-variety, but we mention that the cotangent bundles of smooth affine spherical $G$-varieties are among the main examples. There is also a $\Ggr$-action commuting with the $G$-action on the objects of both sides (cf. \cite[\S 1.4]{BZSV}), whose role we suppress in the introduction. 

Now we turn to the global period conjecture of BZSV. Let $C$ be a smooth projective curve over the complex numbers $\CC$. Let $(G, M) \leftrightarrow (\chG, \chM)$ be dual pairs of split reductive groups and hyperspherical varieties over $\CC$. Assume furthermore, partially for simplicity, that we have presentations $M = T^* X$ and $\chM = T^* \chX$ where $X$ (resp. $\chX$) is a spherical variety for $G$ (resp. $\chG)$. Let $\Bun_G$ be the moduli stack of $G$-bundles on $C$, and let $\LocG$ be the (derived) moduli stack of $\chG$-local systems on $C$. 

\subsubsection{Automorphic period sheaf} From the datum $(G,M=T^* X)$, BZSV define a space $\BunX$ equipped with a map $\piaut \co \BunX \rightarrow \Bun_G$,\footnote{This is an oversimplification: the construction of $\BunX$ depends furthermore on the action of $\Ggr$ on $X$.} and an associated (un-normalized) period sheaf 
\[
\cP_X := \piaut_! (\ul{\CC}_{\BunX} := \text{constant sheaf on }\BunX) \in \Dmod(\Bun_G).
\]

\begin{example}[Relation to period integrals]\label{ex: period sheaf} Viewing $\Bun_G$ as the mapping stack $\Map(C, \rB G)$, $\BunX$ is approximately the mapping stack $\Map(C, X/G)$. Suppose $X = G/H$ is a homogeneous spherical variety. Then it is often the case that 
\[
\BunX \cong \Map(C, X/G) \cong \Map(C, \rB H) = \Bun_H.
\]
In the analogous situation of a curve $C$ over a \emph{finite} field $\F_q$, $\Bun_G(\F_q)$ is canonically identified, by Weil's uniformization theorem, with a double coset groupoid of the form 
\[
[G] := G(F) \bs G(\A_F) / \prod_{c \in |C|} G(\wh{\cO}_{C,c})
\]
for $F  = \F_q(C)$ the function field of $C$. In this case $\cP_X$ categorifies the pushforward of the constant function along $[H] \rightarrow [G]$, which represents the distribution 
\[
f \mapsto \int_{[H]} f \, dh
\]
on functions $f$ on $[G]$. In this sense $\cP_X$ categorifies the period integral associated to $X$. 
\end{example}

\subsubsection{Spectral period sheaf}
From the datum $(\chG, \chM=T^* \chX)$, BZSV define a space $\LocX$ equipped with a map $\pispec \co \LocX \rightarrow \LocG$, and an associated (un-normalized) ``$L$-sheaf''
\[
\cL_{\chX} := \pispec_* (\omega_{\LocX} := \text{dualizing complex on }\LocX)^{\shear} \in \IndCoh_{\Nilp}(\LocG).
\]
Here $(-)^{\shear}$ is the ``shearing operation'', which shifts and twists graded components; we will not explain it in the introduction.\footnote{The shearing depends on the action of $\Ggr$ on $\chX$.} 

In contrast to the period sheaf $\cP_X$, which is clearly related to classical automorphic periods as explained in Example \ref{ex: period sheaf}, it is not immediately obvious how $\cL_{\chX}$ is related to an $L$-function. One of the main points of this paper is to explain this precisely. The punchline is that we identify a piece of $\cL_{\chX}$ as the (relative) \emph{factorization homology} of a certain (relative) \emph{cocommutative factorization coalgebra}, which is a sheaf-theoretic analogue of an Euler product. This will be discussed further in \S \ref{sssec: factorization homology}. 

\begin{remark}
It is essential to consider $\LocX$ and $\LocG$ as \emph{derived} stacks. Morally, one should also consider $\BunX$ as a derived stack, but because we will always work with either $D$-modules or constructible sheaves or \'{e}tale sheaves on the automorphic side, it is possible to formulate the period sheaf without doing this. Nevertheless, we will see in certain calculations that it is still useful to incorporate the derived structure.
\end{remark}

\subsubsection{The duality conjecture} There is a conjectural \emph{Geometric Langlands equivalence}
\begin{equation}
\LL_G \co \Dmod(\Bun_G) \xrightarrow{\sim} \IndCoh_{\Nilp}(\LocG),
\end{equation}
which (if it exists) is characterized by various compatibility properties, discussed in \cite{AG15} and \cite{Gai15}. For $\GL_1$, this has been known since the work of Rothstein \cite{Roth96} and Laumon \cite{Lau96} (although in a slightly different formulation), and recently a proof for general groups has been announced by Gaitsgory--Raskin partially in joint work with Arinkin, Beraldo, Campbell, Chen, Faergeman, Lin, and Rozenblyum \cite{GLC}. 

The global period conjecture of BZSV concerns the behavior of the equivalence $\LL_G$ on specific objects. BZSV define ``normalized'' versions of the period sheaf and $L$-sheaf. For the period sheaf, the normalized version $\cP_X^{\norm}$ is obtained from $\cP_X$ by certain shifts and twists. For the $L$-sheaf, $\cL_X^{\norm}$ is obtained from $\cL_{\chX}$ by shifting, twisting, and tensoring with a certain line bundle. Then BZSV conjecture (see \cite[Conjecture 12.1.1]{BZSV} for the precise statement) that the Geometric Langlands equivalence $\LL_G$ takes\footnote{The version we write here omits the duality involution appearing in \emph{loc. cit.}, which appears for reasons of convention in \cite{BZSV} regarding right/left actions. See \cite[\S A.1]{CV} for a discussion of this point.}
\begin{equation}\label{eq: BZSV conj}
\Dmod(\Bun_G) \ni \cP_X^{\norm} \xmapsto{\LL_G} \cL_{\chX}^{\norm} \in \IndCoh_{\Nilp}(\LocG)
\end{equation}

\subsection{Results}We will study the three examples listed at the beginning of the introduction.\footnote{Actually, we replace the Rankin-Selberg period by a singular subperiod.} In the language of relative Langlands duality, they are the {\sf{Tate}}, {\sf{Hecke}}, and {\sf{singular Rankin-Selberg}} (for $\GL_2$) examples in Figure \ref{fig: dual pairs}. 

Strictly speaking, singular spherical varieties fall outside the scope of \cite{BZSV}. However, we expect that there should be some extension of relative Langlands duality that encompasses singular varieties. Part of our motivation for analyzing the {\sf{singular Rankin-Selberg}} example is to investigate what such an extension might look like. The work of Chen-Venkatesh \cite{CV} also investigates examples of singular spherical varieties. 

\begin{figure}[!h]
\centering
\begin{tabular}{|c|c|c|c|c|}
\hline
& $X$ & $G$ & $\chX$ & $\chG$ \\
\hline
\sf{Tate}  & $\A^1$ & $\G_m$ & $\A^1$ & $\G_m$ \\
\sf{Hecke}  & $\GL_2 / \GL_1$ & $\GL_2$ & $\A^2$ & $\GL_2$ \\
$\substack{\text{\sf{singular}}\\ \text{\sf{Rankin-Selberg}}}$ & $\Ind_{\GL_2}^{\GL_2 \times \GL_2} (\A^2 \setminus 0)$ & $\GL_2 \times \GL_2$ & $(\A^2 \otimes \A^2)^{\rank \leq 1}$ & $\GL_2 \times \GL_2$ \\
\hline 
\end{tabular}
\caption{Examples of dual hyperspherical varieties $(G, M = T^* X) \leftrightarrow (\chG, \chM = T^* \chX)$.}\label{fig: dual pairs}
\end{figure}

Our main results take the following form. In each of the examples of Figure \ref{fig: dual pairs}, we produce exact triangles 
\begin{equation}\label{eq: aut triangle}
\cA \rightarrow \cP_X^{\norm} \rightarrow \cB \in \Dmod(\Bun_G)
\end{equation}
and 
\begin{equation}\label{eq: spec triangle}
\check{\cA} \rightarrow \cL_{\chX}^{\norm} \rightarrow \check{\cB} \in \IndCoh_{\Nilp}(\LocG)
\end{equation}
and we prove a result of the form: 

\begin{thm}\label{thm: intro thm}
In each of the examples of Figure \ref{fig: dual pairs}, assuming the existence of the Geometric Langlands equivalence $\LL_G$ for the respective $G$, it matches $\cA \xmapsto{\LL_G} \check{\cA}$ and $\cB \xmapsto{\LL_G} \check{\cB}$ as indicated in the diagram below, where the rows are exact:
\begin{equation}\label{eq: intro thm}
\begin{tikzcd}
\cA \ar[r] \ar[d, "\LL_G", mapsto]  &  \cP_X^{\norm} \ar[r] \ar[d, "\LL_G?", dashed, mapsto] &  \cB \ar[d, "\LL_G", mapsto]  & \in \Dmod(\Bun_G) \ar[d, "\LL_G"] \\
\check{\cA} \ar[r] &  \cL_{\chX}^{\norm} \ar[r] &  \check{\cB} & \in \IndCoh_{\Nilp}(\LocG)
\end{tikzcd}
\end{equation} 
\end{thm}

See Theorem \ref{thm: Tate duality}, Theorem \ref{thm: Hecke}, and Theorem \ref{thm: RS} for the precise statements (which involve some notation that we do not want to explain here). The existence of a decomposition into two pieces is an artifact of the low rank examples that we study; in general we would expect to find multi-step filtrations of $\cP_X^{\norm}$ and $\cL_{\chX}^{\norm}$ whose associated gradeds can be matched under $\LL_G$ by our methods. 

\begin{remark}
Our results are stated in terms of the Geometric Langlands equivalence (which, as already noted, has been announced in \cite{GLC}). Note that it is known in the case of the {\sf{Tate}} period, since there $G = \GL_1$. For the {\sf{Hecke}} and {\sf{singular Rankin-Selberg}} periods, we only use the Geometric Langlands equivalence for $G = \GL_2$, whose proof was outlined several years ago in \cite{Gai15}. 
\end{remark}

Theorem \ref{thm: intro thm} reduces the BZSV Conjecture \eqref{eq: BZSV conj} (in the cases of Figure \ref{fig: dual pairs}) to an extension problem. It would be interesting to investigate the meaning of the extension classes. In the {\sf{Tate}} case, the group of extensions is relatively easy to describe, and there is a distinguished canonical element on each side that we are able to pin down as the relevant extension classes, thus proving the full relative duality conjecture in that case. 

\begin{thm}
In the {\sf{Tate}} case, we have $\cP_X^{\norm} \xmapsto{\LL_G}  \cL_{\chX}^{\norm} $ as conjectured in \eqref{eq: BZSV conj}. 
\end{thm}

Let us comment on the nature of the exact triangles \eqref{eq: aut triangle} and \eqref{eq: spec triangle}. They are induced by geometric decompositions of the hyperspherical varieties. For example, in the Tate case, one has the stratification of $X$ (resp. $\chX$) by $G$-orbits (resp. $\chG$-orbits): $\A^1 = 0 \sqcup (\A^1 \setminus 0)$. On the automorphic side, this induces an open-closed decomposition of $\BunX$, and \eqref{eq: aut triangle} comes from the excision triangle for this open-closed decomposition.

On the spectral side, the stratification by $\chG$-orbits induces (in all three examples) an open-closed decomposition 
\[
\begin{tikzcd}
V \ar[r, hook, "\text{closed}"] &  \LocX & W \ar[l, hook', "\text{open}"']
\end{tikzcd}
\]
and \eqref{eq: spec triangle} comes from the corresponding exact triangle for cohomology with supports, 
\[
\rR\Gamma_V(\LocX, -) \rightarrow \rR\Gamma(\LocX, -) \rightarrow \rR\Gamma(W, -).
\]
In particular, the term $\check{\cA}$ in \eqref{eq: spec triangle} is the relative local cohomology along the zero-section $\LocG \inj \LocX$ of the dualizing sheaf of $\LocX$. 


In the case of the {\sf{Tate}} period, there is a matching between the $G$-orbits on $X$ and $\chG$-orbits on $\chX$ (reversing the closure relations) such that matching orbits contribute corresponding pieces under the Geometric Langlands equivalence. However, this naive matching does not persist in general: in both the {\sf{Hecke}} and {\sf{singular Rankin-Selberg}} periods, $X$ is homogeneous for the $G$-action, while there are two $\chG$-orbits on $\chX$. In those cases, the exact triangle for $\cP_X^{\norm}$ comes instead from a stratification \emph{in the Fourier dual space}. A space and its Fourier dual are united in microlocal geometry, suggesting that in general, one should look for a decomposition in the hyperspherical variety $M  = T^* X$. Indeed, a closed-open stratification $X = Z \sqcup U$ induces a decomposition of $T^*X$ as
\begin{equation}\label{eq: Lagrangian decomposition}
T^* X =   (Z \times_X T^* X) \sqcup (T^* U) .
\end{equation}
Note that $Z \times_X T^*X$ can be thought of as the embedding in $T^*X$ of the Lagrangian correspondence 
\[
\begin{tikzcd}
T^*Z  & \ar[l, hook'] Z \times_X T^*X \ar[r, twoheadrightarrow] & T^*X
\end{tikzcd}
\]
associated to the map $Z \inj X$. More generally, it happens that $T^*X$ admits a decomposition into Lagrangian correspondences, which may not come from a stratification of $X$. We then see a corresponding decomposition of the (automorphic or spectral) $X$-period, as in \eqref{eq: aut triangle} or \eqref{eq: spec triangle}. See \S \ref{ssec: microlocal} for how this plays out in the example of the {\sf{Hecke}} period. 

This suggests that the duality $M \leftrightarrow \chM$ often preserves more structure on both sides, which would be interesting to codify.

\subsection{Methods} The proof of Theorem \ref{thm: intro thm} for the {\sf{Tate}} case is in some sense by direct computation. For the {\sf{Hecke}} and {\sf{singular Rankin-Selberg}} periods, the classical approach is via the ``unfolding method'' of Rankin-Selberg, and our approach can be seen as a categorification of the Rankin-Selberg method. However, we note some significant differences which make our task more difficult. The Rankin-Selberg method analyzes the period of a single cuspidal Hecke eigenform, while we want to analyze the automorphic period sheaf as a distribution on all of $\Bun_G$. If we were only interested in the pairing of the period sheaf against individual cuspidal Hecke eigensheaves, then the relevant analysis has already been carried out by Lysenko \cite{Lys02}. But our task also includes the issue of analyzing periods of Eisenstein series where the corresponding $L$-function has a pole, a subtle topic which is not understood even in the analysis of classical automorphic forms. Miraculously, we find a very clean answer; in some sense, working at the level of sheaves allows us to bypass convergence problems that would make it difficult to articulate the answer numerically. 

We highlight here some ways in which our calculations make interesting contact with other themes in mathematics. 

\subsubsection{Derived Fourier analysis} The premise of the unfolding method is to ``unfold'' automorphic periods in terms of Whittaker periods. Classically, this means that one should consider the Fourier expansion of a Hecke eigenform. 

But in the geometric context, it turns out that the spaces which arise naturally from trying to geometrize natural vector spaces are not vector bundles, which a priori precludes us from applying the sheaf-theoretic Fourier transform. For example, in the {\sf{Tate}} case one is led to consider the space $\BunX$ over $\Bun_{\GL_1}$ whose fiber over a line bundle $L$ its space of global sections $\rH^0(C; L)$. Thus each fiber is a vector space, but the fibers vary in rank (even on connected components) hence cannot interpolate into a vector bundle. For this reason it is natural to consider \emph{derived vector bundles}, which are derived linear spaces associated to perfect complexes, generalizing how vector bundles are associated to locally free coherent sheaves. For example, one wants to enhance $\BunX$ to a derived vector bundle over $\Bun_{\GL_1}$ whose fibers are, informally speaking, the complexes $\rR\Gamma(C;L)$. In general, derived vector bundles have both ``derived'' and ``stacky'' aspects, which are interchanged under duality. 

In order to execute unfolding at the geometric level, it is therefore necessary to use an expansion of the Fourier transform to derived vector bundles, which we call the \emph{derived Fourier transform}. It turns out that the good properties of the Deligne-Laumon Fourier transform extend from vector bundles to derived vector bundles, although this is technically quite tricky to establish -- see \cite[\S 6 and Appendix A]{FYZ3}. The trickiest aspect is to show the near-involutivity of the Fourier transform, for which the crux is the statement that the derived Fourier transform of the constant sheaf is (an appropriate shift and twist of) the delta-sheaf at the origin of the dual bundle. It is this statement which we use to unfold the {\sf{Hecke}} and {\sf{singular Rankin-Selberg}} period sheaves. 

We arrived at the formulation of derived Fourier analysis when contemplating the \emph{functional equation} for the {\sf{Tate}} period (\S \ref{ssec: functional equation}). The theory was written up in \cite{FYZ3}, which used it to study the modularity of higher theta functions in the sense of \cite{FYZ2}. Not coincidentally, the modularity of theta functions is closely connected to the functional equation for zeta functions.

\subsubsection{Factorization homology}\label{sssec: factorization homology} Our analysis gives a new perspective on the appearance of $L$-functions in periods of automorphic forms. Let $V$ be a representation of $\chG$. Then to a Galois representation $\sigma \co \Gal(F^s/F) \rightarrow \chG(\ol{\Q}_{\ell})$ with Frobenii having rational characteristic polynomials, one forms the $L$-function as an Euler product 
\[
L(V,\sigma, s) := \prod_{v \in |F|} L_v(V, \sigma, s)
\]
where $L_v(V,\sigma, s)$ is some local factor at the place $v \in |F|$ formed using the characteristic polynomial of $\sigma(\Frob_v)$ acting on $V$. 

The geometric counterpart to such Euler products is \emph{factorization homology}. To explain it, let $Y$ be a variety over a field $\CC$. (We will eventually apply this theory to $Y = C$, although there is a version of the theory for any variety over any field.) We write $\pt := \Spec \CC$. Then we have a pair of adjoint functors
\[
\begin{tikzcd}
\Dmod(Y) \ar[r, "\pr_!", bend left]  & \Dmod(\pt)  \ar[l, "\pr^!", bend left] 
\end{tikzcd}
\]
where the derived categories of $D$-modules $\Dmod(Y)$ and $\Dmod(\pt)$ can be equipped with a symmetric monoidal structure with respect to the $!$-tensor product (under the usual identification  of $\Dmod(\pt)$ with the derived category of $\CC$-vector spaces, this is the usual tensor product). The functor $\pr^!$ is symmetric monoidal with respect to these symmetric monoidal structures. Therefore, it induces a functor
\begin{equation}\label{eq: !-pull CAlg}
\CAlg(\Dmod(Y) ) \xleftarrow{\pr^!} \CAlg(\Dmod(\pt)),
\end{equation}
where $\CAlg(-)$ denotes the category of commutative algebras in a symmetric monoidal category. Then the \emph{factorization homology} functor 
\begin{equation}\label{eq: factorization homology}
\rR \Gamma_c^{\Fact}(Y;-) \co \CAlg(\Dmod(Y) ) \rightarrow \CAlg(\Dmod(\pt) )
\end{equation}
is the \emph{left adjoint} to \eqref{eq: !-pull CAlg}. 

\begin{example}
If $Y = \{ y_1, \ldots, y_n\}$ is a finite disjoint union of points, then there is an evident equivalence of categories 
\[
\CAlg(\Dmod(Y )) \cong \bigotimes_{i=1}^n \CAlg(\Dmod(\pt)),
\]
under which identification we have 	
\[
\rR \Gamma_c^{\Fact}( Y;  \cF) \cong \bigotimes_{i=1}^n \cF_{y_i}.
\]
This suggests the intuition, emphasized in the work of Gaitsgory-Lurie \cite{GL16}, that in general $\rR\Gamma^{\Fact}_c(Y; \cF)$ is a ``continuous tensor product of $\cF$ over the points of $Y$''. The analogy to Euler products is evident. 
\end{example}

In the work of Gaitsgory-Lurie \cite{GL16} on Weil's Tamagawa Number Conjecture over function fields, a certain product of zeta functions arose as the trace of Frobenius on a certain factorization homology group. In our analysis here, general $L$-functions (and even more general objects) will be categorified by factorization homology. More precisely, let $\omega_V$ be the dualizing complex of the $\chG$-representation $V$ regarded as an affine space over $\CC$, and consider the local cohomology at the origin $\rR\Gamma_0(V; \omega_V)$; informally, it is the space of distributions on $V$ with set-theoretic support at $0$. Then $\rR\Gamma_0(V; \omega_V)$ is equipped with the natural structure of a cocommutative coalgebra in $\Rep(\chG)$, which allows to form a corresponding $\chG$-equivariant \emph{cocommutative factorization coalgebra} on $C$. The $\chG$-equivariant structure allows to twist this cocommutative factorization coalgebra by a $\wh{G}$-local system $\sigma$ on $C$, and we regard its factorization homology as the categorification of the $L$-function $L(V,\sigma, s)$. 

As for how such objects arise in our calculation, suppose $\chX = V$ is a (linear) $\chG$-representation, as is the case in the {\sf{Tate}} and {\sf{Hecke}} periods. Then the origin $0 \in V$ contributes a closed $\chG$-stratum, and we will stratify $\LocX$ by the induced zero-section of $\LocG$. Then, using the work of Beilinson-Drinfeld on \emph{chiral algebras}, we prove a local-global principle that relates the derived pushforward of $\omega_{\LocX}$, with set-theoretic support along the zero-section, to the relative factorization homology of the cocommutative factorization coalgebra attached to $\rR\Gamma_0(V; \omega_V)$ (and twisted by the universal $\wh{G}$-local system on $C$). More generally, the same argument will show that each $\chG$-fixed point $x \in \chX$ contributes a term to the $L$-sheaf $\cL_{\chX}$ which is the relative factorization homology over $\LocG$ of the local cohomology $\rR\Gamma_x(\chX; \omega_{\chX})$ turned into a cocommutative factorization coalgebra and then twisted by the universal local system. However, if $x$ is not a \emph{smooth} point of $\chX$, as happens in the {\sf{singular Rankin-Selberg}} period, then the resulting factorization homology is analogous to an Euler product which is \emph{not} an $L$-function. 

\begin{remark}
The paper \cite{BZSV} proposes some other ways in which arithmetic local-global principles, such as Euler products, should be categorified by factorization algebras and factorization homology. More precisely, \cite[\S 16]{BZSV} discusses the how the \emph{RTF-algebra} of a spherical variety -- which approximately categorifies the associated relative trace formula -- arises as the factorization homology of the (local) \emph{Plancherel algebra}. Under a suitable factorizable form of the local conjectures in \cite[\S 7,8]{BZSV}, the RTF-algebra would be identified (as an object of the \emph{global Hecke category}, obtained as factorization homology of the factorizable derived Satake category) with the algebra of ``$L$-observables'' obtained from a dual hyperspherical variety on the spectral side, which is examined in \cite[\S 17]{BZSV}. This algebra of $L$-observables is related to the ``$L^2$ version'' \cite[Conjecture 12.8.1]{BZSV} of the global duality conjecture. That version treats more general situations than the one from \cite[Conjecture 12.1.1]{BZSV} formulated in \eqref{eq: BZSV conj}, in that it does not require the dual hyperspherical variety $\chM$ to be polarized, but it is also less precise because it does not pin down a Langlands dual description of the period sheaf $\cP_X$; roughly speaking, it describes the cuspidal part of the ``square'' of the period sheaf.
\end{remark}

\subsection{Acknowledgments} We thank David Ben-Zvi, Lin Chen, Gurbir Dhillon, Vladimir Drinfeld, John Francis, Charles Fu, Quoc Ho, Sam Raskin and Yiannis Sakellaridis for relevant conversations. We especially thank Dennis Gaitsgory and Akshay Venkatesh for many discussions and explanations regarding the material here. We are grateful to David Ben-Zvi, Sanath Devalapurkar, Sergey Lysenko, and two very attentive referees for comments on a draft. TF was supported by the NSF grant DMS-2302520 and JW was supported by NSF grant DMS-1803173.

\section{Notation and conventions}

\subsection{Coefficient fields}
Let $\F$ be an algebraically closed field, which will be the ground field for the ``automorphic side''. 

Let $k$ be an algebraically closed field of characteristic $0$, which will be the ground field for the ``spectral side'', as well as the coefficients for sheaves on the automorphic objects. When discussing the ``de Rham'' setting, we will have a fixed identification $k \cong \CC$. 

\subsection{Reductive groups} 
Let $G$ be a split reductive group over $\F$ and $\chG$ its Langlands dual group, regarded as a split reductive group over $k$.

\subsection{Categories and sheaves} When our categories are linear over a field, that field will be assumed to be $k$, except in the discussion of \S \ref{ssec: non-unital} and \S \ref{ssec: fact homology}. All our functors are derived, e.g., $f_*$ means $\rR f_*$, $\Gamma$ means $\rR \Gamma$, etc. If we need to refer to individual cohomology groups, we will write $\rR^i f_*$ and $\rH^i$, etc. 

Our categories are derived, e.g., $\Vect$ is the derived category of $k$-vector spaces and $\Rep(\chG)$ is the derived category of $\chG$-representations, unless specified otherwise. 


We denote by $\ulk_S$ the constant sheaf on a space $S$ with value $k$. When context is clear, we may omit the subscript. 

When $\F = \CC$, we denote by $\Dmod(S)$ the derived category of $D$-modules on $S$. 

\subsubsection{Categories of sheaves}\label{sssec: 3 sheaves} Recall that there are (at least) three flavors of the Geometric Langlands equivalence, which we refer to as \emph{de Rham} \cite{AG15}, \emph{Betti} \cite{BZN18}, and \emph{finite field} \cite{sixI}. Our final theorems are all about the de Rham setting because this is where the proof of the Geometric Langlands Equivalence is best documented, but many of our calculations are agnostic to the specific sheaf theory, so we introduce a notation that is agnostic to the sheaf theory. For an algebraic stack $S$ over $\F$, we denote by $\Shv(S)$ any of the three flavors of ``topological'' sheaf theories \cite[\S B.4--B.7]{BZSV}:
\begin{enumerate}
\item (de Rham) If $\F = \CC$, then we may take $\Shv(S)$ to be the derived category of $D$-modules on $S$, denoted $\Dmod(S)$. In this case we must work with $k = \CC$.
\item (Betti) If $\F = \CC$, then we may take $\Shv(S)$ to be the derived category of all sheaves on $S(\CC)$  equipped with the analytic topology, in the sense of \cite[Appendix G]{sixI}. There are no restrictions on $k$. 
\item (finite field) If $\F$ has positive characteristic, then we may take $\Shv(S)$ to be the derived category ind-constructible $\ell$-adic \'{e}tale sheaves on $S$. Then $k =\ol{\Q}_{\ell}$ for some prime $\ell$. 
\end{enumerate}
We denote by $\Shv_{\hol}(S) \subset \Shv(S)$ the subcategory of \emph{holonomic} sheaves. In the $D$-module context this means holonomic D-modules in the usual sense, while in the Betti and $\ell$-adic contexts it means constructible objects. 

\subsubsection{Artin-Schreier}\label{sssec: AS} The \emph{Artin-Schreier sheaf} refers to a certain object \cite[Definition 10.5.1]{BZSV} in $\Shv(\A^1)$, defined case-by-case in each flavor of sheaf theory: 
\begin{enumerate}
\item (de Rham) The exponential $D$-module. 
\item (Betti) The locally constant $\CC^\times$-equivariant (via squaring) sheaf on $\CC$ defined as $(j_! k^- \oplus j_* k)[-1]$ where $j \co \CC^\times \inj \CC$ and $k,k^-$ are the trivial and non-trivial rank one local systems on $\CC^\times/\CC^\times \cong \rB \mu_2$. 
\item (Finite Field) The \'{e}tale Artin-Schreier sheaf induced by a choice of non-trivial additive character $\psi \co \F_q \rightarrow k^\times$. 
\end{enumerate}

\subsubsection{Shifting and twisting}

For $n \in \Z$, the endofunctor $\cK \mapsto \cK(n)$ of $\Shv(S)$ is the $n$th Tate twist. This only has meaning in the finite field context, where we will choose a square root of the cyclotomic character to make sense of half Tate-twists. We use the notation $\tw{n} := \Pi^n [n](n/2)$ as in \cite[(2.4)]{BZSV}. Note that  ``$\tw{n}$'' in \cite{FYZ3} is instead called $\tw{2n}$ in this paper in order to be consistent with \cite{BZSV}.

\subsubsection{Correspondences}
Consider a correspondence
\[
\begin{tikzcd}
A_1 & A^\flat \ar[l, "c_1"'] \ar[r,"c_2"] & A_2
\end{tikzcd}
\]
Our convention is that ``the functor induced by $\cK \in \Shv(A^\flat)$'' is 
\[
c_{2!} (\cK \otimes c_1^* (-)) \co \Shv(A_1) \dashrightarrow \Shv(A_2).
\]
Note that in the de Rham case this is not a priori defined on all of $\Shv(A_1)$, which is why we use the dashed arrow; we will only apply it to objects on which it is defined.

\subsubsection{Graded categories}
For a monoid $\Lambda$ and a category $\msf{C}$, we denote by $\msf{C}^{\Lambda}$ the category of $\Lambda$-graded objects in $\msf{C}$.

\subsection{Geometry}
Let $C$ be a smooth projective curve over $\F$. For $d \geq 0$, we write $C^{(d)} = \Sym^d C$ for the $d$th symmetric power of $C$, with the convention $C^{(0)} := \Spec \F$.

\subsubsection{Spin structure}\label{ssec: spin} Fix a square root $\Omega_C^{1/2}$ of the canonical line bundle $\Omega_C$. The Conjectures of Ben-Zvi--Sakellaridis--Venkatesh are formulated in terms of such a choice -- see \cite[\S 10.1.2]{BZSV}.

\subsubsection{Derived mapping stacks}\label{sssec: derived mapping stacks} If $X$ is a proper scheme and $Y$ is a stack locally of finite presentation over a field, then the derived mapping stack $\Map(X,Y)$ is constructed in \cite[\S 3.6]{TV05} and \cite[\S 2.2.6.3]{TV08} (the first reference constructs the internal hom on stacks on any site, and the second establishes its geometricity properties under the given assumptions). On objects, $\Map(X,Y)$ sends an animated ring $R$ to the anima of morphisms $X_R \rightarrow Y$. 

Suppose $Y$ is equipped with a map $\pi \co Y \rightarrow X$. Then we have the derived space of sections $\hSect(X,Y)$, defined as the derived fibered product 
\[
\begin{tikzcd}
\hSect(X,Y) \ar[r] \ar[d] & \Map(X,Y) \ar[d, "\pi"]  \\
\{\Id\} \ar[r] 
& \Map(X,X)
\end{tikzcd}
\]

\subsubsection{Moduli of $G$-bundles} By definition, for $\rB G$ the classifying stack of $G$, we have 
\[
\Bun_G := \Map(C, \rB G).
\]
Although defined a priori as a derived stack, this is in fact a classical Artin stack \cite[Example 5.3]{FYZ2}. 

\subsubsection{Moduli of $\chG$-local systems} Parallel to the sheaf theory, there are several different flavors of moduli spaces of local system -- see \cite[\S 11.1.2]{BZSV} for discussion of them. Our calculations will be somewhat specific to the de Rham setting, so we recall that only. 

Recall the \emph{de Rham stack} $C_{\dR}$ associated to $C$ from \cite[\S 3.1]{Gai15}. For $\rB \chG$ the classifying stack of $\chG$, we define 
\[
\LocG := \Map(C_{\dR}, \rB \chG).
\]
We abbreviate $\Loc_n := \Loc_{\GL_n}$ and $\Bun_n := \Bun_{\GL_n}$. 
	
\subsubsection{Tangent and cotangent complexes} 
Let $f \co X \rightarrow Y$ be a map of derived Artin stacks. We denote by $\bL_f$ the cotangent complex of $f$, and by $\bT_f$ the tangent complex of $f$ if $\bL_f$ is a perfect complex (in which case $\bT_f$ is defined to be the dual of $\bL_f$). 

\subsubsection{Classical truncation} For a derived stack $S$, we denote by $S_{\cl}$ its classical truncation, which comes equipped with a canonical closed embedding $S_{\cl} \inj S$. 

\subsection{Coherent sheaves} We follow \cite{Gai13} in our conventions on (ind-)coherent sheaves. For a map of derived Artin stacks $f \co X \rightarrow Y$, we denote by $\omega_{X/Y} := f^! \cO_Y$ the relative dualizing sheaf. For $Y$ a point, we abbreviate $\omega_X := \omega_{X/Y}$. 

\subsubsection{Perfect complexes} For a derived stack $S$, we denote by $\Perf(S)$ the category of perfect complexes on $S$. This is equivalent to complexes of finite tor-dimension, defined in \cite[Tag 0652]{stacks-project}, which we follow in using \emph{cohomological grading} to define tor-amplitude. 

For $\cE \in \Perf(S)$, we denote by $\cE^\vee := \cRHom(\cE, \cO_S)$ the dual perfect complex.

\subsubsection{Coherent singular support}\label{ssec: singular support}
On the spectral side, the period sheaves $\pispec_* (\omega_{\LocX})$ naturally live in $\IndCoh(\LocG)$ but typically do not lie in the full subcategory $\IndCoh_{\Nilp}(\LocG)$ which is the domain of the Geometric Langlands equivalence. The embedding of this full subcategory is left adjoint to a co-localization functor $\IndCoh(\LocG) \rightarrow \IndCoh_{\Nilp}(\LocG)$. For an object $\cF \in \IndCoh(\LocG)$, we write $\cF \in \IndCoh_{\Nilp}(\LocG)$ for its image under this co-localization. That is, we use the specification of the ambient category in order to indicate the application of a singular support co-localization functor.  

\part{Tools}

\section{Derived Fourier analysis}\label{sec: FT}

\subsection{Derived Fourier transform}\label{ssec: Derived FT}

In this subsection we recall the ``Derived Fourier analysis'' developed in \cite{FYZ3}. While \emph{loc. cit.} worked in the setting of $\ell$-adic \'etale shaves, essentially the same constructions and arguments go through more generally for holonomic sheaves in the sense of \S \ref{sssec: 3 sheaves}. We shall give the statements, omitting proofs since these are identical to the $\ell$-adic setting. 

\subsubsection{Derived vector bundles}\label{sssec: derived vector bundles}
Let $S$ be a derived Artin stack. There is a functor $\Tot_S$ from the category $\Perf(S)$ of perfect complexes on $S$ to the category of derived stacks over $S$, which extends the usual construction of a vector bundle from a locally free coherent sheaf. We normalize this construction as in \cite[\S 6]{FYZ3}. For $\cE \in \Perf(S)$, we will call $E := \Tot_S(\cE)$ the associated \emph{derived vector bundle} associated to $E$. The \emph{virtual rank} of $E$, denoted $\rank(E)$, is the locally constant function on $S$ given by $s\mapsto \chi(\cE_{s})$, the Euler characteristic of the fiber of $\cE$ at a geometric point $s$.

\begin{defn}
For a derived vector bundle $E \rightarrow S$, we denote by $z_E \co S \rightarrow E$ the zero-section. Note that $z_E$ need not be a closed embedding. We write $\delta_E := z_{E!} (\ulk_S ) \in \Shv(E)$ and call it the \emph{delta-sheaf of $E$}.  
\end{defn}

\subsubsection{The derived Fourier transform for derived vector bundles}

Let $S$ be a derived Artin stack, $\cE \in \Perf(S)$. For $\cE^\vee \in \Perf(S)$ the linear dual of $\cE$, we have a tautological evaluation pairing $\cE \otimes \cE^\vee \rightarrow \cO_S$. Setting $E := \Tot_S(\cE)$ and $E^\vee := \Tot_S(\cE^\vee)$, this induces on total spaces a map 
\[
\ev \co E \times_S E^\vee \rightarrow \A^1.
\]
Let $\AS \in \Shv(\A^1)$ be the Artin-Schreier sheaf (cf. \S \ref{sssec: AS}), which we note is holonomic. 

The \emph{derived Fourier transform}
\[
\FT_E \co \Shv_{\hol}(E) \rightarrow \Shv_{\hol}(E^\vee)
\]
is the functor 
\[
\cK \mapsto \pr_{2!} (\pr_1^* (\cK) \otimes \ev^* \AS)[\rank(E)]
\]
where the maps are as in the diagram 
\[
\begin{tikzcd}
& E \times_S E^\vee \ar[dl, "\pr_1"'] \ar[dr, "\pr_2"] \ar[r, "\ev"] & \A^1 \\
E & & E^\vee
\end{tikzcd}
\]

We now tabulate some basic properties of the derived Fourier transform. Below we let $r$ be the virtual rank of $E \rightarrow S$.

\subsubsection{Base change}\label{sssec: FT base change} Let $h \co \wt{S} \rightarrow S$ be a map of derived stacks. For a derived vector bundle $E \rightarrow S$, let $\wt{E} \rightarrow E$ be its base change along $h$. So we have derived Cartesian squares 
\[
\begin{tikzcd}
\wt{E} \ar[r, "h^E"] \ar[d] & E \ar[d] \\
\wt{S} \ar[r, "h"] & S
\end{tikzcd} 
\hspace{1cm}
\begin{tikzcd}
\wt{E}^\vee \ar[r, "h^{E^\vee}"] \ar[d] & E^\vee \ar[d] \\
\wt{S} \ar[r, "h"] & S
\end{tikzcd}
\]

Then there are canonical natural isomorphisms of functors $\Shv_{\hol}(E) \rightarrow \Shv_{\hol}(\wt{E}^\vee)$
\begin{equation}\label{eq: FT bc h^*}
\FT_{\wt{E}} \circ (h^E)^* \cong (h^{E^\vee} )^* \circ \FT_E
\end{equation}
\begin{equation}\label{eq: FT bc h^!} 
\FT_{\wt{E}} \circ (h^E)^! \cong (h^{E^\vee} )^! \circ \FT_E
\end{equation}
and canonical natural isomorphisms of functors $\Shv_{\hol}(\wt{E}) \rightarrow \Shv_{\hol}(E^\vee)$
\begin{equation}\label{eq: FT bc h_!}
\FT_{E} \circ (h^E)_!  \cong (h^{E^\vee} )_! \circ \FT_{\wt{E}}
\end{equation}
\begin{equation}\label{eq: FT bc h_*}
\FT_{E} \circ (h^E)_*  \cong (h^{E^\vee} )_* \circ \FT_{\wt{E}}.
\end{equation}

\subsubsection{Involutivity} In the ``de Rham'' and ``Finite Field'' cases, there is a canonical natural isomorphism $\FT_{E^\vee} \circ \FT_{E} \cong \mrm{mult}_{-1}^* (-r)$ of functors $\Shv_{\hol}(E) \rightarrow \Shv_{\hol}(E)$, where $\mrm{mult}_{-1}$ is multiplication by $-1$ on $E$. 

In the ``Betti'' case, there is also such a natural isomorphism for the full subcategory of \emph{$\G_m$-equivariant sheaves}, which is all that is needed for the purposes of this paper. In fact, the theory of derived Fourier analysis for $\G_m$-equivariant sheaves is developed in a uniform way (for general sheaf-theoretic contexts) in \cite[\S 8]{FK24}. 

\begin{remark} 
The construction of this natural isomorphism is the central point of derived Fourier analysis. It is significantly more involved than in the situation of classical vector bundles. 
\end{remark}

\subsubsection{Functoriality}\label{sssec: FT functoriality} Let $f \co E' \rightarrow E$ be a linear map of derived vector bundles having virtual ranks $r', r$ respectively. This induces a morphism $f^\vee \co E^\vee \rightarrow E'^\vee$ of dual derived bundles. Then we have canonical natural isomorphisms of functors $\Shv_{\hol}(E') \rightarrow \Shv_{\hol}(E^\vee)$: 
\begin{enumerate}
\item\label{item: functoriality 1} $f^{\vee *} \circ \FT_{E'} \cong \FT_E \circ f_! [r'-r]$,
\item\label{item: functoriality 1.25} $f^{\vee !} \circ \FT_{E'} \cong  \FT_E \circ f_*[r-r'](r-r')$,
\end{enumerate}
and canonical natural isomorphisms of functors $\Shv_{\hol}(E) \rightarrow \Shv_{\hol}(E'^\vee)$: 
\begin{enumerate}[resume]
\item\label{item: functoriality 1.5} $\FT_{E'} \circ f^* \cong f^\vee_! \circ \FT_E [r-r'] (r-r')$,
\item\label{item: functoriality 2} $\FT_{E'} \circ f^! \cong f^\vee_* \circ \FT_E[r'-r]$. 
\end{enumerate}

\begin{example}\label{ex: FT delta}
Natural isomorphism \eqref{item: functoriality 1} gives an isomorphism
\[
\FT_E(\delta_E) \cong \ulk_{E^\vee} [r] \in \Shv(E^\vee)
\]
and natural isomorphism \eqref{item: functoriality 1.5} gives an isomorphism
\[
\FT_{E^\vee} (\ulk_{E^\vee}) \cong  \delta_E [-r](-r) \in \Shv(E).
\]
\end{example}

\section{Chiral algebras and factorization homology}\label{sec: factorization algebras}
In this section, we recall or develop some tools in the theory of \emph{chiral algebras}. This subject was developed by Beilinson-Drinfeld \cite{BD04}, although our presentation instead follows the references \cite{FG12, GL16, Ho17, Ho21a, Ho21b}. Most of this section just 	reviews material from the literature: factorization (co)algebra theory is reviewed in \S \ref{ssec: fact homology} and \S \ref{ssec: graded factorization algebras}, and derived infinitesimal geometry is reviewed in \S \ref{ssec: derived infinitesimal geometry}. The only ``new'' material is \S \ref{ssec: factorization of distribution coalgebras}, which studies factorization homology of cocommutative coalgebras of distributions. Our work suggests that these are the correct categorification of $L$-functions. 

\subsection{Non-unital algebras}\label{ssec: non-unital} Let $(\msf{C}, \otimes)$ be a stable symmetric monoidal category. There is an equivalence between the notion of augmented unital commutative algebras in $(\msf{C}, \otimes)$, by which we mean algebras $A$ equipped with an augmentation $\epsilon \co A \rightarrow \mathbf{1}$ (morphisms are required to be compatible with the augmentation) and \emph{non-unital} algebras $\ol{A}$. The equivalence sends an augmentated $k$-algebra $A$ to $\ol{A} := \ker(\epsilon)$, and in the other direction sends a non-unital algebra $\ol{A}$ to $\mathbf{1} \oplus \ol{A}$.

In this section, we denote by $\CAlg(\msf{C}, \otimes)$ the category of \emph{non-unital} commutative algebras in the symmetric monoidal category $(\msf{C}, \otimes)$. 

A similar discussion applies to augmented cocommutative coalgebras. We denote by $\CoCAlg(\msf{C}, \otimes)$ the category of \emph{non-unital} cocommutative coalgebras in the symmetric monoidal category $(\msf{C}, \otimes)$.

\subsection{Factorization algebras}\label{ssec: fact homology}


\subsubsection{The Ran space}
Suppose $C$ is a proper variety over $\F$. Let $\FS$ be the category of finite sets, with morphisms being \emph{surjective} maps. The \emph{Ran space} of $C$ is the prestack 
\[
\Ran(C)  := \colim_{I \in \FS^{\op}} C^I.
\]

The category of sheaves $\Shv(\Ran C)$ is defined by descent along $!$-pullback. We have a ``union'' map
\[
\Ran(C) \times \Ran(C) \xrightarrow{\un} \Ran(C).
\]
Since the transition maps defining $\Ran C$ are proper, $\un$ has a notion of $!$-direct image.

\subsubsection{The convolution tensor structure} The map 
\[
\un_! \co  \Shv(\Ran(C) \times \Ran(C)) \rightarrow \Shv(\Ran(C))
\]
induces a symmetric monoidal structure on $\Shv(\Ran(C))$, which we denote $\otstar$. We call this the \emph{convolution tensor structure}.

\begin{remark}
By its construction, for the map $\pr_{\Ran} \co \Ran C \rightarrow \pt := \Spec \F$ the proper pushforward $\pr_{\Ran !}$ \emph{is symmetric monoidal}. Consider the commutative diagram 
\[
\begin{tikzcd}
\CAlg(\Shv(C), \sotimes) & \CAlg(\Shv(\Ran C), \otstar) \ar[l, "\Delta^!"']  \\
& \CAlg(\Shv(\pt)) \ar[ul, "\pr^!"] \ar[u, "\pr^!_{\Ran}"'] 
\end{tikzcd}
\]
It turns out that the left adjoint to $\Delta^!$ exists (and will be discussed below): it is the formation of the associated \emph{commutative factorization algebra}. Therefore, forming the left adjoint diagram shows that the left adjoint to $\pr^!$, which was called \emph{factorization homology} in \S \ref{sssec: factorization homology}, agrees with $\pr_{\Ran !}$ composed with the the left adjoint to $\Delta^!$. 
\end{remark}

\subsubsection{Commutative factorization algebras} Recall that $\CAlg(\Shv(\Ran C), \otstar)$ denotes the category of commutative (non-unital) algebra objects in the convolution tensor structure. In particular, an object of $\CAlg(\Shv(\Ran C), \otstar)$ entails $\cA \in \Shv(\Ran C)$ plus a map 
\[
\un_! (\cA \boxtimes \cA) \rightarrow \cA.
\]
By adjunction, such data is equivalent to $\cA \boxtimes \cA \rightarrow \un^! \cA$. We often abuse notation by referring to an object of $\CAlg(\Shv(\Ran C) ,\otstar)$ by ``$\cA$'', suppressing the multiplication map and higher coherence data.

The disjoint subspace $(\Ran C \times \Ran C)_{\disj}$ is defined as 
\[
(\Ran C \times \Ran C)_{\disj}  = \colim_{\substack{I_1 \in \FinSurj^{\op} \\ I_2 \in \FinSurj^{\op}}} (C^{I_1} \times C^{I_2})_{\disj}
\]
where $(C^{I_1} \times C^{I_2})_{\disj}$ parametrizes pairs $(D_1 \in C^{I_1}, D_2 \in C^{I_2})$ such that $D_1$ and $D_2$ have disjoint support. 

\begin{defn}\label{def: factorization algebra}
A \emph{factorization algebra} is an $\cA \in \CAlg(\Shv(\Ran C), \otstar)$ such that the structure map 
\[
\cA \boxtimes \cA \rightarrow \un^! \cA
\]
restricts to an isomorphism on the open subspace $(\Ran C \times \Ran C)_{\disj}$. The full subcategory of factorization algebras is denoted $\CAlg_{\Fact}^\star(\Shv(\Ran C)) \subset \CAlg(\Shv(\Ran C), \otstar) $.
\end{defn}

\subsubsection{Diagonal embedding}

The diagonal map $\Delta \co C \rightarrow \Ran C$ induces a symmetric monoidal functor 
\[
\Delta^! \co (\Shv(\Ran C), \otstar) \rightarrow (\Shv(C), \sotimes)
\]
because of the Cartesian square 
\begin{equation}\label{eq: diag Ran}
\begin{tikzcd}
C \ar[r, "\Delta \times \Delta"] \ar[d, " \Id"] & \Ran C \times \Ran C \ar[d, "\un"] \\
C \ar[r, "\Delta"] & \Ran C
\end{tikzcd}
\end{equation}

\begin{thm}[{\cite[Theorem 5.6.1]{GL16}}]\label{thm: commutative fact alg} The functor $\Delta^!$ admits a fully faithful left adjoint 
\[
\Fact \co \CAlg(\Shv(C), \sotimes) \rightarrow \CAlg (\Shv(\Ran C), \otstar),
\]
whose essential image is $\CAlg_{\Fact}^{\star}(\Shv(\Ran C))$. In particular, $\Delta^!$ restricts to an equivalence 
\[
\Delta^! \co \CAlg_{\Fact}^\star(\Shv(\Ran C)) \xrightarrow{\sim} \CAlg(\Shv(C), \sotimes).
\]
\end{thm}

\begin{example}[Free factorization algebras]\label{ex: free fact algebra}
The symmetric monoidality of $\Delta^!$ implies that we have a commutative diagram 
\[
\begin{tikzcd}
\Shv(C) &  \CAlg(\Shv(C), \sotimes) \ar[l, "\oblv"] \\
\Shv(\Ran C) \ar[u, "\Delta^!"] & \CAlg(\Shv(\Ran C), \otstar) \ar[l, "\oblv"] \ar[u, "\Delta^!"]  
\end{tikzcd}
\]
Hence the diagram of left adjoints commutes:
\begin{equation}\label{eq: free fact alg}
\begin{tikzcd}
\Shv(C) \ar[r, "\Sym^!"]  \ar[d, "\Delta_!"] &  \CAlg(\Shv(C), \sotimes)   \ar[d, "\Fact"]   \\
\Shv(\Ran C) \ar[r, "\Sym^\star"] & \CAlg(\Shv(\Ran C), \otstar) 
\end{tikzcd}
\end{equation}
where $\Sym^?$ refers to the formation of the symmetric algebra with respect to the respective symmetric monoidal structure. Let $M \in \Shv(C)$ and $\cA := \Sym^!(M) \in \CAlg(\Shv(C), \sotimes)$. Then \eqref{eq: free fact alg} gives a natural isomorphism $\Fact(\cA) \cong \Sym^{\star} (\Delta_! M)$. Explicitly, 
\begin{align*}
\Sym^\star(\Delta_! M) = \bigoplus_{n>0} \un_! (((\Delta_! M)^{\boxtimes n})_{\Sigma_n}) = \bigoplus_{n > 0} \Delta_{n!} ((M^{\boxtimes n})_{\Sigma_n}),
\end{align*}
where $\Delta_n$ is the composition 
\[
\begin{tikzcd}
C^n \ar[r, "\Delta^n"] \ar[dr, "\Delta_n"'] & (\Ran C)^n \ar[d, "\un"] \\
& \Ran C
\end{tikzcd}
\]
and the coinvariants $(-)_{\Sigma_n}$ are for the natural action of the symmetric group $\Sigma_n$ on $n$ elements. See \cite{Ho17}, especially \S 3.5 of \emph{loc. cit.}, for more discussion on free factorization algebras. 
\end{example}

\subsubsection{Chiral algebras} 
Consider the correspondence
\[
\begin{tikzcd}
& (\Ran C \times \Ran C)_{\disj} \ar[dl] \ar[dr, "\un"] \\
\Ran(C) \times \Ran(C) & & \Ran(C)
\end{tikzcd}
\]
Then $!$-pullback and $*$-pushforward through this diagram equips $\Shv(\Ran(C))$ with another symmetric monoidal structure called the \emph{chiral tensor structure}, which we denote $\otch$. For us, $\otch$ is an intermediary technical device.

As usual, in any symmetrical monoidal category we have an operadic notion of associative/commutative (co)algebras and Lie algebras. 
Below we summarize the results on ``chiral Koszul duality'' from \cite{FG12} that we will need. First, we have the following algebraic structures. 
\begin{itemize}
\item (\emph{Lie$^\star$ algebras on $\Ran C$}) We denote by $\Lie^\star(\Ran C)$ the category of Lie algebras in the symmetric monoidal category $(\Shv(\Ran C), \otstar)$. 
\item (\emph{Lie$^\star$ algebras on $C$}) We denote by $\Lie^\star(C) \subset \Lie^\star(\Ran C)$ the full subcategory of sheaves supported on $\Delta(C) \subset \Ran(C)$. 
\item (\emph{Chiral Lie algebras on $\Ran C$}) We denote by $\Lie^{\mrm{ch}}(\Ran C)$ the category of Lie algebras in $(\Shv(\Ran C), \otch)$.
\item (\emph{Commutative algebras on $\Ran C$}) We abbreviate by $\CAlg^\star(\Ran C)$ the category of commutative algebras in $(\Shv(\Ran C), \otstar)$. We denote by $\CAlg^{\mrm{ch}}(\Ran C)$ the category of commutative algebras in $(\Shv(\Ran C), \otch)$.
\item (\emph{Cocommutative coalgebras on $\Ran C$}) We denote by $\CoCAlg^\star(\Ran C)$ the category of cocommutative coalgebras in $(\Shv(\Ran C), \otstar)$. We denote by $\CoCAlg^{\mrm{ch}}(\Ran C)$ the category of cocommutative coalgebras in $(\Shv(\Ran C), \otch)$.
\end{itemize}

\begin{remark}\label{rem: diagonal Lie-*} It follows from the Cartesian square \eqref{eq: diag Ran} that $\Delta^!$ induces an equivalence between $\Lie^\star(C) $ and $\Lie(C) := \Lie(\Shv(C), \sotimes)$. See \cite[Proposition 4.2.3]{Ho17} for the details. 
\end{remark}

The identity functor on $\Shv(C)$ promotes to a lax symmetric monoidal functor 
\[
(\Shv(\Ran C), \otch)  \rightarrow 
(\Shv(\Ran C), \otstar).
\]
This formally induces functors
\begin{equation}
\CAlg^{\star} (\Ran C) \xleftarrow{\oblv} \CAlg^{\ch}(\Ran C),
\end{equation}
\begin{equation}\label{eq: oblv Lie}
\Lie^{\star}(\Ran C)  \xleftarrow{\oblv} \Lie^{\ch}(\Ran C),
\end{equation}
and
\begin{equation}
\CoCAlg^\star(\Ran C) \xrightarrow{\oblv} \CoCAlg^{\mrm{ch}}(\Ran C).
\end{equation}
The left adjoint to \eqref{eq: oblv Lie} is a functor
\[
\Lie^\star(\Ran C) \xrightarrow{\Ind_{\Lie}^{\star \rightarrow \mrm{ch}}} \Lie^{\mrm{ch}}(\Ran C). 
\]
The notion of \emph{cocommutative factorization coalgebras} in $\CoCAlg^{\mrm{ch}}(\Ran C)$ is defined in \cite[\S 2.4.6]{FG12}\footnote{Francis-Gaitsgory call it ``commutative factorization coalgebras''.}. The full subcategory spanned by such is denoted $\CoCAlg^{\mrm{ch}}_{\Fact}(\Ran C)$. Let us describe the definition informally for intuition. For a non-empty finite set $I$ and $\cF \in \Shv(\Ran C)$, let $\cF_I$ be the $!$-restriction of $\cF$ along $C^I \rightarrow \Ran C$. If $\cA \in \CoCAlg^{\mrm{ch}}(\Ran C)$, then $\cA$ has a system of maps $\cA_I \rightarrow \cA_{I_1} \sotimes \ldots \sotimes \cA_{I_n}$ for every (disjoint) partition $I = I_1 \sqcup \ldots \sqcup I_n$. We say that $\cA$ is a \emph{cocommutative factorization coalgebra} if all such maps restrict to isomorphisms on the open subspaces $(C^{I_1} \times \ldots \times C^{I_n})_{\disj}$.




\subsection{Graded factorization (co)algebras}\label{ssec: graded factorization algebras} Let $\Lambda$ be a free abelian monoid of finite rank, identified with $\Z_{\geq 0}^r$ via a choice of basis $\alpha_1, \ldots, \alpha_r$. Let $\Lambdapos := \Lambda \setminus \{0\}$. We will consider $\Lambdapos$-graded factorization algebras. These turn out to have more elementary incarnations, as explained in \cite[\S 4.1]{G16} and \cite[\S 4, 5]{GL19}. 

\subsubsection{Colored divisors} For $\lambda = \sum_i d_i \alpha_i \in \Lambda$, write
\[
C^{(\lambda)} := C^{(d_1)} \times C^{(d_2)} \times \ldots \times C^{(d_n)}.
\]
We define the space of \emph{$\Lambda$-colored divisors} $\Div^{\Lambdapos}(C)$ to be 
\[
\Div^{\Lambdapos}(C) = \coprod_{\lambda \in \Lambdapos} C^{(\lambda)}.
\]
It can also be described as the moduli space of $(D \in \Div C, \phi \co D \rightarrow \Lambdapos)$. We think of $\phi$ as equipping $D$ with a ``coloring'' by $\Lambdapos$, and sometimes refer to a point of $\Div^{\Lambdapos}(C)$ as a ``$\Lambdapos$-colored divisor''.

\subsubsection{Convolution and chiral tensor structures} Addition of divisors induces maps
\[
\add_{\lambda_1, \lambda_2} \co C^{(\lambda_1)} \times C^{(\lambda_2)} \rightarrow C^{(\lambda_1+\lambda_2)}
\]
which assemble into a map 
\[
\add \co \Div^{\Lambdapos}(C)  \times \Div^{\Lambdapos}(C) \rightarrow \Div^{\Lambdapos} C.
\]
This induces a symmetric monoidal structure on $\Shv(\Div^{\Lambdapos} C)$ that we call the \emph{convolution tensor structure} $\otstar$, with 
\[
\cF \otstar \cG := \add_!(\cF \boxtimes \cG).
\]

The disjoint locus $(\Div^{\Lambdapos} C \times \Div^{\Lambdapos} C)_{\disj}$ is defined as the open subspace of $(D_1, D_2)\in \Div^{\Lambdapos} C \times \Div^{\Lambdapos} C$ such that $D_1$ and $D_2$ have disjoint support. Then $!$-pull and $*$-push through the diagram 
\[
\begin{tikzcd}
& (\Div^{\Lambdapos} C \times \Div^{\Lambdapos} C)_{\disj} \ar[dl, hook] \ar[dr, "\un"] \\
\Div^{\Lambdapos} C \times \Div^{\Lambdapos} C & &\Div^{\Lambdapos} C
\end{tikzcd}
\]
induces another symmetric monoidal structure on $\Shv(\Div^{\Lambdapos} C)$ that we call the \emph{chiral tensor structure} $\otch$.

As explained in \cite[\S 4]{RasII}, these structures allow to imitate the theory of factorization algebras on $\Div^{\Lambdapos} C$ (instead of $\Ran C$). Moreover, while the pro-nilpotence of the chiral symmetric monoidal structure on $\Shv(\Ran(C))$ was a non-trivial theorem in \cite{FG12}, it is obvious in $\Shv(\Div^{\Lambdapos} C)$ due to the nature of the grading. Hence all the results of \cite{FG12} pass over to the graded setting, with the same proofs. 


\subsubsection{Graded factorization algebras}The notion of \emph{$\Lambdapos$-graded factorization algebras} is developed in \cite[\S 5]{GL19}. Informally, this consists of 
\begin{equation}\label{eq: graded fact algebra 0}
\cA = \{\cA^{\lambda}\}_{\lambda \in \Lambdapos} \in \Shv(\Div^{\Lambdapos} C) \cong \prod_{\lambda \in \Lambdapos} \Shv(C^{(\lambda)})
\end{equation}
equipped with a ``homotopy-compatible'' system of identifications 
\begin{equation}\label{eq: graded factorization algebra}
\cA^{\lambda_1+\lambda_2} |_{(C^{(\lambda_1)} \times C^{(\lambda_2)})_{\mrm{disj}}} \cong \cA^{\lambda_1} \boxtimes \cA^{\lambda_2}|_{(C^{(\lambda_1)} \times C^{(\lambda_2)})_{\mrm{disj}}}
\end{equation}
for all $\lambda_1, \lambda_2 \in \Lambdapos$. 

\begin{remark}[Ranification]
This definition looks rather different to the notion of factorization algebra in Definition \ref{def: factorization algebra}. The translation goes through a ``Ranification'' procedure explained in \cite[\S 4.7]{RasII}, which turns a factorization algebra $\cA = \{\cA^\lambda\}_{\lambda \in \Lambdapos}$ in the sense of \eqref{eq: graded factorization algebra} into a factorization algebra $\mrm{Ranification}(\cA) \in \Shv(\Ran C)^{\Lambdapos}$ in the sense of Definition \ref{def: factorization algebra}; we note that 
\[
\Delta^! \mrm{Ranification}(\cA)   = \prod_{\lambda \in \Lambdapos} \Delta_\lambda^! \cA_{\lambda} \in \Shv(C)^{\Lambdapos}. 
\] 
Thanks to Ranification, when working with $\Lambdapos$-graded factorization algebras we can replace $\Ran(C)$ with $\Div^{\Lambdapos} C$. 
\end{remark}

\subsubsection{Commutative graded factorization algebras}\label{sssec: commutative graded fact algebra}
We say that a $\Lambdapos$-graded factorization algebra $\cA := \{\cA^{\lambda}\}$ is \emph{commutative} if for all $\lambda_1, \lambda_2 \in \Lambdapos$ the maps 
\[
\cA^{\lambda_1} \boxtimes \cA^{\lambda_2}|_{(C^{(\lambda_1)} \times C^{(\lambda_2)})_{\disj}} \cong \cA^{\lambda_1+\lambda_2} |_{(C^{(\lambda_1)} \times C^{(\lambda_2)})_{\disj}
}
\] 
extend to maps\footnote{which may not be unique or ``canonical''} 
\begin{equation}\label{eq: commutative graded fact algebra}
\cA^{\lambda_1} \boxtimes \cA^{\lambda_2} \rightarrow \add^!_{\lambda_1, \lambda_2} \cA^{\lambda_1+\lambda_2}.
\end{equation}

For each $\lambda \in \Lambdapos$, let $\Delta_\lambda \co C \rightarrow C^{(\lambda)}$ be the diagonal map. Let
\[
\Delta^! := \prod_{\lambda \in \Lambdapos} \Delta_\lambda^! \co \Shv(\Div^{\Lambdapos} C) \rightarrow \Shv(C)^{\Lambdapos}.
\]

Parallel to Theorem \ref{thm: commutative fact alg}, we have: 

\begin{thm}\label{thm: commutative graded fact alg}
The functor $\cA \mapsto \Delta^! \cA$ induces an equivalence between the category of commutative factorization algebras in $\Shv^\star(\Div^{\Lambdapos} C)$ and commutative algebras in $\Shv(C)^{\Lambdapos}$. 
\end{thm}

\subsubsection{Graded factorization coalgebras}\label{sssec: graded factorization coalgebras} 
The case of (cocommutative) coalgebras graded by $\Lambdapos$ will play an important role for us. A reference is \cite{GWA}, although Gaitsgory calls ``cocommutative factorization algebras'' what we call ``cocommutative factorization coalgebras''. Indeed, a \emph{factorization coalgebra} on $\Div^{\Lambdapos} C$ is by definition the same as a factorization algebra: informally speaking, a collection 
\[
\cA = \{\cA^{\lambda}\}_{\lambda \in \Lambdapos} \in \Shv(\Div^{\Lambdapos} C) \cong \bigoplus_{\lambda \in \Lambdapos} \Shv(C^{(\lambda)})
\]
equipped with a ``homotopy-compatible'' system of identifications 
\begin{equation}\label{eq: graded factorization coalgebra}
\cA^{\lambda_1+\lambda_2} |_{(C^{(\lambda_1)} \times C^{(\lambda_2)})_{\mrm{disj}}} \cong \cA^{\lambda_1} \boxtimes \cA^{\lambda_2}|_{(C^{(\lambda_1)} \times C^{(\lambda_2)})_{\mrm{disj}}}.
\end{equation}
(The ``moral'' difference between graded factorization algebras and coalgebras is that the right side of \eqref{eq: graded factorization coalgebra} is $\add_{\lambda_1, \lambda_2}^*$ for factorization coalgebras, while the analogous expression for factorization algebras should use $\add_{\lambda_1, \lambda_2}^!$, but $\add_{\lambda_1, \lambda_2}$ is \'{e}tale when restricted to the disjoint locus, so these coincide.) 

We say that a factorization coalgebra $\cA$ is \emph{cocommutative} if the maps come by adjunction from maps
\begin{equation}\label{eq: cocommutative graded factorization coalgebra}
\cA^{\lambda_1 + \lambda_2} \rightarrow \cA^{\lambda_1} \otstar \cA^{\lambda_2}
\end{equation}
on $C^{(\lambda_1+\lambda_2)}$. If $\cA$ is holonomic, so that it lies in the domain of the partially defined adjoint $
\add^*$, then \eqref{eq: cocommutative graded factorization coalgebra} can be reformulated as saying that each \eqref{eq: graded factorization coalgebra} extends to a map
\[
\add_{\lambda_1, \lambda_2}^* \cA^{\lambda_1+\lambda_2}  \rightarrow \cA^{\lambda_1} \boxtimes \cA^{\lambda_2}
\]
on $C^{(\lambda_1)} \times C^{(\lambda_2)}$. 

\begin{remark}\label{rem: Verdier dual commutative graded fact alg} It is evident that Verdier duality induces an equivalence between holonomic commutative $\Lambdapos$-graded factorization algebras and holonomic cocommutative $\Lambdapos$-graded factorization coalgebras. Let us say that an object of $\Shv(C)^{\Lambdapos}$ is \emph{graded holonomic} if its component in each grading degree $\lambda \in \Lambdapos$ is holonomic. In turn, Theorem \ref{thm: commutative graded fact alg} gives an equivalence between holonomic commutative $\Lambdapos$-graded factorization algebras and graded holonomic commutative algebras in $\Shv(C)^{\Lambdapos}$.
\end{remark}

The map 
\begin{equation}
\Delta^* := \prod_{\lambda \in \Lambdapos} \Delta_\lambda^* \co \Shv_{\hol}(\Div^{\Lambdapos} C) \rightarrow \Shv_{\hol}(C)^{\Lambdapos}
\end{equation}
is symmetric monoidal for $\otstar$-monoidal structure on $\Shv_{\mrm{hol}}(\Div^{\Lambdapos} C)$ and the usual monoidal structure $\otimes$ on $\Shv_{\mrm{hol}}(C)^{\Lambdapos}$. It is clear that Verdier duality on $\Shv(\Div^{\Lambdapos} C)$ induces an equivalence between commutative algebras in $\Shv_{\hol}(C)^{\Lambdapos}$ and cocommutative coalgebras in $\Shv_{\hol}(C)^{\Lambdapos}$. Applying Remark \ref{rem: Verdier dual commutative graded fact alg} to Theorem \ref{thm: commutative graded fact alg} then gives: 

\begin{thm}\label{thm: cocommutative graded fact coalg} The functor $\cA \mapsto \Delta^* \cA$ induces an equivalence between the category of cocommutative factorization coalgebras in $\Shv_{\mrm{hol}}^\star(\Div^{\Lambdapos} C)$ and cocommutative cocolgebras in $(\Shv_{\hol}(C)^{\Lambdapos}, \otimes)$. 
\end{thm}


\begin{example}\label{ex: free graded fact algebra}
Let $V$ be a vector space and let $\ol{\Sym}(V) = \bigoplus_{n>0} \Sym^n(V)$ be the free graded (non-unital) commutative symmetric algebra on $V$. (The grading corresponds to the standard scaling $\G_m$-action on $V$.) By Example \ref{ex: free fact algebra}, the associated factorization algebra $\Fact(\ol{\Sym}(V)) = \{ \cA^{(n)} \in \Shv(C^{(n)})\}$ has 
\begin{equation}
 \cA^{(n)} := (\pi_{n!} \cA^{\boxtimes n})_{\Sigma_n}
\end{equation}
where $\pi_n \co C^n \rightarrow C^{(n)}$ is the quotient map and the coinvariants $(-)_{\Sigma_n}$ are for the natural action of the symmetric group $\Sigma_n$ on $n$ elements. 

Note that the graded dual of $\ol{\Sym}(V)$ is $\ol{\Sym}(V^*)$, which has a graded cocommutative coalgebra structure. We invert the grading so that $\ol{\Sym}(V^*)$ is graded by $\Z_{>0}$. Then the associated $\Z_{>0}$-graded cocommutative factorization coalgebra (under Theorem \ref{thm: cocommutative graded fact coalg}) is given levelwise by the Verdier dual of $\cA^{(n)}$.
\end{example}

\subsubsection{Factorization homology} Let $\cA  = \{\cA^\lambda\}$ be a $\Lambdapos$-graded factorization algebra on $\Div^{\Lambdapos} C$. Then we define its \emph{factorization homology} to be 
\[
\rR \Gamma_c^{\Fact}(C, \cA) := \bigoplus_{\lambda \in \Lambdapos} \rR \Gamma_c(C^{(\lambda)}, \cA^{\lambda}) \in \CAlg(\Vect_k).
\]
If $\cA$ is a commutative factorization algebra, then this is compatible with the previous definition of factorization homology \eqref{eq: factorization homology} under the equivalence of Theorem \ref{thm: commutative graded fact alg}. We will also apply this to cocommutative factorization coalgebras.


\subsection{Recollections on derived infinitesimal geometry}\label{ssec: derived infinitesimal geometry}

We will now study some constructions related to formal geometry. The theory of formal completion is complicated in the derived setting. For a closed embedding of \emph{classical} schemes $Z \inj Y$, one has the notion of ``$n$th infinitesimal neighborhood'' of $Z$ in $Y$, which is the closed subscheme of $Y$ defined by the $n$th power of the ideal sheaf of $Z$. Then the formal completion of $Z$ in $Y$ is the colimit over $n$ of the $n$th infinitesimal neighborhoods. However, in derived algebraic geometry it does not make sense to take the ``$n$th power of an ideal sheaf'', so it is subtler to form the $n$th infinitesimal neighborhood. 

In \cite[Chapter 9, \S 5]{GRII}, this problem is solved: for a map $f \co Z \rightarrow Y$ of derived stacks, Gaitsgory-Rozenblyum define a sequence of derived stacks
\begin{equation}\label{eq: derived thickenings}
Z = Z^{(0)} \rightarrow  Z^{(1)} \rightarrow  Z^{(2)}  \rightarrow  \ldots \rightarrow	 Z^{(n)}  \rightarrow  	\ldots \rightarrow Y.
\end{equation}
When $f$ is a closed embedding, then $Z^{(n)}$ is the derived analogue of the ``$n$th infinitesimal neighborhood of $Z$ in $Y$''. The extensions are produced using the theory of the cotangent complex, and the derived deformation to the normal cone. 

\begin{remark}\label{rem: non-LCI infinitesimal thickening}
Suppose $f$ is a closed embedding of classical schemes. If $f$ is a \emph{regular} embedding, then $Z^{(n)}$ is the usual $n$th infinitesimal neighborhood. In general, if $f$ is not a regular embedding then $Z^{(n)}$ can have non-trivial derived structure, even though $f$ is a map of classical schemes. In particular, the $Z^{(n)}$ constructed above may not agree with the classical notion of $n$th infinitesimal neighborhood. 
\end{remark}

\begin{defn}
The \emph{formal completion} of $Y$ along $f$ is $\colim_n Z^{(n)}$, the colimit formed in the category of derived stacks. We will denote it by $Y_f^{\wedge}$, or sometimes just $Y_Z^{\wedge}$ by abuse of notation. 
\end{defn}

\subsubsection{Ind-coherent sheaves on the formal completion}\label{sssec: sheaves on formal completion} For the notions of ind-coherent sheaves on derived ind-stacks and formal completions, we refer to \cite[\S 6,7]{GR14}. Suppose $Z \inj Y$ is a closed embedding of stacks. Let $\IndCoh_Z(Y) \subset \IndCoh(Y)$ be the full subcategory  of ind-coherent sheaves with set-theoretic support on $Z$. This admits a co-localization functor $\IndCoh(Y) \rightarrow \IndCoh_Z(Y)$. Letting $U = Y \setminus Z $ be the open complement of $Z$ and $j \co U \inj Y$ be its open embedding, the co-localization functor sends an ind-coherent sheaf $\cF$ to the derived fiber of $\cF \rightarrow j_* j^* \cF$. 

Let $\cZ := Y_Z^\wedge = \colim_n Z^{(n)}$ be the formal completion of $Y$ along $Z$. Let $\wh{i} \co \cZ \rightarrow Y$ be the colimit of the $i^{(n)} \co Z^{(n)} \rightarrow Y$. Then $\wh{i}^!$ factors through the co-localization map $\IndCoh(Y) \rightarrow \IndCoh_Z(Y)$, and induces an equivalence (cf. \cite[\S 7.4]{GR14})
\[
\begin{tikzcd}
\IndCoh(\cZ)  \ar[r, bend left, "\wh{i}_*"'] & \IndCoh_Z(Y) \ar[l, bend left, "\wh{i}^!"'] 
\end{tikzcd}
\]
where the left adjoint $\wh{i}_*$ is the colimit of the $\wh{i}^{(n)}_*$. The co-localization functor $\IndCoh(Y) \rightarrow \IndCoh_Z(Y) \cong \IndCoh(\cZ) $ is monoidal (with respect to $!$-tensor product). In particular, it sends the unit $\omega_Y \in \IndCoh(Y)$ to $\omega_{\cZ} \in \IndCoh(\cZ)$.

\subsubsection{Description as modules over a Lie algebroid} Let $f \co Z \rightarrow Y$ be any map of derived stacks. Then \cite[Chapter 8, Section 3.2]{GRII} constructs a \emph{Lie algebroid structure} on the relative tangent complex $\bT_{Z/Y}$. We denote this Lie algebroid by $\cL(Z/Y)$. We will only use the notion of Lie algebroid as a formal device used to access $\IndCoh(Y_f^\wedge)$, so we will not elaborate on it.

Recall that every derived stack has an underlying reduced sub-stack, and a \emph{nil-isomorphism} is a map that is an isomorphism on underlying reduced stacks \cite[Chapter 1, Definition 8.1.2]{GRI}. Equivalently, $f$ is representable in derived schemes and is an isomorphism on reduced subschemes after any base change from a derived scheme $S$ to $Y$.

Let $i \co Z\rightarrow \cZ$ be a nil-isomorphism, for example the natural map $Z \rightarrow \cZ := Y_Z^{\wedge}$ induced by a map $f \co Z \rightarrow Y$. Then by \cite[Chapter 5, \S 2.2.6, and Chapter 8, \S 4.1]{GRII} we have an equivalence
\begin{equation}\label{eq: lie-algebroid equivalence}
\cL(Z/\cZ)\Mod(\IndCoh(Z))  \cong \IndCoh(\cZ),
\end{equation}
such that the following diagram commutes: 
\[
\begin{tikzcd}
 \cL(Z/\cZ)\Mod(\IndCoh(Z)) \ar[d, "\mrm{oblv}"]  &   \IndCoh(\cZ)\ar[d, "i^!"]   \ar[l, "\sim"] \\
\IndCoh(Z)  \ar[r, equals] & \IndCoh(Z) 
\end{tikzcd}
\]

Since $\cL(Z/\cZ)$ is a Lie algebroid over $Z$, there is the \emph{Lie algebra homology} functor
\[
\Chev(\cL(Z/\cZ), -) \co \cL(Z/\cZ)\Mod(\IndCoh(Z))  \rightarrow \IndCoh(Z).
\]
It is the left adjoint of the functor $\triv \co \IndCoh(Z) \rightarrow \cL(Z/\cZ)\Mod(\IndCoh(Z))$ equipping an ind-coherent sheaf on $Z$ with the trivial $\cL(Z/\cZ)$-module structure. The functor $\Chev(\cL(Z/\cZ), -)$ can be modeled explicitly in terms of the bar construction. 

Suppose we have a retraction $\pi \co \cZ \rightarrow Z$. Then from the construction of the equivalence \eqref{eq: lie-algebroid equivalence}, the following diagram commutes: 
\begin{equation}
\begin{tikzcd}
 \cL(Z/\cZ)\Mod(\IndCoh(Z)) \ar[d, "{\Chev(\cL(Z/\cZ), -)}"', bend right] &   \IndCoh(\cZ) \ar[d, "\pi_*"', bend right]   \ar[l, "\sim"] \\
\IndCoh(Z)  \ar[u, bend right, "\triv"'] \ar[r, equals] &  \IndCoh(Z) \ar[u, bend right, "\pi^!"']
\end{tikzcd}
\end{equation}

\begin{example}\label{ex: distributions at a point}
Consider $\omega_Z \in \IndCoh(Z)$. Then $\omega_{\cZ} \cong \pi^! (\omega_Z) \in \IndCoh(\cZ)$ corresponds under \eqref{eq: lie-algebroid equivalence} to $\omega_Z \in \IndCoh(Z)$ equipped with the trivial action of $\cL(Z/\cZ)$, and then $\pi_*(\omega_{\cZ}) \in \IndCoh(Z)$ is naturally isomorphic to the Lie algebra homology $\Chev(\cL(Z/\cZ), \omega_{Z})$. 
\end{example}



\subsection{Factorization homology of distribution coalgebras}\label{ssec: factorization of distribution coalgebras} In this subsection we prove the key local-global principle that will be used to calculate spectral periods. In this subsection we focus on the de Rham setting, where in particular $\F = k$.

\subsubsection{Twisted mapping stack} Let $Y$ be an affine scheme acted upon by $\chG$ and $L$ be a $\chG$-local system on $C$, thought of as a $\chG$-bundle $L \rightarrow C_{\dR}$. Then we may form the space $\pi \co Y   \times^G L \rightarrow C_{\dR}$. The derived space of sections $\hSect(C_{\dR}, Y \times^{\chG} L)$ is defined in \S \ref{sssec: derived mapping stacks}; informally speaking, it parametrizes flat maps from $C$ to $Y$ twisted by the local system $L$. 

\subsubsection{Distribution coalgebras}\label{sssec: distributions}
For a closed subscheme $Z \inj Y$ we may form $\Gamma_Z(Y; \omega)$, the derived global sections of the dualizing sheaf of $Y$ with set-theoretic support in $Z$. The formation of $\Gamma_Z(Y; \omega)$ is functorial in maps of pairs $(Y,Z)$ which are proper in the second factor. 

In particular, let $y \in Y$ and consider $\Gamma_y(Y, \omega)$. Informally, it is the space of ``distributions on $Y$ supported at $y$''. The diagonal map $Y \xrightarrow{\Delta} Y \times Y$ equips $\Gamma_y(Y, \omega)$ with the structure of a unital cocommutative coalgebra, the unit coming from $  (Y,y) \rightarrow (\pt, \pt)$. We denote by $\ol{\Gamma}_y(Y; \omega)$ the corresponding non-unital cocommutative coalgebra (cf. \S \ref{ssec: non-unital}). 

Suppose $y \in Y$ is a $\chG$-fixed point. Then the $\chG$-action on $Y$ equips $\ol{\Gamma}_y(Y; \omega)$ with an action of $\chG$. Since the diagonal map is $\chG$-equivariant (as is the projection to a point), the commutative coalgebra structure is compatible with the $\chG$-action. This promotes $\ol{\Gamma}_y(Y; \omega)$ to a cocommutative coalgebra in $\Rep(\chG)$.

Suppose furthermore that $\ol{\Gamma}_y(Y; \omega)$ has a $\Z_{>0}$-grading.\footnote{For example, this happens  if $Y$ has a commuting $\G_m$-action for which $y$ is an attracting point (i.e., $y$ is a $\G_m$-fixed point and in some neighborhood $U \ni y$ the action extends to $U \times \A^1 \rightarrow U$ and $y$). In the applications to spherical varieties, the $\G_m$-action on $Y$ will be the $\Ggr$-action on $\chX$.} We may view $\ol{\Gamma}_y(Y; \omega)$ as a constant cocommutative coalgebra in $\Shv(C)^{\Z_{>0}}$ by $*$-pullback, and then by Theorem \ref{thm: cocommutative graded fact coalg} there is an associated $\Z_{>0}$-graded cocommutative factorization coalgebra, which we denote $\Fact(\ol{\Gamma}_y(Y; \omega)) \in \CAlg_{\Fact}^{\star}(\Div^{\Z_{>0}} C)$. 

\subsubsection{Twisting a factorization algebra by local systems}\label{sssec: twist FA by local system} 
Recall that $\rB \chG$ is the classifying stack $[\Spec k/\chG]$. In \cite[\S 4.2]{Gai15} the space $(\rB\chG)_{\Ran C}$ is defined. For a derived affine scheme $S$, the $S$-points of $(\rB\chG)_{\Ran C}$ is the groupoid of pairs of  $\ul{c} \in \Ran(C)_{\dR}(S)$ and a map 
\[
(\cD_{\ul{c}})_{\dR} \times_{S_{\dR}} S \rightarrow \rB \chG
\]
where $\cD_{\ul{c}}$ is the formal completion of $S \times C$ along its closed subscheme $\Gamma_{\ul{c}}$, defined as the union of the graphs of the maps $(S_{\cl})_{\red} \rightarrow C$ comprising $\ul{c}$. 

There is a convolution symmetric monoidal structure $\otstar$ on $\QCoh((\rB \chG)_{\Ran C})$, defined in \cite[\S 4.2.3]{Gai15}. In particular, when $\chG$ is the trivial group, so that $(\rB \chG)_{\Ran C} = \Ran C$, it specializes to the convolution symmetric monoidal $\otstar$ defined in \S \ref{ssec: fact homology} . 

A \emph{$\chG$-equivariant factorization algebra} is $\cA \in \CAlg(\QCoh((\rB \chG)_{\Ran C}), \otstar)$ such that 
\[
\oblv(\cA) \in \CAlg( \QCoh((\Ran C)_{\dR}), \otstar) \cong \CAlg( \Dmod(\Ran C), \otstar)
\]
is a factorization algebra. A $\chG$-equivariant factorization algebra can be twisted by a $\chG$-local system on $C$. To explain this, recall from \cite[\S 4.3]{Gai15} that there is an evaluation map 
\begin{equation}\label{eq: ran evaluation}
\ev \co \Ran(C)_{\dR} \times \LocG \rightarrow (\rB\chG)_{\Ran C}.
\end{equation}
For a derived affine scheme $S$, the map \eqref{eq: ran evaluation} sends $(\ul{c} \in \Ran(C)_{\dR}(S), L \in \LocG(S))$ to the pair of $\ul{c}$ and the composition of the maps 
\[
(\cD_{\ul{c}})_{\dR} \times_{S_{\dR}} S \xrightarrow{\pr_1} (\cD_{\ul{c}})_{\dR} \rightarrow C_{\dR} \xrightarrow{L} \rB \chG.
\] 
We denote by $\cA^{\univ}$ the $!$-pullback of $\cA \in \QCoh((\rB\chG)_{\Ran C})$ to $\Ran(C_{\dR}) \times \LocG$ along $\ev$. Informally, it is the universal twisting of $\cA$ by local $\chG$-local systems on $C$; restricting to $L \in \LocG$ gives the twist of $\cA$ by a particular $\chG$-local system $L$. 

\begin{example}
Let $y \in Y$ be as in \S \ref{sssec: distributions}. Since $\ol{\Gamma}_y(Y; \omega)$ is a cocommutative coalgebra in $\Rep(\chG)$, this construction applies to give $\Fact(\ol{\Gamma}_y(Y; \omega))^{\univ} \in \QCoh(\Ran(C)_{\dR} \times \LocG)$. 
\end{example}

\subsubsection{Local-global principle for distributions} 


Let $\LocY$ be the derived mapping stack $\Map(C_{\dR}, Y/\chG)$. The tautological map $Y/\chG \rightarrow \rB \chG$ induces a map 
\[
\pi \co \LocY \rightarrow \LocG. 
\]
The $\chG$-fixed point $y \in Y$ induces a distinguished section $\sigma_y  \co \LocG \rightarrow \LocY$ to $\pi$. Let $j \co U \inj \LocY$ be the complement of $\sigma_y$. A key formula for us will be the following ``local-global principle'' for calculating 
\[
\pi_{\sigma_y *}^\wedge (\omega_{\LocY}) := \fib(\pi_{*} \omega_{\LocY} \rightarrow \pi_* j_* \omega_U),
\]
the direct image of $\omega_{\LocY}$ under $\pi$, with supports along the section $\sigma_y$. 

\begin{prop}\label{prop: distribution algebra} Assume that the cotangent complex of $Y$ is perfect (in particular, has bounded tor-dimension). Let $\pr_2 \co \Ran(C)_{\dR} \times \LocG \rightarrow \LocG$ be the second projection map. Then there is a natural isomorphism
\begin{equation}\label{eq: distribution algebra}
\ol{\pi_{\sigma_y *}^\wedge (\omega_{\LocY})}  \cong \pr_{2*}(\Fact(\ol{\Gamma}_y(Y; \omega))^{\univ}) \otimes \omega_{\LocG} \in \QCoh(\LocG)
\end{equation}
where on the left side $\ol{(-)}$ refers to removing the unit (\S \ref{ssec: non-unital}).
\end{prop}



\begin{remark}
The factor of $\omega_{\LocG}$ on the right side of \eqref{eq: distribution algebra} comes from the fact that tensoring with $\omega_{\LocG}$ induces a fully faithful symmetric monoidal functor $(\QCoh(\LocG), \otimes) \rightarrow \IndCoh(\LocG, \sotimes)$. 
\end{remark}

\begin{proof}

Recall that $\bT_{y/Y} \cong \bT_yY[-1]$ has a canonical Lie algebra structure. This can be thought of as coming from the group structure of the derived loop space (i.e., inertia group) $y \times_Y y$, whose associated Lie algebra is $\bT_y Y[-1]$. Moreover, the $\chG$-action on $(Y,y)$ equips $\bT_yY[-1]$ with the structure of a Lie algebra in $\Rep(\chG)$. 

Regarding $\bT_{y/Y} \in \Shv(\pt)$, this equips its $!$-pullback $\ul{\bT_{y/Y}} \in \Shv(C)$ with the structure of a $\chG$-equivariant Lie algebra in $(\Shv(C), \sotimes)$, so we can view it as a $\chG$-equivariant Lie$^\star$ algebra on $C$ by Remark \ref{rem: diagonal Lie-*}. 


Let $\pr_2 \co C_{\dR} \times \LocG \rightarrow \LocG$. The relative tangent complex of the map $\sigma_y \co \LocG \rightarrow \LocY$ is
\[
\bT_{\sigma_y} \cong \bT_{\pi}[-1]  \cong    \pr_{2*}(\ul{\bT_{y/Y}^{\univ}})		
\]
where $\ul{\bT_{y/Y}^{\univ}} \in \QCoh(C_{\dR} \times \LocG)$ is the twist of $\bT_{y/Y}$ by the universal family of $\chG$-local systems. The $\chG$-equivariant Lie algebra structure on $\ul{\bT_{y/Y}}$ equips $\pr_{2*}(\ul{\bT_{y/Y}^{\univ}})$ with a Lie algebra structure in $\QCoh(\LocG)$. Then according to Example \ref{ex: distributions at a point}, we have a natural isomorphism 
\begin{equation}\label{eq: dist 1}
\pi_{\sigma_y *}^\wedge (\omega_{\LocY})  \cong \Chev(\pr_{2*}(\ul{\bT_{y/Y}^{\univ}}) ) .
\end{equation}

For sanity of notation we will present the rest of the proof over the trivial local system only. The full version over $\LocG$ can be done by repeating the argument below (and the theory of chiral algebras) with $C_{\dR}$ replaced by $C_{\dR} \times \LocG$, $\Ran(C)_{\dR}$ replaced by  $\Ran(C)_{\dR} \times \LocG$, etc. 

A formula of Beilinson-Drinfeld, reproved in \cite[Proposition 6.3]{FG12} via chiral Koszul duality\footnote{The two actually differ by the equivalence of removing the unit, cf. \cite[Remark 6.4.6]{FG12}. We are using the version of Francis-Gaitsgory.}, gives for any $\Lie^\star$ algebra $\cL$ on $C$ a natural isomorphism\footnote{The formulation below is the one of \cite[Proposition 4.1.11]{Ho17}.}
\[
\Chev(\Gamma_{\dR}(C;\cL)) \cong \rR \Gamma_c(\Ran C; \Chev(\Ind_{\Lie}^{\star \rightarrow \mrm{ch}} \cL)).
\]
Applying this to $\cL = \ul{\bT_{y/Y}}$, we obtain a natural isomorphism 
\begin{equation}\label{eq: dist 2} 
\Chev( \Gamma_{\dR}(C;  \ul{\bT_{y/Y}})) \cong \rR \Gamma_c(\Ran C; \Chev(\Ind_{\Lie}^{\star \rightarrow \mrm{ch}} \cL)).
\end{equation}
Another application of Example \ref{ex: distributions at a point} to $y \in Y$ gives a natural isomorphism between the cocommutative coalgebra $\ol{\Gamma}_y(Y; \omega)$ and $\Chev(\bT_{y/Y})$. Comparing \eqref{eq: dist 1} and \eqref{eq: dist 2}, we see that it suffices to construct a natural isomorphism $\Chev(\Ind_{\Lie}^{\star \rightarrow \mrm{ch}} \cL) \cong \Fact(\Chev(\bT_{y/Y}))$. This follows from \cite[Lemma 6.2.6]{FG12}, which we note is dual to \cite[Proposition 2.3.5]{GWA} (which is proved by the same general argument of \cite[Lemma 6.2.6]{FG12}). 


\end{proof}

\section{Geometric Langlands equivalence}

Let $G$ be a split reductive group. In the de Rham setting, the Geometric Langlands equivalence is an equivalence of categories 
\begin{equation}\label{eq: GLC}
\LL_G \co \Dmod(\Bun_G) \xrightarrow{\sim} \IndCoh_{\Nilp}(\LocG)
\end{equation}
characterized by various properties \cite[Conjecture 3.4.2]{Gai15}. Some of the properties (including all those that we will use later) are listed in this section. For $\GL_1$ this has been known since the work of Rothstein \cite{Roth96} and Laumon \cite{Lau96} (although in a slightly different formulation), and recently a proof for general groups has been announced at \cite{GLC}.

We will use the existence of \eqref{eq: GLC}, satisfying the properties mentioned in this section, as a black box -- that is to say, we will not use any internal details of the proof. Furthermore, we will only use \eqref{eq: GLC} for $G = \GL_1$ (where it is classical) and $G = \GL_2$ (where the proof was outlined several years ago in \cite{Gai15}). 

\begin{remark}
For all the compatibilities below, 
analogous formulations are true in the Betti \cite{BZN18} and \'etale \cite{sixI} contexts, though perhaps not as precisely documented in some respects. In order to save ourselves the trouble of having three separate discussions for each point, we will focus on the formulations in the \emph{de Rham} context in this section.  
\end{remark}

\subsection{Whittaker compatibility}\label{sssec: Whittaker compatibility}
Fix a Borel subgroup $B <G$ with unipotent radical $N$, and write $T := B/N$. This induces a natural map $\Bun_B \rightarrow \Bun_T$, and for $\cF_T \in \Bun_T(\F)$ we write $\Bun_N^{\cF_T}$ for the fibered product 
\[
\begin{tikzcd}
\Bun_N^{\cF_T} \ar[r] \ar[d] & \Bun_B \ar[d] \\
\pt \ar[r, "\cF_T"] & \Bun_T
\end{tikzcd}
\]
Let $\{\alpha_i\}_{ i \in \cI}$ be the simple (positive) roots corresponding to $N$. We write $\cF_T^{\alpha_i}$ for the line bundle obtained from $\cF_T$ via the homomorphism $\alpha_i \co T \rightarrow \G_m$. There is an evaluation morphism 
\[
\ev \co \Bun_N^{\cF_T} \rightarrow \prod_{i \in \cI} \ul{\rH^1(C, \cF_T^{\alpha_i})}
\]
constructed in \cite[\S 4.1.1]{FGV01}, where for a vector space $V/\F$ we write $\ul{V}$ for the associated affine space regarded as an $\F$-scheme. 

Let $X := G/U$ regarded as a $G \times T$-space via $(g,t) \co Ux \mapsto Ut^{-1} x g$, with the $\Ggr$ acting through the character $\G_m \xrightarrow{(1, e^{-2 \rho})  } G \times T$. Now regard $X$ as a $G$-space with the same action of $\Ggr$. Then $\BunX \rightarrow \Bun_G$, defined in \cite[\S 10.2]{BZSV}, can be identified with $\Bun_N^{\cF_T}$ such that $\cF_T^{\alpha_i} \cong \Omega_C$ for each $i$. Consider the composite map
\begin{equation}\label{eq: whit ev}
\BunX \xrightarrow{\ev} \prod_{i \in \cI} \underline{\rH^1(C, \Omega)} \xrightarrow{\mrm{add}} \underline{ \rH^1(C, \Omega)} \cong \A^1.
\end{equation}

\begin{example}\label{ex: whit GL_1} 
For $G = \GL_1$, the classical truncation of the stack $\BunX$ has $R$-points the groupoid of pairs 
\[
(\cL \in \Bun_G(R), t \co  \cO_C \boxtimes R \rightarrow \cL).
\] 
\end{example}

\begin{example}\label{ex: whit GL_2}
For $G = \GL_2$, the classical truncation of the stack $\BunX$ has $R$-points the groupoid of extensions
\[
\Omega_C^{1/2} \boxtimes R \rightarrow \cF \rightarrow \Omega_C^{-1/2} \boxtimes R.
\]
The map $\ev \co \BunX \rightarrow \A^1$ sends such an extension to its class in $\rH^1(C_R,\Omega_C \boxtimes R) \cong R$. 
\end{example}

Recall the Artin-Schreier sheaf on $\A^1$ from \S \ref{sssec: AS}. We let $\Psi  \in \Shv(\BunX)$ be the $*$-pullback of the Artin-Schreier sheaf along \eqref{eq: whit ev}. The \emph{(automorphic) Whittaker sheaf} $\sW \in \Shv(\Bun_G)$ is then the $!$-pushforward\footnote{In the de Rham case, note that this is defined because $\Psi$ is holonomic \cite[p.47]{Gai15}.} of $\Psi$ along the natural map $\BunX \rightarrow \Bun_G$. 


Up to twist, the Whittaker sheaf corresponds to the dualizing sheaf of $\LocG$ under $\LL_G$. More precisely, we normalize the Geometric Langlands Equivalence as in \cite[\S 12.2.1]{BZSV} so that for $\rho$ the half-sum of positive roots of $G$ and $\check{\rho}$ the half-sum of positive roots for $\chG$, we have 
\begin{equation}\label{eq: whittaker normalization}
\Shv(\Bun_G) \ni \sW\tw{(g-1) (\dim N - \tw{2\rho, 2 \check{\rho}})} \xmapsto{\LL_G} \omega_{\LocG} \tw{-(g-1) \dim G} \in \QCoh(\LocG) \subset \IndCoh_{\Nilp}(\LocG).
\end{equation}
Note that $\dim \Bun_G = (g-1) \dim G$.

\subsection{Hecke compatibility} We summarize \cite[\S 4]{Gai15}, which formulates the ``Hecke compatibility'' property of the Geometric Langlands Equivalence $\LL_G$. 

\subsubsection{Spectral localization}We will define a certain \emph{spectral localization} functor 
\[
\Loc^{\spec} \co \QCoh((\rB \chG)_{\Ran C}) \rightarrow \QCoh(\LocG)
\]
following \cite[\S 4.3]{Gai15}. Recall from \S \ref{sssec: twist FA by local system} the ``evaluation map'' $\ev \co \Ran(C)_{\dR} \times \LS \rightarrow (\rB\check{G})_{\Ran C}$. The diagram 
\[
\begin{tikzcd}
\Ran(C)_{\dR} \times \LS  \ar[r, "\mrm{ev}"] \ar[d, "\pr_2"] &  (\rB\check{G})_{\Ran C} \\
\LS 
\end{tikzcd}
\]
defines the functor
\[
\Loc^{\spec}  = \pr_{2*} \ev^!  \co \QCoh((\rB\check{G})_{\Ran C}) \rightarrow \QCoh(\LS).
\]

\subsubsection{Hecke operators} The Ran version of the Hecke stack, denoted $\Hk_{G,\Ran C}$, has $R$-points the groupoid of tuples $(\ul{c}, \cF_1, \cF_2, \beta)$ where $\ul{c} \in \Ran C (R)$, $\cF_i \in \Bun_G(R)$, and $\beta$ is an isomorphism $\cF_1|_{C_R \setminus \ul{c}} \cong \cF_2|_{C_R \setminus \ul{c}}$. There is a correspondence diagram 
\begin{equation}\label{diag: Hk Ran}
\begin{tikzcd}
& \Hk_{G,\Ran C} \ar[dl, "h_1"'] \ar[dr, "\pr \times h_2"] \\
\Bun_G & & \Ran C \times \Bun_G
\end{tikzcd}
\end{equation}

Convolution defines a map 
\begin{equation}\label{eq: Hk Ran conv}
\Hk_{G,\Ran C}\times_{h_2, \Bun_G , h_1} \Hk_{G,\Ran C} \rightarrow \Hk_{G,\Ran C}.
\end{equation}
Pull-push via the correspondence diagram 
\[
\begin{tikzcd}
\Hk_{G,\Ran C} \times_{h_2, \Bun_G , h_1} \Hk_{G,\Ran C} \ar[r] \ar[d, "\eqref{eq: Hk Ran conv}"] &\Hk_{G,\Ran C} \times \Hk_{G,\Ran C} \\
\Hk_{G,\Ran C} \times_{\Ran C} (\Ran C \times \Ran C) \ar[d,"\Id \times \un"]\\
\Hk_{G,\Ran C}
\end{tikzcd}
\]
defines a monoidal structure on $\Shv(\Hk_{G,\Ran C})$. 

The correspondence diagram \eqref{diag: Hk Ran} defines an action of the monoidal category $\Shv(\Hk_{G,\Ran C})$ on $\Shv(\Bun_G)$. We refer to this as the ``action of Hecke operators''. 

\subsubsection{Geometric Satake}
The Geometric Satake equivalence \cite{MV07} induces a functor 
\begin{equation}\label{eq: ran geometric satake}
\QCoh((\rB \chG)_{\Ran C})  \xrightarrow{\Sat} \Shv(\Hk_{G, \Ran C}).
\end{equation}

The \emph{Hecke compatibility property} of the Geometric Langlands equivalences stipulates that $\LL_G$ intertwines the action of $\QCoh((\rB \chG)_{\Ran C})$ on $\Dmod(\Bun_G)$ (via \eqref{eq: ran geometric satake} and Hecke operators) and on $\IndCoh_{\Nilp}(\LS)$ (via spectral localization).

\subsubsection{Graded version} By the discussion in \S \ref{ssec: graded factorization algebras}, for a $\Lambdapos$-graded object in $\QCoh((\rB \chG)_{\Ran C})$, there is an alternative formulation of spectral localization in terms of graded configuration space $\Div^{\Lambdapos} C$ instead of $\Ran C$. 

There is a parallel story for the Hecke action, using $\Hk_{G, \Div^{\Lambdapos} C}$ instead of $\Hk_{G, \Ran C}$. The stack $\Hk_{G, \Div^{\Lambdapos} C}$ has $R$-points the groupoid of tuples $(D, \cF_1, \cF_2, \beta)$ where $D \in \Div^{\Lambdapos} C(R)$, $\cF_i \in \Bun_G(R)$, and $\beta$ is an isomorphism $\cF_1|_{C_R \setminus D} \cong \cF_2|_{C_R \setminus D}$. The graded version of \eqref{eq: ran geometric satake} is a functor
\begin{equation}\label{eq: graded ran geometric satake}
\QCoh((\rB \chG)_{\Div^{\Lambdapos} C}) \xrightarrow{\Sat} \Shv(\Hk_{G, \Div^{\Lambdapos} C}).
\end{equation}
Hecke compatibility then says that $\LL_G$ intertwines the action of $\QCoh((\rB \chG)_{\Div^{\Lambdapos} C}) $ on $\IndCoh_{\Nilp}(\LocG)$, via spectral localization and tensoring, with its action on $\Dmod(\Bun_G)$ via \eqref{eq: graded ran geometric satake} and Hecke operators. 

\begin{remark}\label{rem: unit}
By definition, the spectral localization functor $\Loc^{\spec}$ takes a factorization (co)algebra in $\QCoh((\rB \chG)_{\Div^{\Lambdapos} C})$ to its relative factorization homology over $\LocG$. Since the factorization homology takes the unit to the unit (cf. \cite[Remark 6.4.6]{FG12}), it will be convenient to formally define factorization homology on \emph{unital augmented} (co)algebras by taking the direct sum with a copy of the unit (see \S \ref{ssec: non-unital}). 
\end{remark}

\begin{example}[Hecke compatibility for the standard representation]\label{ex: Hecke compatibility}
Let $G = \GL_r$, so $\chG = \GL_r$. We denote by $\Std \in \Rep(\chG)$ the standard representation. Consider 
\begin{equation}
\ol{\Sym}^{\bu}(\Std) := \bigoplus_{n>0} \Sym^n(\Std) \in \Rep(\chG),
\end{equation}
the free (non-unital) graded cocommutative coalgebra in $\Rep(\chG)$ generated by $\Std$ in graded degree 1. There is a corresponding (under Theorem \ref{thm: cocommutative graded fact coalg}) cocommutative factorization coalgebra 
\begin{equation}
\Fact(\ol{\Sym}^{\bu}(\Std))  = \left\{
\Fact(\ol{\Sym}^{\bu}(\Std)) \in  \QCoh((\rB \chG)_{C^{(n)}}) \right\}_{n>0}.
\end{equation}
By Example \ref{ex: free graded fact algebra}, it is given explicitly by 
\begin{equation}
\Fact(\ol{\Sym}^{\bu}(\Std))_n \cong \pi_{n!} ( \Std^{\boxtimes n})_{\Sigma_n}
\end{equation}
where $\pi_n \co C^n \rightarrow C^{(n)}$ is the natural quotient and $\Sigma_n$ is the symmetric group on $n$ elements. 

Then the universally twisted (cf. \S \ref{sssec: twist FA by local system}) version
\begin{equation}
\Fact(\ol{\Sym}^{\bu}(\Std))^{\univ}  =: \{
\Fact(\ol{\Sym}^{\bu}(\Std))^{\univ}_n \in  C^{(n)}_{\dR} \times \LocG \}_{n >0}
\end{equation}
is described as follows. We have a universal $\chG$-local system on $C \times \LocG$, and composing it with the standard representation gives $L^{\univ} \in \QCoh(C_{\dR} \times \LocG)$. Taking the $n$-fold exterior product of $L^{\univ}$ in the $C_{\dR}$ factor gives $L^{\univ}_n \in \QCoh(C^n_{\dR} \times \LocG)$. Then for $\pi_n \co C^n \rightarrow C^{(n)}$ the obvious projection, we have 
\begin{equation}
\Fact(\ol{\Sym}^{\bu}(\Std))^{\univ}_n  \cong (\pi_{n!} L^{\univ}_n)_{\Sigma_n} \in \QCoh(C^{(n)}_{\dR} \times \LocG).
\end{equation}

Note that under the Geometric Satake equivalence for $\GL_r$, each $\Sym^d(\Std) \in \Rep(\chG)$ corresponds to a shift of the \emph{constant sheaf} on its support in $\Hk_{G}$ for every $d \in \Z_{\geq 0}$. (This is well-known; the key point is that each non-zero weight space of $\Sym^d(\Std)$ is one-dimensional.) Consequently, $L^{\univ}_n$ corresponds to the constant sheaf on $\Hk_{G, C^n}$ normalized to be relatively perverse (in the sense of \cite{HS23}) over $C^n$, namely $k\tw{n(r-1)}$. Hence \eqref{eq: graded ran geometric satake} takes 
\begin{equation}\label{eq: non-unital free fact algebra}
\left\{  \Fact(\ol{\Sym}^{\bu}(\Std)   )_n \in  \QCoh((\rB \chG)_{C^{(n)}})  \right\}_{n>0} \mapsto \left\{k\tw{n(r-1)} \in \Shv(\Hk_{G, C^{(n)}}) \right\}_{n>0}.
\end{equation}

Now let $\Sym^{\bu}(\Std) \in \Rep(\chG)$ be the \emph{unital} cocommutative coalgebra $k \oplus \ol{\Sym}^{\bu}(\Std)$. Then we may formally set 
\begin{equation}
\Fact(\Sym^{\bu}(\Std)) = \{\Fact(\Sym^{\bu}(\Std))_n  \in \QCoh ((\rB\chG)_{C^{(n)}})\}_{n \geq 0}
\end{equation}
to be 
\begin{equation}\label{eq: augment spectral fact algebra}
\Fact(\Sym^{\bu}(\Std))_n := \begin{cases} \pi_{n!} ( \Std^{\boxtimes n})_{\Sigma_n}   \in \QCoh((\rB\chG)_{C^{(n)}}) & n>0, \\ \cO_{\pt} \in \QCoh(\pt) & n=0.\end{cases}
\end{equation}

We also formally define the universally twisted version 
\begin{equation}
\Fact(\Sym^{\bu}(\Std))^{\univ} = \{\Fact(\Sym^{\bu}(\Std))_n^{\univ}  \in \QCoh (C_{\dR}^{(n)} \times \LocG )\}_{n \geq 0}
\end{equation}
to be 
\begin{equation}\label{eq: augment spectral universal fact algebra}
\Fact(\Sym^{\bu}(\Std))_n^{\univ} = \begin{cases} (\pi_{n!} L^{\univ}_n)_{\Sigma_n}   \in \QCoh(C^{(n)}_{\dR} \times \LocG) & n>0, \\ \cO_{\LocG} \in \QCoh(\LocG) & n=0.\end{cases}
\end{equation}
By definition, we have
\begin{equation}\label{eq: Loc of Sym}
\Loc^{\Spec}(\Fact(\Sym^{\bu}(\Std))) := \pr_{2*} \left( \prod_{n \geq 0} \Fact(\Sym^{\bu}(\Std))_n^{\univ} \right) \in \QCoh(\LocG).
\end{equation}

Let $\Hk_{G, C^{(n)}}^{\Std}$ be the stack of tuples $(D, \cF_1, \cF_2, \beta)$ where $\beta \co \cF_1 \inj \cF_2$ is an upper modification such that $\coker(\beta)$ is a line bundle over $D$. Let 
\begin{equation}\label{eq: Hecke Std}
\Hk_{G, \Sym^{\bu} C} := \bigcup_{n \geq 0} \Hk_{G, C^{(n)}}^{\Std}.
\end{equation}
The object of $\Shv( \Hk_{G, \Sym^{\bu} C})$ corresponding to \eqref{eq: augment spectral fact algebra} under the Geometric Satake equivalence is
\begin{equation}\label{eq: augment automorphic fact algebra}
\Sat(
\Fact(\Sym^{\bu}(\Std))_n ) \cong \begin{cases} 
k\tw{n(r-1)} \in \Shv(\Hk_{G, C^{(n)}}^{\Std}) & n>0, \\ 
k \in \Shv( \Hk_{G, C^{(0)}} \cong \Bun_G) & n=0.
\end{cases} 
\end{equation}
Hecke compatibility then stipulates that the tensoring action of $\Loc^{\spec} (\Fact(\Sym^{\bu}(\Std))) \in \QCoh(\LocG)$ on $\IndCoh_{\Nilp}(\LocG)$ is intertwined under $\LL_G$ with the action of the kernel sheaf 
\begin{equation}\label{eq: Hk Std kernel sheaf}
\prod_{n \geq 0} \ulk \tw{n(r-1)} \in \Shv(\bigcup_{n \geq 0} \Hk_{G, C^{(n)}}^{\Std})
\end{equation}
on $\Dmod(\Bun_G)$ via convolving on the correspondence
\begin{equation}
\begin{tikzcd}
& \Hk_{G, \Sym^{\bu} C} := \bigcup_{n \geq 0} \Hk_{G, C^{(n)}}^{\Std} \ar[dl] \ar[dr] \\
\Bun_G & & \Bun_G
\end{tikzcd}
\end{equation}
\end{example}

\subsection{Eisenstein compatibility}\label{ssec: Eisenstein compatibility} Let $P \subset G$ be a standard parabolic subgroup with Levi quotient $M$. This induces a standard parabolic subgroup $\chP \subset \chG$ with Levi quotient $\chM$. The Eisenstein compatibility property of $\LL_G$ is discussed in \cite[\S 6]{Gai15}. We will briefly summarize it. 

With reference to the correspondence diagram 
\[
\begin{tikzcd}
& \Bun_P \ar[dr, "p"] \ar[dl, "q"'] \\
\Bun_M & & \Bun_G
\end{tikzcd}
\]
the \emph{Eisenstein functor} \cite[\S 6.3]{Gai15}
\[
\Eis_P \co \Shv(\Bun_M) \rightarrow \Shv(\Bun_G)
\]
is defined to be $\Eis_P := p_! q^*$. 

With reference to the correspondence diagram 
\[
\begin{tikzcd}
& \Loc_{\chP} \ar[dr, "\check{p}"] \ar[dl, "\check{q}"']  \\
\Loc_{\chM} & & \Loc_{\chG}
\end{tikzcd}
\]
the spectral Eisenstein functor 
\[
\Eis_{\chP}^{\spec} \co \IndCoh_{\Nilp}(\Loc_{\chM}) \rightarrow \IndCoh_{\Nilp}(\Loc_{\chG})
\]
is defined to be $\Eis_{\chP}^{\spec}:= \check{p}_* \circ \check{q}^!$. Note that our normalization is slightly different from \cite[\S 6.4]{Gai15}, in that we use $\check{q}^!$ instead of $\check{q}^*$; the two differ by tensoring by an invertible sheaf (namely, a shift of the determinant of $\bT_{\check{q}}$) since $\check{q}$ is quasi-smooth. 

Then the diagram 
\[
\begin{tikzcd}
\Dmod(\Bun_M) \ar[d, "\Eis_{\chP}"] \ar[r, "\LL_M", "\sim"'] & \IndCoh_{\Nilp}(\Loc_{\chM}) \ar[d, "\Eis_{\chP}^{\spec}"] &   \\
 \Dmod(\Bun_G)  \ar[r, "\LL_G"] & \IndCoh_{\Nilp}(\Loc_{\chG}) 
\end{tikzcd}
\]
commutes up to tensoring by an explicit invertible sheaf in either column; the details are pinned down in \cite[\S 6.4.8]{Gai15}. 

\begin{notat}Objects which do not lie in $\IndCoh_{\Nilp}(\LocG)$ arise naturally (an example that will come up for us is the IndCoh-pushforward of $\cO_{\pt}$ along the map $\pt \rightarrow \LocG$ corresponding to the trivial local system). There is a co-localization functor $\IndCoh(\LocG) \rightarrow \IndCoh_{\Nilp}(\LocG)$, which is right adjoint to the natural inclusion $\IndCoh_{\Nilp}(\LocG) \inj \IndCoh(\LocG) $. By convention, whenever we write $\Eis^{\spec}(\ldots) \in \IndCoh_{\Nilp}(\LocG)$, we are automatically applying this co-localization functor (if necessary). This has the same effect as applying the co-localization to the input object of $\IndCoh(\Loc_{\chM})$.
\end{notat}

\part{Examples}

\section{Tate period}

\subsection{Setup}We study an anzsatz that we call the ``Tate period'', as it corresponds to the integral representations in Tate's thesis \cite{Tate67} for the $L$-functions attached to Hecke characters. 

\subsubsection{Automorphic side} Let $X = \A^1$, with $G = \GL_1$ acting by the standard representation.

\subsubsection{Spectral side} Let $\chX = \A^1$, with $\chG = \GL_1$ acting by the standard representation.

\subsubsection{$\Ggr$-action} According to \cite[Example 12.6.7]{BZSV}, the ``neutral'' action of $\Ggr$ is the standard scaling action on both $X$ and $\chX$, so the ``un-normalized'' action of $\Ggr$ is trivial on both $X$ and $\chX$ (cf. \cite[Example 12.6.7]{BZSV}). Below we will take the un-normalized action of $\Ggr$.

\subsubsection{Automorphic period} 
By definition, $\BunX$ is the (derived) mapping stack 
\[
\BunX := \Map(C, \A^1/\G_m)
\]
with $R$-points the groupoid of pairs $(\cL, s)$ where: 
\begin{itemize}
\item $\cL$ is a line bundle on $C_R$. 
\item $s$ is a section of $\cL$ on $C_R$. 
\end{itemize} 

\begin{remark}
The derived structure plays no role until we study the functional equation in \S \ref{ssec: functional equation}, so we will pass to the classical truncation until then. 
\end{remark}


The map $\piaut \co \BunX \rightarrow \Bun_G$ forgets the datum of $s$. The automorphic period sheaf is
\[
\PX := \piaut_! (k_{\BunX}) \in \Shv(\Bun_G).
\]
In the language of \cite{BZSV}, this is the un-normalized period sheaf for the un-normalized action of $\Ggr$. 


\subsubsection{Spectral period} By definition, $\LocX$ is the derived mapping stack
\[
\LocX := \Map(C_{\dR}, \A^1/\G_m)
\]
is the derived stack with $R$-points the groupoid of pairs $(\cL, s)$ where:  
\begin{itemize}
\item $L \in \Loc_1(R)$ is a line bundle on $C_R$ with a flat connection along $C$ (informally speaking, an $R$-family of rank 1 local systems on $C$). 
\item $s \in \hSect(C_{\dR}, L)(R)$, or informally speaking a flat section of $L$. 
\end{itemize} 
The map $\pispec_*  \co \LocX \rightarrow \LocG$ forgets the datum of $s$. The (co-localized) spectral period sheaf is
\[
\LX := ( \pispec_* \omega_{\LocX})^{\shear} \in \QCoh(\LocG).
\]
Note that $\pispec_* \omega_{\LocX} = (  \pispec_* \omega_{\LocX})^{\shear}$ in this case, since the action of $\Ggr$ is trivial. In the language of \cite{BZSV}, this is the un-normalized $\cL$-sheaf for the un-normalized action of $\Ggr$. 

\subsubsection{Duality} The main theorem of this section is the following:

\begin{thm}[Geometric Langlands duality for the Tate period]\label{thm: Tate duality}
The Geometric Langlands equivalence for $G = \GL_1$ takes 
\[
\Dmod(\Bun_G) \ni  \PX \tw{g-1} \xmapsto{\LL_G}  \LX^{\shear} \in \IndCoh(\LocG).
\]
\end{thm}

It is explained in \cite[\S 12.2.3]{BZSV} that Theorem \ref{thm: Tate duality} is equivalent to \cite[Conjecture 12.1.1]{BZSV} for the neutral action of $\Ggr$. (The duality involution in \emph{loc. cit.} disappears under the translation from right actions to left actions -- see \cite[Appendix A]{CV}.)

\subsection{Spectral side}
Let $Z \inj \LocX$ be the zero section of $\pispec$ and $U$ be the open complement of $Z$. Note that $U$ is isomorphic to $\Map(C_{\dR}, \G_m/\G_m) = \pt$, while $Z$ is isomorphic to $\LocG$. We consider the exact triangle in $\QCoh(\LocG)$ coming from filtering by sections with set-theoretic support in $Z$: 
\begin{equation}\label{eq: Tate spectral triangle}
(\pispec)_{Z*}^{\wedge} (\omega_{\LocX}) \rightarrow \pispec_* (\omega_{\LocX}) \rightarrow \pispec_*|_U ( \omega_{\LocX}|_U)
\end{equation}
where we are using the notational conventions of \S \ref{ssec: singular support}. 

\subsubsection{Open stratum}\label{sssec: spectral open tate} We analyze the rightmost term of \eqref{eq: Tate spectral triangle}. It is evidently the direct image along $\pt \xrightarrow{\triv} \Loc_{1}$ at the trivial local system, i.e., the skyscraper sheaf $\cO_{\mrm{triv}}$ of the trivial local system.

\subsubsection{Closed stratum}\label{sssec: spectral closed tate} Before we analyze the leftmost term of \eqref{eq: Tate spectral triangle}, we establish some generalities. 

Recall that an affine derived scheme $S$ is \emph{$n$-coconnective} if its homotopy groups vanish in degrees greater than $n$, and \emph{eventually coconnective} if it is $n$-coconnective for some $n$ \cite[Chapter 2, \S 1.2]{GRI}. There is also a definition of eventually coconnective derived stacks, although it is a bit subtle: see \cite[Chapter 2, \S 2.6]{GRI}. The significance for us is that if a derived stack $S$ is eventually coconnective, then $\cO_S \in \Coh(S)$, so $\Perf(S) \subset \Coh(S)$, which ind-extends to a natural embedding $\Xi \co \QCoh(S) \inj \IndCoh(S)$. In particular, for a finite type map of derived stacks $f \co Y \rightarrow S$ such that $S$ is eventually coconnective, the \emph{relative dualizing complex} $\omega_{Y/S} := f^! (\cO_S) \in \IndCoh(Y)$ is defined.

\begin{lemma}\label{lemma: omega closed support}
Let $S$ be a locally finite type derived stack over a field and let $f \co Y_0 \inj Y$ be a closed embedding of finite type derived stacks over $S$. Let $\pi \co Y \rightarrow S$ the structure map to $S$. Denote by $\cY_0 := Y_{Y_0}^\wedge$ the formal completion of $Y$ along $Y_0$ and $\pi \co \cY_0 \rightarrow S$ be the restriction of $\pi$. Then we have isomorphisms 
\[
\pi_{Y_0*}^\wedge (\omega_{Y}) \cong \pi_* (\omega_{\cY_0}) \in \IndCoh(S).
\]
If $S$ is eventually coconnective, so that $\omega_{\cY_0/S} \in \IndCoh(Y)$ is defined, then we also have an isomorphism
\[
\pi_{Y_0*}^\wedge (\omega_{Y/S}) \cong \pi_* (\omega_{\cY_0/S}) \in \IndCoh(S).
\]
\end{lemma}

\begin{proof}By definition, $\pi_{Y_0*}^\wedge (\omega_{Y})$ is the composition of $\pi_*$ with the counit $\IndCoh(Y) \surj \IndCoh_{Y_0}(Y) \inj \IndCoh(Y)$. Letting $\wh{i} \co \cY_0 \rightarrow Y$ be the map induced by formal completion, the equivalence $\IndCoh(\cY_0) \xrightarrow{\sim} \IndCoh_{Y_0}(Y)$ is implemented by the adjunction $(\wh{i}_*, \wh{i}^!)$, cf. \S \ref{sssec: sheaves on formal completion}. Under this equivalence, we have 
\[
\pi_{Y_0*}^\wedge (\omega_{Y})   \cong \pi_* \wh{i}_* \wh{i}^! (\omega_{Y}) \cong \pi_* (\omega_{\cY_0}) \in \IndCoh(S)
\]
since $\wh{i}^! \omega_{Y} \cong \omega_{\cY_0}$, and similarly for the relative dualizing complex if $S$ is eventually coconnective. 
\end{proof}

We will use Proposition \ref{prop: distribution algebra}, with $Y := \chX$ and $y := 0 \inj \A^1 = Y$, to calculate $(\pispec)_{Z*}^{\wedge} (\omega_{\LocX}) $. For this we first need to compute the unital cocommutative coalgebra $\Gamma_0(\chX; \omega_{\chX})$. It will be convenient to establish a more general statement for future use.

\begin{lemma}\label{lem: formal thickening of dvb} Let $S$ be an eventually coconnective locally finite type derived stack over $k$ (hence of characteristic $0$), $\chG$ be a reductive group, and let $\pi \co E \rightarrow S$ be a $\chG$-equivariant derived vector bundle associated to a perfect complex $\cE \in \Perf(S/\chG)$ with tor-amplitude in $[0, \infty)$, where $\chG$ acts trivially on $S$. We view $S \inj E$ as the zero section. Equipping $E$ with the scalar $\G_m$-action, we have a natural isomorphism
\[
\pi_{S*}^{\wedge} (\omega_{E/S}) \cong \Sym^{\bu}(\cE) \in \QCoh(S/\chG)
\]
as $\G_m$-graded unital cocommutative coalgebras. Here $\Sym^{\bu}(\cE):= \bigoplus_{n \geq 0} \Sym^n(\cE)$ is the free graded $\chG$-equivariant unital cocommutative coalgebra over $S$ generated by $\cE$ (placed in graded degree $1$). 
\end{lemma}

\begin{remark}
Contrast the statement of Lemma \ref{lem: formal thickening of dvb} with the tautological statement 
\begin{equation}\label{eq: push O}
\pi_*(\cO_E) \cong \Sym^{\bu}(\cE^\vee)
\end{equation}
the free graded unital commutative algebra in $\Rep(\chG)$ generated by $\cE^\vee$ (placed in graded degree $-1$), which has a dual description. The appearance of this duality is somehow crucial. 
\end{remark}

\begin{proof} Let $z \co S \inj E$ be the zero-section. Let
\[
\begin{tikzcd}
S = S^{(0)} \ar[r, hook] \ar[d, "\pi^{(0)}"] &  S^{(1)} \ar[r, hook] \ar[d, "\pi^{(1)}"] &  \ldots  \ar[r, hook] \ar[d, "\pi^{(2)}"] &  E_S^\wedge \ar[r, hook] \ar[d] & E \ar[d, "\pi"] \\
S \ar[r, equals] & S \ar[r, equals]   & S\ar[r, equals]  & S \ar[r, equals]  &     S
\end{tikzcd}
\]
be the system of derived infinitesimal neighborhoods of $z \co S \inj E$ as in  \eqref{eq: derived thickenings}, with $E_S^\wedge := \colim_n S^{(n)}$. Then by Lemma \ref{lemma: omega closed support} we have $\pi_{S*}^\wedge (\omega_{E/S}) \cong \pi_* (\omega_{E_S^\wedge/S})$, and by definition of ind-coherent sheaves on ind-schemes we have (cf. \S \ref{sssec: sheaves on formal completion})
\begin{equation}\label{eq: omega colimit}
\pi_*(\omega_{E_S^\wedge/S}) = \colim_n \pi^{(n)}_*(\omega_{S^{(n)}/S}) \in \IndCoh(S).
\end{equation}
Since each $\pi^{(n)}$ is proper, Grothendieck duality applies to give a trace map $\pi^{(n)}_*(\omega_{S^{(n)}/S}) \rightarrow \cO_S$. Taking the colimit in $n$, we get a trace map 
\begin{equation}\label{eq: trace map}
\pi_{S*}^{\wedge} (\omega_{E_S^\wedge/S})   \rightarrow \cO_S.
\end{equation}
This induces a pairing 
\begin{equation}\label{eq: perfect pairing}
\pi_{S*}^{\wedge} (\omega_{E_S^\wedge/S}) \otimes_{\cO_S} \pi_* (\cO_E) \rightarrow 
\pi_{S*}^{\wedge} (\omega_{E_S^\wedge/S})  \xrightarrow{\eqref{eq: trace map}} \cO_S,
\end{equation}
which we claim is perfect in each degree of the grading. This means that the map induced by  \eqref{eq: perfect pairing},
\begin{equation}\label{eq: perfect duality}
\pi_* (\cO_E) \rightarrow \pi_{S*}^{\wedge} (\omega_{E_S^\wedge/S})^*,
\end{equation}
is an isomorphism, where we recall that the RHS denotes the $\cO_S$-linear dual of $\pi_{S*}^{\wedge} (\omega_{E_S^\wedge/S})$. In each graded degree, \eqref{eq: perfect duality} becomes a map of perfect complexes, so to check that it is an isomorphism on $S$, we can (by the finiteness hypotheses) check that it is an isomorphism after derived base change to every closed point of $S$. Over a point, $\cE$ splits, so we may assume that $\cE \cong \bigoplus_{i \geq 0} \cE^i[-i]$ where each $\cE^i$ is a locally free coherent sheaf, which induces a splitting $E \cong \bigoplus_{i \geq 0} E^i[-i]$. This then reduces to the case where $E$ is concentrated in a single degree $i \geq 0$. 

If $i > 0$, then $\pi \co E \rightarrow S$ is proper, as the classical truncation of $\pi$ is an isomorphism, and the claim follows directly from Grothendieck duality. We therefore reduce to the case $i=0$, where $E$ is a classical vector bundle over $S$. Returning to the situation of \eqref{eq: omega colimit}, since each $\pi^{(n)}$ is proper, Grothendieck duality applies to give a natural isomorphism  
\[
\pi^{(n)}_*(\omega_{S^{(n)}/S}) \cong (\pi^{(n)}_* \cO_{S^{(n)}})^{\vee} \in \IndCoh(S).
\]
Since $z$ is a regular embedding, $ \pi^{(n)}_* \cO_{S^{(n)}}$ is the quotient of $\Sym^{\bu}(\cE^\vee)$ by the $n$th power of its augmentation ideal. Dualizing this and taking the colimit over $n$ completes the proof that \eqref{eq: perfect pairing} is a perfect pairing. 

Finally we observe that the commutative coalgebra structure on $
\pi_{S*}^{\wedge} (\omega_{E_S^\wedge/S})$ is graded dual to the commutative algebra structure on $\pi_* (\cO_E)$ under \eqref{eq: perfect pairing}, so we conclude from \eqref{eq: push O}.  
\end{proof}

\begin{cor}\label{cor: tate spectral distribution}
The cocommutative coalgebra $\Gamma_0(\chX; \omega_{\chX}) \in \Rep(\GL_1)$ identifies with $\Sym^{\bu}(\Std)$, the free unital cocommutive coalgebra in $\Rep(\GL_1)$ on $\Std$ placed in graded degree $1$. 
\end{cor}

\begin{proof}
Apply Lemma \ref{lem: formal thickening of dvb} with $\chG = \GL_1$, $S = \pt$, $X = \chX$. 
\end{proof}

\begin{prop}\label{prop: Tate Loc}
Let $\Fact(\Sym^{\bu}(\Std))$ be the graded unital cocommutative factorization coalgebra \eqref{eq: augment spectral fact algebra}. Then we have an isomorphism 
\[
(\pispec)_{Z*}^{\wedge} (\omega_{\LocX}) \cong \Loc^{\spec} (\Fact(\Sym^{\bu} (\Std)))  \otimes \omega_{\LocG} \in \QCoh(\LocG).
\]
\end{prop}

\begin{proof} Applying Proposition \ref{prop: distribution algebra} with $y$ being the origin inside the $\chG$-scheme $\Std$, placed in graded degree 1 (corresponding to the standard scaling action of $\Ggr$), we obtain an isomorphism
\begin{equation}\label{eq: Tate non-unital Loc 1} 
\ol{(\pispec)_{Z*}^{\wedge} (\omega_{\LocX})} \cong (\pr_{2*} \Fact(\ol{\Gamma}_0(\chX; \omega_{\chX}))^{\univ} ) \otimes \omega_{\LocG} \in \QCoh(\LocG).
\end{equation}
In Corollary \ref{cor: tate spectral distribution} we saw that $\ol{\Gamma}_0(\chX; \omega_{\chX}) \cong \ol{\Sym}^{\bu} (\Std) \in \Rep(\chG)^{\Z_{>0}}$. Inserting this into \eqref{eq: Tate non-unital Loc 1} yields an isomorphism
\begin{equation}\label{eq: Tate non-unital Loc 2} 
\ol{(\pispec)_{Z*}^{\wedge} (\omega_{\LocX})} \cong (\pr_{2*} \Fact(\ol{\Sym}^{\bu} (\Std))^{\univ})  \otimes \omega_{\LocG} \in \QCoh(\LocG).
\end{equation}
Recall from \S\ref{ssec: non-unital} that we may pass between non-unital cocommutative coalgebras and augmented unital graded cocommutative coalgebras by taking the direct sum with the unit. Applying this to \eqref{eq: Tate non-unital Loc 2}, we obtain an isomorphism 
\begin{equation}\label{eq: Tate unital Loc} 
(\pispec)_{Z*}^{\wedge} (\omega_{\LocX}) \cong (\pr_{2*} \Fact(\Sym^{\bu} (\Std))^{\univ}) \otimes \omega_{\LocG} \in \QCoh(\LocG)
\end{equation}
where the meaning of $\Fact(\Sym^{\bu} (\Std))^{\univ}$ is as in \eqref{eq: augment spectral universal fact algebra}. We conclude by using that 
\[
\pr_{2*} \Fact(\Sym^{\bu} (\Std))^{\univ} \cong \Loc^{\spec} (\Fact(\Sym^{\bu} (\Std))) \in \QCoh(\LocG)
\]
by definition (cf. \eqref{eq: Loc of Sym}).


\end{proof}

\subsubsection{Extension class}\label{ssec: Tate spectral summary} Putting together \S \ref{sssec: spectral open tate} and \S \ref{sssec: spectral closed tate}, we have produced an exact triangle in $\IndCoh(\LocG)$,
\begin{equation}\label{eq: tate spectral side}
\Loc^{\spec}(\Fact(\Sym^{\bu} (\Std)))\otimes \omega_{\LocG} \rightarrow \pi_{\spec *} (\omega_{\LocX}) \rightarrow \cO_{\triv}.
\end{equation}

Next we want to understand the connecting homomorphism 
\[
\cO_{\triv} \rightarrow \Loc^{\spec}(\Fact(\Sym^{\bu} (\Std))) \otimes \omega_{\LocG}[1] \cong \pi_{Z*}^{\wedge} (\omega_{\LocX})   [1].	
\]
Let us factor $\triv \co \pt \rightarrow \LocG$ as a composition 
\[
\begin{tikzcd}
\pt  \ar[d, "q"'] \ar[dr, "\triv"] \\
\rB \G_m \ar[r, hookrightarrow, "\iota"] & \LocG
\end{tikzcd}
\]
Since $\iota$ is a closed embedding, we have an adjunction $(\iota_*,\iota^!)$, hence 
\begin{align}\label{eq: spectral ext group 1}
\Hom_{\LocG}(\cO_{\triv} , (\pispec)_{Z*}^{\wedge} (\omega_{\LocX}) [1]) & \cong \Hom_{\LocG}(\iota_* q_* \cO_{\pt} , (\pispec)_{Z*}^{\wedge}\omega_{\LocX}  [1])  \nonumber  \\
& \cong  \Hom_{\rB \G_m}(q_* \cO_{\pt}, \iota^! (\pispec)_{Z*}^{\wedge} \omega_{\LocX} [1]).
\end{align}

We have a derived Cartesian square
\begin{equation}\label{eq: Tate spec ext diag 1}
\begin{tikzcd}
V \ar[d, "\pi"] \ar[r, "\iota"] & \LocX \ar[d, "\pispec"] \\
\rB\G_m \ar[r, "\iota"] & \LocG
\end{tikzcd}
\end{equation}
where $V := \Spec (\Sym^{\bu} ( \Hdr(C; k)^\vee))$ is the derived vector bundle associated to the perfect complex $\Hdr(C; k)$ over $k$. (In fact, the argument of \cite[Proposition 5.34]{FYZ2} shows that $\LocX$ is the derived vector bundle over $\LocG$ associated to the perfect complex $L \mapsto  \Hdr(C; L)$.) We may further base change along its zero-section $Z \inj \LocX$ as well as the formal completion 
$\cZ \rightarrow \LocX$ thereof, giving a tower of derived Cartesian squares
\begin{equation}\label{eq: Tate spec ext diag 2}
\begin{tikzcd}
\cV \ar[d] \ar[r, "\iota"] & \cZ  \ar[d] \\ 
V \ar[d, "\pi"] \ar[r, "\iota"] & \LocX \ar[d, "\pispec"] \\
\rB\G_m \ar[r, "\iota"] & \LocG
\end{tikzcd}
\end{equation}
We abuse notation by calling all the horizontal maps $\iota$. In particular, from the outer derived Cartesian square 
\begin{equation}
\begin{tikzcd}
\cV \ar[d, "\pi"] \ar[r, "\iota"] & \cZ  \ar[d, "\pispec"] \\ 
\rB\G_m \ar[r, "\iota"] & \LocG
\end{tikzcd}
\end{equation}
we get a base change natural isomorphism 
\begin{equation}
\iota^! \pispec_* (\omega_{\cZ}) \cong \pi_* \iota^! (\omega_{\cZ}) \cong \pi_*(\omega_{\cV}) \in \QCoh(\rB \G_m).
\end{equation}
Inserting this into \eqref{eq: spectral ext group 1} and recalling from Lemma \ref{lemma: omega closed support} that $\pispec_* (\omega_{\cZ}) \cong (\pispec)_{Z*}^\wedge (\omega_{\LocX})$, we have a natural isomorphism
\begin{equation}\label{eq: Tate spec ext eq 1}
\Hom_{\LocG}(\cO_{\triv} ,  (\pispec)_{Z*}^{\wedge} (\omega_{\LocX})   [1])  \cong \Hom_{\rB\G_m}(q_* \cO_{\pt}, \pi_* \omega_{\cV}[1]).
\end{equation}
Recall that the functor $q^*$ induces $\QCoh(\rB\G_m) \xrightarrow{\sim} \Rep(\G_m)$. Note that $q^! \cong q^*[1]$ since $q$ is a $\G_m$-torsor (hence smooth of relative dimension $1$, with trivial relative dualizing sheaf), so by base change and Lemma \ref{lem: formal thickening of dvb} we have 
\[
q^* \pi_* (\omega_{\cV}) \cong \Sym^{\bu} (\Hdr(C; k)) [-1] 
\]
where the $\G_m$-action is via the natural scaling action on $\Hdr(C, k)$. Putting this into \eqref{eq: Tate spec ext eq 1} gives 
\begin{align}\label{eq: Tate spectral ext}
\Hom_{\LocG}(\cO_{\triv} , (\pispec)_{Z*}^{\wedge} (\omega_{\LocX})   [1])  & \cong  \Hom_{\rB\G_m}(q_* \cO_{\pt}, \pi_* \omega_{\cV}  [1]) \nonumber \\
& \cong  \Hom_{\G_m}(\cO(\G_m), \Sym^{\bu} (\Hdr(C;k))[-1][1]) \nonumber \\
& \cong \prod_{n \geq 0} \rH^0_{\dR}(C; k) .
\end{align}
Hence the isomorphism class of the underlying lying object $\pispec_* (\omega_{\LocX} )\in \IndCoh_{\Nilp}(\LocG)$ is determined by the extension class $[\pispec_* (\omega_{\LocX})]$ of \eqref{eq: Tate spectral triangle} in each of the 1-dimensional spaces $\rH^0_{\dR}(C; k)$ indexed by $n \geq 0$ in \eqref{eq: Tate spectral ext}, modulo the action of $k^{\times}$. In other words, we have pinned down the isomorphism class of $\pispec_* (\omega_{\LocX}) \in \IndCoh_{\Nilp}(\LocG)$ once we know whether or not the extension \eqref{eq: Tate spectral triangle} is split for each $n$. This is addressed by the following Lemma. 

\begin{lemma}\label{lem: tate spectral connecting}
For each $n \geq 0$, the extension \eqref{eq: Tate spectral triangle} maps to the class of $1 \in \rH^0_{\dR}(C; k)$ in \eqref{eq: Tate spectral ext}. 
\end{lemma}

\begin{proof}
The statement can be proved after base change along $\iota \co \rB\G_m \rightarrow \LocG$. Examining the computation of extensions, we see that it also suffices to prove for the classical truncation $V_{\cl}$ of $\cV$ (since the derived part only contributes in higher degrees), which is $\A^1/\G_m$ with the standard action of $\G_m$. Its relative dualizing sheaf over $\rB\G_m$ is $\omega_{V_{\cl}} \cong \cO_{V_{\cl}}[1]$. Then the statement is that in the exact triangle
\[
\pi_{0*}^\wedge  (\omega_{V_{\cl}}) \rightarrow \pi_* (\omega_{V_{\cl}}) \rightarrow   \cO_{\triv}  \in \QCoh(\rB \G_m)
\]
the map $\rR^0 \pi_* ( \cO_{\triv})  \rightarrow \rR^1 \pi_{0*}^{\wedge} (\omega_{V_{\cl}})$ is an isomorphism in the $n$ graded component for each $n \geq 0$; and in fact is the identity map with respect to the canonical identification $\rR^1 \pi_{0*}^{\wedge} (\omega_{V_{\cl}}) \cong k$. Applying $q^!$ and using base change, this reduces to the explicit computation of the exact triangle 
\[
\Gamma_0(\A^1; \omega_{\A^1} ) \rightarrow  \Gamma(\A^1; \omega_{\A^1}) \rightarrow \Gamma(\G_m; \omega_{\G_m}).
\]
Since $\A^1$ and $\G_m$ are smooth, we have compatible identifications $ \omega_{\A^1} \cong \cO_{\A^1}[1]$ and $\omega_{\G_m} \cong \cO_{\G_m}[1]$. The assertion that the map is an isomorphism is clear from the explicit identification of $\Gamma(\A^1; \cO_{\A^1}) \rightarrow \Gamma(\G_m; \cO_{\G_m})$, which also provides the canonical identification $\rR^1\Gamma_0(\A^1; \omega_{\A^1} ) \cong \Std$, and then the statement that the extension class is the identity is tautological.  
\end{proof}

\subsection{Automorphic side}\label{ssec: Tate automorphic side} We let $i \co Z \inj \BunX$ be the zero section and $j \co U \inj \BunX$ its open complement.

We filter $\piaut_! (\ul{k}_{\BunX})$ using the open-closed exact triangle in $\Shv(\BunX)$,
\begin{equation}
j_! \ul{k}_U \rightarrow \ul{k}_{\BunX} \rightarrow i_* \ul{k}_Z,
\end{equation}
which after applying $\piaut_!$ gives an exact triangle in $\Shv(\Bun_G)$, 
\begin{equation}\label{eq: Tate aut exact triangle}
\piaut_! j_! (\ul{k}_{U}) \rightarrow \piaut_! (\ul{k}_{\BunX}) \rightarrow \piaut_! i_*  (\ul{k}_Z).
\end{equation}

\subsubsection{Closed stratum}
Since the $\piaut$ restricts to an isomorphism $Z \xrightarrow{\sim} \Bun_G$, we have 
\begin{equation}\label{eq: Tate aut closed}
\piaut_! (\ul{k}_Z) \cong \ul{k}_{\Bun_G} \in \Shv(\Bun_G).
\end{equation}

\subsubsection{Open stratum} Recall that for $G = \GL_1$, we have $\Bun_G = \coprod_{d \in \Z} \Bun_{G}^d$ where $\Bun_{G}^d$ is the connected component parametrizing line bundles of degree $d$. 

The intersection $U \cap \Bun_G^d$, which parametrizes degree $d$ line bundles plus a non-zero section, is isomorphic to $C^{(d)}$, by the map sending a section to its divisor. Hence we have
\[
U \cong \Sym^{\bu} C = \coprod_{d \geq 0} C^{(d)}.
\]
The restriction $\piaut|_U \co U \rightarrow \Bun_G$ is a version of the Abel-Jacobi map: it sends a divisor $D$ to the line bundle $\cO(D)$. Therefore we denote it $\AJ \co U \rightarrow \Bun_G$, and we write $\AJ^d$ for the restriction to the connected component $C^{(d)} \subset U$, which lands in $\Bun_G^d$. Hence we have 
\[
\piaut_!(\ul{k}_U) \cong \AJ_! (\ul{k}_U ) \cong \AJ_!( \ul{k}_{\Sym^{\bu} C}).
\]

\subsubsection{Extension class}\label{sssec: Tate aut summary} We have just produced an exact triangle in $\Shv(\Bun_G)$: 
\begin{equation}\label{eq: Tate aut side}
\AJ_! (\ul{k}_{\Sym^{\bu} C}) \rightarrow \piaut_! ( \ul{k}_{\BunX}) \rightarrow \ul{k}_{\Bun_G} .
\end{equation}
To analyze the extension class of \eqref{eq: Tate aut exact triangle}, we will use a trick that was explained to us by Akshay Venkatesh. Consider the exact triangle in $\Shv(\BunX)$,
\begin{equation}\label{eq: Tate dual excision}
i_! i^! (\ul{k}_{\BunX}) \rightarrow \ul{k}_{\BunX} \rightarrow j_* j^* (\ul{k}_{\BunX}).
\end{equation}
Now apply $\piaut_*$ to \eqref{eq: Tate dual excision}. Since $i$ is a closed embedding, we have $i_! = i_*$, so we get an exact triangle in $\Shv(\Bun_G)$,
\begin{equation}\label{eq: tate aut ext 1}
\piaut_*  i_* i^!  (\ul{k}_{\BunX})  \rightarrow \piaut_* (\ul{k}_{\BunX}) \rightarrow \piaut_* j_*  j^* (\ul{k}_{\BunX}).
\end{equation}
Under the identification $\piaut \co Z \cong \Bun_G$, we have $\piaut \circ i =  \Id$. There is a $\G_m$-action on $\BunX$ contracting the fibers to the zero section, which by the Contraction Principle \cite[Theorem C.5.3]{DG15} implies that the natural transformation $i^! \rightarrow \piaut_!$ is an isomorphism of functors from $\G_m$-monodromic sheaves on $\BunX$ to sheaves on $\Bun_G$. Putting this into \eqref{eq: tate aut ext 1} gives the exact triangle
\begin{equation}\label{eq: tate aut ext 2}
\piaut_!  (\ul{k}_{\BunX} ) \rightarrow \piaut_* (\ul{k}_{\BunX} ) \rightarrow \piaut_* j_* j^* (\ul{k}_{\BunX}).
\end{equation}
Now, again by the Contraction Principle the natural transformation $\piaut_* \rightarrow i^*$ is an isomorphism of functors from $\G_m$-monodromic sheaves on $\BunX$ to sheaves on $\Bun_G$, giving in particular an isomorphism $\piaut_* (\ul{k}_{\BunX} ) \cong \ulk_{\Bun_G}$. Finally, observe that 
\[
\piaut_*\circ j_* \cong (\piaut \circ j)_*  = \AJ_*
\]
which transforms \eqref{eq: tate aut ext 2} into the exact triangle 
\begin{equation}\label{eq: tate aut ext 3}
\piaut_! (\ul{k}_{\BunX}) \rightarrow \ulk_{\Bun_G} \rightarrow \AJ_*(j^* \ul{k}_{\BunX}).
\end{equation}

Note that we have a natural isomorphism $\AJ_* \cong \AJ_![1]$ because $\AJ$ is the composition of a $\G_m$-torsor and a projective space bundle. Unraveling the construction, we have a commutative diagram 
\[
\begin{tikzcd}
 \ulk_{\Bun_G} \ar[r, "\eqref{eq: Tate aut side}"] \ar[d, equals] & \AJ_! (j^* \ulk_{\BunX})[1]	  \ar[d, "\sim"] \\
 \ulk_{\Bun_G} \ar[r, "\eqref{eq: tate aut ext 3}"] &  \AJ_*(j^* \ul{k}_{\BunX})
 \end{tikzcd}
\]
where the top row is the connecting homomorphism for \eqref{eq: Tate aut side}, which therefore lies in the space 
\begin{align}\label{eq: Tate aut ext}
\Hom_{\Bun_G}(\ul{k}_{\Bun_G} , \AJ_* (\ul{k}_{\Sym^{\bu} C})) & \cong  \Hom_{\Sym^{\bu} C}(\ul{k}_{\Sym^{\bu} C}, \ul{k}_{\Sym^{\bu} C}) \nonumber \\
& \cong \prod_{d \geq 0} \Hom_{\Sym^d C} (\ul{k}_{\Sym^d C}, \ul{k}_{\Sym^d C}) \cong \prod_{d \geq 0} k.
\end{align}

\begin{lemma}\label{lem: tate aut connecting}
For each $d \geq 0$, under our identifications the connecting map for \eqref{eq: Tate aut side} is adjoint to the identity map in $\Hom_{\Sym^d C} (\ul{k}_{\Sym^d C}, \ul{k}_{\Sym^d C})$ in \eqref{eq: Tate aut ext}. 
\end{lemma}

\begin{proof}
The map is given by a scalar for each $d$. To compute the scalar, we may base change to a point $\cL  \in \Bun_G(\F)$. Let 
\[
\begin{tikzcd}
V \ar[d, "\pi"] \ar[r, "l"]  & \BunX \ar[d, "\piaut"] \\
\cL \ar[r, "l"] & \Bun_G
\end{tikzcd}
\]
Note that $l^* \piaut_* \cong \pi_* l^*$ on $\G_m$-monodromic sheaves, by hyperbolic localization. Here $V$ is the derived vector bundle associated to $\Gamma(C; \cL)$. Write $i \co 0 \inj V$ for the origin and $j$ for its open complement. The map  $\piaut_* (\ulk_{\BunX}) \rightarrow \AJ_*(j^* \ul{k}_{\BunX})$ from \eqref{eq: tate aut ext 2} pulls back under $l^*$ to the map given by the unit of $j$ in the left side of
\[
\Hom(\pi_* \ulk_V, \pi_* j_* j^* \ulk_V)  \cong \Hom(j^* \pi^* \pi_* \ulk_V, j^* \ulk_V),
\]
which transports to the counit of $\pi$ on the right side, which is evidently the identity upon identifying $\pi^* \pi_* \ulk_V = \ulk_V$. The hyperbolic localization is implemented by the map $\pi_* \ulk_V \rightarrow \pi_* i_* i^* \ulk_V$, and under the identifications $\pi_* \ulk_V = \ulk_L$ and $i^* \ulk_V = \ulk_L$, it is the identity map. Hence we are done. 
\end{proof}



\subsection{Comparison} We will now complete the proof of Theorem \ref{thm: Tate duality}. We note that in this case $\IndCoh_{\Nilp}(\LocG) \cong \QCoh(\LocG)$ because $\chG$ is a torus. 

We first collect some normalizations. The Whittaker compatibility \eqref{eq: whittaker normalization} specializes in this case to the statement that for the trivial bundle $\triv \in \Bun_G$, we have 
\begin{equation}\label{eq: Tate whittaker normalization}
\Dmod(\Bun_1) \ni \delta_{\triv} \xmapsto{\LL_G} \omega_{\Loc_1}\tw{-(g-1)} \in \QCoh(\Loc_1).
\end{equation}
On the other hand, as $\dim \Bun_1 = (g-1)$, we have for $\triv \in \Loc_1$ the trivial local system, 
\begin{equation}\label{eq: Tate perverse normalization}
\Dmod(\Bun_1) \ni \ulk\tw{g-1} \xmapsto{\LL_G} \cO_{\triv} \in \QCoh(\Loc_1).
\end{equation}

\begin{lemma}\label{lem: Tate HkWhit}
We have 
\[
\Dmod(\Bun_1) \ni  \AJ_! (\ul{k}_{\Sym^{\bu} C})  \tw{g-1} \xmapsto{\LL_G}  \Loc^{\spec}(\Fact(\Sym^{\bu} \Std)) \otimes \omega_{\LocG} \in \QCoh(\Loc_1).
\]
\end{lemma}

\begin{proof}
By the special case of Hecke compatibility of $\LL_G$ in Example \ref{ex: Hecke compatibility}, $\LL_G$ intertwines the tensoring action of $ \Loc^{\spec}(\Fact(\Sym^{\bu} (\Std))) \in \QCoh(\LocG)$ with the Hecke convolution action of $\prod_{n \geq 0} \ulk_{\Hk_{G, C^{(n)}}^{\Std}} \tw{n(r-1)} \in \Shv(\bigcup_{n \geq 0} \Hk_{G, C^{(n)}})$ on $\Dmod(\Bun_G)$. In this case, note that $n(r-1) = 0$ because $r=1$, and $\Hk_{G, C^{(n)}}^{\Std}$ can be identified with the stack of tuples $(D, \cL_1, \cL_2, \cL_2 \cong \cL_1(D))$. Applying this statement to \eqref{eq: Tate whittaker normalization}, we deduce that $\LL_G$ takes 
\begin{equation}\label{eq: Tate HkWhit}
\Dmod(\Bun_G) \ni \left( \prod_{n \geq 0} \ulk_{\Hk_{G, C^{(n)}}^{\Std}}  \right)\star \delta_{\triv} \xmapsto{\LL_G}  \Loc^{\spec}(\Fact(\Sym^{\bu} (\Std))) \otimes \omega_{\LocG} \tw{-(g-1)}\in \QCoh(\LocG).
\end{equation}
Finally, unraveling the definition of the Hecke convolution reveals that 
\[
 \ulk_{\Hk_{G, C^{(n)}}^{\Std}} \star \delta_{\triv} \cong \AJ_!^n (\ul{k}_{C^{(n)}})
 \]
 so that the left side of \eqref{eq: Tate HkWhit} is identified with $\AJ_! (\ul{k}_{\Sym^{\bu} C})$. Inserting this into \eqref{eq: Tate HkWhit} and twisting by $\tw{g-1}$ completes the proof.
\end{proof}

Inspecting \eqref{eq: tate spectral side} and \eqref{eq: Tate aut side}, and using Proposition \ref{prop: Tate Loc}, we have a diagram where $\LL_G$ intertwines the indicated objects:
\begin{equation}
\begin{tikzcd}[column sep = tiny
]
\AJ_! (\ul{k}_{\Sym^{\bu} C}) \tw{g-1}   \ar[r] \ar[d, "\LL_G", "\text{Lemma \ref{lem: Tate HkWhit}}"', mapsto] &  \piaut_! (\ul{k}_{\BunX} ) \tw{g-1} \ar[d, "\LL_G?", dashed] \ar[r] &  \ul{k}_{\Bun_G} \tw{g-1}\ar[d, "\LL_G", "\eqref{eq: Tate perverse normalization}"', mapsto]  & \in \Dmod(\Bun_G)  \ar[d, "\LL_G"] \\ 
\Loc^{\spec}(\Fact(\Sym^{\bu} (\Std)))\otimes \omega_{\LocG} \ar[r] &   \pispec_* (\omega_{\LocX}) \ar[r] &   
\cO_{\mrm{triv}}  & \in \QCoh(\LocG)
\end{tikzcd}
\end{equation}
To verify that $\LL_G$ takes $\piaut_! (\ul{k}_{\BunX} )\tw{g-1}$ to $\pispec_* (\omega_{\LocX}) $, it therefore suffices to verify that $\LL_G$ takes the connecting map $\ul{k}_{\Bun_G} \tw{g-1} \rightarrow 
\AJ_! (\ul{k}_{\Sym^{\bu} C}) \tw{g-1}[1]$ to the connecting map $\cO
_{\mrm{triv}} \rightarrow \Loc^{\spec}(\Fact(\Sym^{\bu} (\Std))) \otimes \omega_{\LocG}[1]$. The functor $\LL_G$ is the identity on the extension groups in question with respect to the identifications \eqref{eq: Tate spectral ext} and \eqref{eq: Tate aut ext}, and then we conclude by using Lemma \ref{lem: tate spectral connecting} and Lemma \ref{lem: tate aut connecting} to identify the connecting maps on each side.\qed 







\subsection{Functional equation}\label{ssec: functional equation} We explain a categorification of the \emph{functional equation} for Tate's zeta functions, which compares the automorphic periods for the standard action of $\G_m$ on $\A^1$ and its inverse. To emphasize the distinction, we write $X := \Std$ for $\A^1$ with the standard scaling action of $G = \G_m$, and $X' := \Std^\vee$ for $\A^1$ with the inverse of the standard action. We equip $X$ and $X'$ with the standard scaling $\Ggr$-action (i.e., the neutral actions).

\subsubsection{Normalized period sheaf for $X$} Let $\cE := \ul{\RHom(\cO_C, \cL^{\univ} \otimes \Omega_C^{1/2})}$ be the perfect complex on $\Bun_1$, whose pullback to $\cL \in \Bun_1(R)$ is $\Gamma(C_R; \cL)$ viewed as an animated $R$-module. Let $\pi \co E \rightarrow \Bun_1$ be the associated derived vector bundle. We define the derived stack $\Bun_G^{X}$ as in \cite[(10.6)]{BZSV} but taking the \emph{derived} fibered product and forming mapping spaces as \emph{derived stacks}:
\[
\begin{tikzcd}
\Bun_G^{X} \ar[r] \ar[d] & \Map(C, \frac{X}{G \times \Ggr}) \ar[d] \\
\Bun_G \ar[r, "\Id \boxtimes \Omega_C^{1/2}"] & \Bun_{G \times \Ggr}
\end{tikzcd}
\]
By \cite[Proposition 5.34]{FYZ2}, $\Bun_G^{X}$ is canonically identified as a derived stack with $E$, in a manner intertwining the projection $\pi$ with $\piaut \co \Bun_G^{X} \rightarrow \Bun_G$. Hence the un-normalized period sheaf for $X$ is $\cP_X \cong \pi_! \ulk_E$. To normalize $\cP_X$, we note that $\eta_X \co G \rightarrow \G_m$ is the standard character (i.e., the identity). We pull back the degree function $\deg \co \Bun_{\G_m} \rightarrow \Z$ to $\Bun_G$ via $\eta_X$, so it assigns to $\cL \in \Bun_G(R)$ the integer $\deg \cL$. Then the \emph{normalized} period sheaf \cite[(10.12)]{BZSV} is 
\begin{equation}\label{eq: tate X normalized period}
\cP_X^{\norm} \cong \pi_! \ulk_E \tw{\deg + (g-1)} \in \Shv(\Bun_G). 
\end{equation}

\subsubsection{Normalized period sheaf for $X'$} Let $\cE' := \ul{\RHom(\cO_C, (\cL^{\univ})^\vee \otimes \Omega_C^{1/2})}$ be the perfect complex on $\Bun_1$, which over $\cL \in \Bun_1(R)$ is $\Gamma(C_R; \cL^\vee \otimes \Omega_C^{1/2})$. Let $\pi' \co E' \rightarrow \Bun_1$ be the associated derived vector bundle. Again, the derived stack $\Bun_G^{X'}$ is canonically identified with $E'$, in a manner intertwining the projection $\pi'$ with $\piaut \co \Bun_G^{X'} \rightarrow \Bun_G$. Hence the un-normalized period sheaf for $X$ is $\cP_{X'} \cong \pi'_! \ulk_{E'}$. Now, in this case $\eta_{X'} \co G \rightarrow \G_m$ is the \emph{inverse} of the standard character. Therefore, the degree function $\deg \co \Bun_{\G_m} \rightarrow \Z$ pulls back to $-\deg \co \Bun_G \rightarrow \Z$ under $\eta_{X'}$. Hence the \emph{normalized} period sheaf is 
\begin{equation}\label{eq: tate X' normalized period}
\cP_{X'}^{\norm} \cong \pi'_! \ulk_{E'} \tw{-\deg + (g-1)} \in \Shv(\Bun_G).
\end{equation}

\subsubsection{Categorification of the functional equation} The following result is the categorification of the functional equation for the $L$-functions attached to Hecke characters of $\GL_1$. 

\begin{prop}\label{prop: tate functional equation}
With the notation above, there is a canonical isomorphism 
\[
\cP_X^{\norm} \cong  \cP_{X'}^{\norm} \in \Shv(\Bun_1).
\]
\end{prop}

\begin{proof} Let $r_E :=  \rank(E)$ be the virtual rank of $E$ over $\Bun_1$. Note that $r_E = \deg$ by Riemann-Roch. Regarding $\Bun_1$ as the total space of the zero vector bundle over itself, we have from \S \ref{sssec: FT functoriality} natural isomorphisms 
\begin{equation}\label{eq: FE 1}
\pi_{!} \ulk_{E} = \FT_0(\pi_{!} \ulk_{E}) \cong z_{E^\vee}^* \FT_E(\ulk_E)[-r_E] \cong z_{E^\vee}^*  (\delta_{E^\vee}) \tw{-2r_E}
\end{equation}
where $z_{E^\vee}$ is the zero-section of $E^\vee$. Note that Serre duality provides an isomorphism of perfect complexes over $R$,
\[
\Gamma(C_R; \cL \otimes \Omega_C^{1/2})^\vee \cong \Gamma(C_R; \cL^\vee \otimes \Omega_C^{1/2}[1])
\]
so that $E$ is dual to $E'[1]$. Hence we have a derived Cartesian square
\begin{equation}\label{diag: FE}
\begin{tikzcd}
E' \ar[r] \ar[d, "\pi'"] \ar[r, "\pi'"] & \Bun_1  \ar[d, "z_{E^\vee}"] \\ 
\Bun_1 \ar[r, "z_{E^\vee}"] & E^\vee
\end{tikzcd}
\end{equation}
Proper base change applied to \eqref{diag: FE} gives a natural isomorphism 
\begin{equation}\label{eq: FE 2}
z_{E^\vee}^*  (\delta_{E^\vee}) = z_{E^\vee}^*  z_{E^\vee!} (\ulk_{\Bun_1}) \cong \pi'_! (\pi')^* (\ulk_{\Bun_1} )\cong \pi'_! (\ulk_{E'}) \in \Shv(\Bun_1),
\end{equation}
and the result follows from combining \eqref{eq: FE 1} and \eqref{eq: FE 2}, and then comparing them to \eqref{eq: tate X normalized period} and \eqref{eq: tate X' normalized period}. 
\end{proof}

\begin{remark}
On quasicompact open substacks of $\Bun_1$, an argument for Proposition \ref{prop: tate functional equation} without using derived Fourier analysis appears in \cite[Lemma 7.3.6]{BG02} and is credited there to Drinfeld. Drinfeld's argument is based on an ingenious construction, which we invite the reader the admire. Nevertheless, we feel that the perspective of the derived Fourier transform puts the statement in its proper general context. In fact, trying to understand the functional equation is what led the authors to the derived Fourier transform, the details of which were then written up in \cite{FYZ3}.
\end{remark}

\section{Hecke period}\label{sec: hecke}

\subsection{Setup}
We study an anzsatz that we call the ``Hecke period'', as it corresponds to Hecke's construction \cite{Hecke37} of $L$-functions for modular forms (i.e., automorphic forms on $\GL_2$) by Mellin transform. 
 
\subsubsection{Automorphic side:} Let $A := \begin{pmatrix} 1 \\ & * \end{pmatrix} \cong \G_m  \inj \GL_2$ and $X = \GL_2/A$, with $G = \GL_2$ acting by left translation. 

\subsubsection{Spectral side:} Let $\chX = \A^2$ with $\chG = \GL_2$ acting by the standard representation. 

\subsubsection{$\Ggr$-action}According to \cite[Example 12.6.7]{BZSV}, the neutral $\Ggr$-action on $X$ is trivial, while the neutral $\Ggr$-action on $\chX$ is the standard scaling action. The un-normalized $\Ggr$-action on $X$ is via the \emph{inverse} of the standard character into the center $\G_m \cong Z(\GL_2)$, while the un-normalized $\Ggr$-action on $\chX$ is still the standard scaling action. Below we will take the un-normalized actions of $\Ggr$.

\subsubsection{Automorphic period sheaf} Examining the recipe of \cite[\S 10.2]{BZSV}, we see that $\BunX$ can be identified with $\Bun_1$ in such a way that the map $\piaut\co \BunX \rightarrow \Bun_G$ is identified with the map sending
\[
\cL \mapsto (\cL \otimes \Omega_C^{1/2}) \oplus \Omega_C^{1/2}
\]
where $\Omega_C^{1/2}$ is the fixed spin structure (\S \ref{ssec: spin}). The automorphic period sheaf is
\[
\PX := \piaut_! (\ulk_{\BunX}) \in \Shv(\Bun_G).
\]
In the language of \cite{BZSV}, this is the un-normalized period sheaf for the un-normalized action of $\Ggr$. 

\subsubsection{Spectral period sheaf}By definition, $\LocX$ is the derived mapping stack
\[
\LocX = \Map(C_{\dR}, \A^2/\GL_2)
\]
with $R$-points the groupoid of pairs $(E,s)$, where: 
\begin{itemize}
\item $E \in \Loc_2(R)$ is a rank 2 vector bundle on $C_R$ with a flat connection along $C$,
\item $s \in \hSect(C_{\dR}, E)(R)$, or informally speaking an $R$-family of flat sections of $E$. 
\end{itemize}
The map $\pispec_*  \co \LocX \rightarrow \LocG$ forgets the datum of $s$. The (co-localized) spectral period sheaf is 
\[
\LX \cong (\pispec_* \omega_{\LocX})^{\shear} :=  \bigoplus_n (\pispec_* \omega_{\LocX})_n \tw{n} \in \IndCoh_{\Nilp}(\LocG)
\]
where $(\pispec_* \omega_{\LocX})_n$ is the $n$th graded piece for the $\Ggr$-action. In the language of \cite{BZSV}, this is the un-normalized $\cL$-sheaf for the un-normalized action of $\Ggr$. 

\subsubsection{Duality} In this case, the projection of \cite[Conjecture 12.1.1]{BZSV} to the formulation of Geometric Langlands in \cite{AG15} is equivalent to 
\[
\Dmod(\Bun_G) \ni \PX \tw{5(g-1)} \xmapsto{\LL_G} \LX \in \IndCoh_{\Nilp}(\LocG).
\]
Indeed, \cite[Proposition 12.6.4]{BZSV} says that \cite[Conjecture 12.1.1]{BZSV} is equivalent to $\PX \tw{r} \xmapsto{\LL_G} \LX^d$ where 
\[
r = \beta_{X'} + \beta_{\chX'} - (g-1) \tau. 
\]
Here the $(-)'$ refers to un-normalized $\Ggr$-action. The duality involution disappears under the translation from right actions to left actions (cf. \cite[Appendix A]{CV}). Here $\tau = \langle \eta, \check{\eta} \rangle = 0$ where $\eta$ is the character by which $G$ scales the volume form on $X$; it vanishes because $\GL_2$ and $A$ are unimodular, so $\GL_2/A$ admits a left $\GL_2$-invariant volume form. Then $\beta_{X'} = (g-1)(\dim G - \dim X'  + \gamma_{X'}) $, but $\gamma_{X'} = 0$ since the action of $\G_m$ is through $G$. Also $\beta_{\chX'}  = (g-1)(\dim \chG - \dim \chX'  +  \gamma_{\chX'})$ and $\gamma_{\chX'} = 2$ since the $\Ggr$-action scales the volume form by the square of the standard character. All in all, we see that 
\[
r = (g-1) (4-3) + (g-1) (4-2 + 2) = 5(g-1).
\]

\subsection{Spectral side} We let $Z \inj \LocX$ be the zero section and $U$ be the open complement of $Z$. We consider the exact triangle in $\IndCoh_{\Nilp}(\LocG)$ coming from filtering by sections with set-theoretic support in $Z$: 
\begin{equation}\label{eq: Hecke spectral triangle}
(\pispec)_{Z*}^{\wedge} (\omega_{\LocX}) \rightarrow \pispec_* (\omega_{\LocX}) \rightarrow \pispec_*|_U ( \omega_{\LocX}|_U)
\end{equation}
where we have implicitly co-localized to the category of sheaves with nilpotent singular support (\S \ref{ssec: singular support}).

\subsubsection{Open stratum}\label{ssec: Hecke spectral open} We analyze the rightmost term of \eqref{eq: Hecke spectral triangle}. Since a flat section of a family of local systems is non-zero at one point if and only if it is non-zero at every point, $U$ is isomorphic to the derived mapping stack $\Map(C_{\dR}, (\A^2 \setminus 0)/\GL_2)$ over $\LocG$. The action of $\GL_2$ on $\begin{pmatrix} 1 \\ 0 \end{pmatrix}$ identifies $\A^2 \setminus 0$ with $\GL_2/\Mir_2$, where 
\[
\Mir_2 := \left\{\begin{pmatrix} 1 & * \\ 0 & * \end{pmatrix}\right\} \inj \GL_2
\]
is the mirabolic subgroup. This induces a presentation of $U$ as the derived mapping stack 
\[
U \cong  \Map(C_{\dR}, \rB\Mir_2)
\]
such that the map $\pispec|_U \co U \rightarrow \Loc_{\chG}$ is the one induced by the inclusion $\Mir_2 \inj \GL_2$. 

Therefore, with respect to the diagram
\[
\begin{tikzcd}
& \Loc_{\Mir_2} \ar[dl, "q_{\Mir}"'] \ar[dr, "p_{\Mir}"] \\
\Loc_{1}  & &  \Loc_{2}
\end{tikzcd}
\]
we have  
\[
\pispec_*|_U ( \omega_U)  \cong  p_{\Mir *} q_{\Mir}^!(\omega_{\Loc_1}) \in \IndCoh(\Loc_2).
\]

In order to fit the formulation of Eisenstein compatibility in \S\ref{ssec: Eisenstein compatibility}, we rewrite this in terms of usual spectral Eisenstein functors (for parabolic subgroups). Let $\chB \subset \chG$ be the standard Borel subgroup containing $\Mir_2$, and $T$ its Levi quotient. We consider the correspondence diagram 
\[
\begin{tikzcd}
&  \Loc_{\check{B}} \ar[dl, "q"'] \ar[dr, "p"]  \\
\Loc_{\check{T}} &  & \LocG
\end{tikzcd}
\]
We then have a commutative diagram 
\[
\begin{tikzcd}
& \Loc_{\Mir_2} \ar[dl, "q_{\Mir}"'] \ar[d, "f"] \ar[ddr, "p_{\Mir}", bend left] \\
\Loc_{1} \ar[d, "g"] &  \Loc_{\check{B}} \ar[dl, "q"] \ar[dr, "p"']  \\
\Loc_{\chT} &  & \LocG
\end{tikzcd}
\]
where the upper left square is Cartesian. By base change for that Cartesian square, we get isomorphisms
\begin{equation}\label{eq: Hecke open 1}
p_{\Mir *} (\omega_U ) \cong p_{\Mir*}q_{\Mir}^! (\omega_{\Loc_1}) \cong p_* q^! g_* (\omega_{\Loc_1}) \in \IndCoh(\LocG).
\end{equation}
Now, $g$ can be identified with the map 
\[
\pt \times \Loc_{1} \xrightarrow{\mrm{triv}\times \Id} \Loc_{1} \times \Loc_{1}
\]
so that
\[
g_* (\omega_{\Loc_1}) \cong  \cO_{\triv} \boxtimes \omega_{\Loc_1} \in \IndCoh(\Loc_1 \times \Loc_1).
\]
Inserting this into \eqref{eq: Hecke open 1}, and co-localizing to the category with nilpotent singular support, we have established: 

\begin{prop} We have an isomorphism (using the notational convention of \S \ref{ssec: singular support})
\[
\pi_{\spec}|_{U*} ( \omega_U) \cong \Eis^{\spec}_{\chB}( \cO_{\triv} \boxtimes \omega_{\Loc_{1}}) \in \IndCoh_{\Nilp}(\LocG). 
\]
\end{prop}

\subsubsection{Closed stratum}\label{ssec: Hecke spectral closed}  Next we analyze the leftmost term in \eqref{eq: Hecke spectral triangle}, with a similar approach to that of \S \ref{sssec: spectral closed tate}. The origin $0 \inj \Std = \chX$ is the only $\chG$-fixed point. Let $\cV$ be its formal completion in $\chX$. By Lemma \ref{lemma: omega closed support}, we have an isomorphism of unital cocommutative coalgebras
\begin{equation}\label{eq: Hecke spectral 0 nbhd}
\Gamma_0(\chX; \omega_{\chX}) \cong \Gamma(\cV; \omega_{\cV}) \in \Rep(\chG \times \Ggr).
\end{equation}
The $\Ggr$-action equips $\Gamma_0(\chX; \omega_{\chX})$ with a non-negative grading. 

\begin{lemma}\label{lem: Hecke spectral factorization algebra}
As a commutative coalgebra in $\Rep(\chG)^{\Z_{>0}}$, $\Gamma_0(\chX; \omega_{\chX})$ identifies with $\Sym^{\bu}(\Std)$ where $\Std$ is the standard representation of $\chG$ in graded degree $1$. 
\end{lemma}

\begin{proof}
This follows from Lemma \ref{lem: formal thickening of dvb}.
\end{proof}

This gives the following description of the leftmost term in $\eqref{eq: Hecke spectral triangle}$. 

\begin{prop}\label{prop: Hecke Loc}
Let $\Fact(\Sym^{\bu}(\Std))$ be the graded unital cocommutative factorization coalgebra \eqref{eq: augment spectral fact algebra}. Then we have an isomorphism
\[
(\pispec)^{\wedge}_{Z*} (\omega_{\LocX}) \cong \Loc^{\spec}(\Fact(\Sym^{\bu} \Std)) \otimes \omega_{\LocG}\in \QCoh(\LocG).
\]
\end{prop}

\begin{proof} Applying Proposition \ref{prop: distribution algebra} to the origin inside $\Std$ viewed as a $\chG$-scheme, equipped with the standard scaling action of $\Ggr$, we obtain an isomorphism
\begin{equation}\label{eq: Hecke non-unital Loc 1} 
\ol{(\pispec)_{Z*}^{\wedge} (\omega_{\LocX})} \cong \pr_{2*} (\Fact(\ol{\Gamma}_0(\chX; \omega_{\chX}))^{\univ} ) \otimes \omega_{\LocG}\in \QCoh(\LocG).
\end{equation}
In \eqref{eq: Hecke spectral 0 nbhd} and Lemma \ref{lem: Hecke spectral factorization algebra} we saw that $\ol{\Gamma}_0(\chX; \omega_{\chX}) \cong \ol{\Sym}^{\bu} \Std \in \Rep(\chG)^{\Z_{>0}}$. Inserting this into \eqref{eq: Hecke non-unital Loc 1} yields an isomorphism
\begin{equation}\label{eq: Hecke non-unital Loc 2} 
\ol{(\pispec)_{Z*}^{\wedge} (\omega_{\LocX})} \cong \pr_{2*} (\Fact(\ol{\Sym}^{\bu} \Std)^{\univ}) \otimes \omega_{\LocG} \in \QCoh(\LocG)
\end{equation}
where the meaning of $\Fact(\ol{\Sym}^{\bu} \Std)^{\univ}$ is as in Example \ref{ex: Hecke compatibility}. The rest follows from adding a copy of the unit, as in the proof of Proposition \ref{prop: Tate Loc}. 

\end{proof}


\subsubsection{Summary}
Putting together \S \ref{ssec: Hecke spectral open} and \S \ref{ssec: Hecke spectral open}, we have produced an exact triangle in $\IndCoh_{\Nilp}(\LocG)$: 
\begin{equation}\label{eq: hecke triangle 1}
\Loc^{\spec}(\Fact(\Sym^{\bu} \Std)) \otimes \omega_{\LocG} \rightarrow  \pispec_* (\omega_{\LocX}) \rightarrow  
\Eis^{\spec}_{\check{B}}( \cO_{\triv} \boxtimes \omega_{\Loc_1})
\end{equation}
where we are invoking the notational convention of \S \ref{ssec: singular support}. 

\subsection{Automorphic side}\label{ssec: Hecke aut side} Twisting by $\Omega_C$, we may identify $\piaut \co \BunX \rightarrow \Bun_2$ with the map $\Bun_1 \rightarrow \Bun_2$ sending $\cL \mapsto (\cL \otimes \Omega_C^{-1/2} )\oplus \Omega_C^{1/2}$. We will write $\piaut \co \Bun_1 \rightarrow \Bun_2$ when using the latter perspective. 
 
 We seek a filtration of $\cP_X$ corresponding to the filtration of $\cL_{\chX}$ in \eqref{eq: RS triangle 1}, which came from the stratification of $\chX$ by $\chG$-orbits. But note that unlike in the example of the Tate period, here there is no obvious ``open-closed'' decomposition of $\piaut_! (\ul{k}_{\Bun_1})$ coming from a stratification of $X$ to match the stratification of $\chX$, since the $G$-action on $X$ is homogeneous. (However, see \S \ref{ssec: microlocal} for matching stratifications in the hyperspherical varieties.) We will instead look for a decomposition in the \emph{Fourier dual} space.

\subsubsection{Unfolding}\label{sssec: unfolding}
On $\Bun_1$ there is a perfect complex $\ul{\RHom(\cL^{\univ}, \Omega_C[1])}$ whose pullback to $\cL \in \Bun_1(R)$, for any animated $\F$-algebra $R$, is $\RHom_{C_R}(\cL,  \Omega_C \boxtimes R[1])$ regarded as an animated $R$-module. We denote its total space by $E \rightarrow \Bun_1$. By Serre duality, the dual derived vector bundle $E^\vee \rightarrow \Bun_1$ can be identified with the total space of the perfect complex $\ul{\RHom(\cO_C, \cL^{\univ})}$, whose pullback to $\cL \in \Bun_1(R)$ is $\RHom_{C_R}(\cO \boxtimes R, \cL)$ regarded as an animated $R$-module. 
 
In fact, $E$ is a classical stack, with $R$-points the groupoid of extensions
\begin{equation}\label{eq: E points}
\Omega_C \boxtimes R \rightarrow  \cF \rightarrow \cL
\end{equation}
where $\cL \in \Bun_1(R)$. The zero-section $\Bun_1 \rightarrow E $ sends $\cL \in \Bun_1(R)$ to the split extension $\cL \oplus (\Omega_C \boxtimes R)$. 
 
 Let $\deg$ be the locally constant function on $\Bun_1$ assigning to a family of line bundles the fiberwise degree. The virtual rank of $E$ as a derived vector bundle over $\Bun_1$ is 
 \[
 r_E :=  \chi(\ul{\RHom(\cL^{\univ}, \Omega_C[1])}) = 1-g+\deg.
 \]
Hence by Example \ref{ex: FT delta} we have an isomorphism 
\begin{equation}\label{eq: Hecke fourier} 
\delta_{E} \cong \FT_{E^\vee}(\ul{k}_{E^\vee} [r_E](r_E)) \in \Shv(E). 
\end{equation}
We will then analyze the decomposition of $\delta_E$ coming from the stratification of $E^\vee$ into its zero section and the complement. To unravel what this means, consider the commutative diagram
 \begin{equation}\label{eq: Hecke aut commutative diagram}
 \begin{tikzcd}
& &  E \times_{\Bun_1} E^{\vee} \ar[dl, "\pr_1"'] \ar[dr, "\pr_2"] \\
 & E \ar[dr] \ar[dl, "f"'] & & E^\vee \ar[dl] \\
 \Bun_2 & & \Bun_1
 \end{tikzcd}
 \end{equation}
 where $f$ sends \eqref{eq: E points} to $\cF \otimes \Omega_C^{-1/2} \in \Bun_2(R)$, and the middle square is derived Cartesian. The classical truncation $(E \times_{\Bun_1} E^\vee)_{\cl}$ has $R$-points the groupoid of diagrams
\begin{equation}\label{eq: hecke period diagram}
\left\{ \begin{tikzcd}
&   \Omega_C \boxtimes R \ar[d]  \\
& \Cal{F} \ar[d] \\
\Cal{O}_{C} \boxtimes R \ar[r, "s"] & \Cal{L}
\end{tikzcd} \right\}
\end{equation}
where 
\begin{itemize}
\item $\cF$ is a rank $2$ vector bundle on $C_R$,
\item $\cL$ is a line bundle on $C_R$, and 
\item $s \co \cO_{C} \boxtimes R \rightarrow \cL$ is any section (possibly zero). 
\end{itemize}
The map 
\[
\ev \co E \times_{\Bun_1} E^\vee \rightarrow \A^1 
\]
assigns to a diagram \eqref{eq: hecke period diagram} the extension class of the pullback of $\cF$ along $s$, which is an element of $\rH^1(C_R, \Omega^1 \boxtimes R) \cong R$. 

\begin{lemma}\label{lem: hecke period fourier}
Recall that $\AS$ is the Artin-Schreier sheaf on $\A^1$. With the notation above, we have an isomorphism 
\[
\piaut_! (\ul{k}_{\Bun_1}) = f_! \delta_E  \cong f_! \pr_{1!} (\ev^* \AS \tw{2r_E}) \in \Shv(\Bun_2).
\]
\end{lemma}

\begin{proof}This follows from \eqref{eq: Hecke fourier} and the commutative diagram \eqref{eq: Hecke aut commutative diagram} by writing out the definition of $\FT$. 
\end{proof}

Consider the derived fibered product
\[
\begin{tikzcd}
 Z \ar[r, hook, "i"] \ar[d] &  E \times_{\Bun_1} E^\vee \ar[d, "\pr_1"] \\ 
  \Bun_1 \ar[r, "0"] & E 
 \end{tikzcd}
\]
Informally speaking, $i \co Z \inj E \times_{\Bun_1} E^\vee$ is the closed substack where the section $s$ is zero. Let $j \co U \rightarrow E \times_{\Bun_1} E^\vee$ be its open complement. Then we have an exact triangle in $\Shv(E \times_{\Bun_1} E^\vee)$, 
\[
j_! j^* (\ev^* \AS) \rightarrow \ev^* \AS \rightarrow i_* i^* (\ev^* \AS). 
\]
Then applying $f_! \pr_{1!} $ and using Lemma \ref{lem: hecke period fourier} gives an exact triangle in $\Shv(\Bun_2)$: 
\begin{equation}\label{eq: Hecke aut exact triangle}
f_!  \pr_{1!}j_! j^*  (\ev^* \AS) \tw{2r_E}  \rightarrow  \piaut_! (\ul{k}_{\BunX})  \rightarrow f_!   \pr_{1!} i_* i^*  (\ev^* \AS)  \tw{2r_E}.
\end{equation}
Henceforth we will have no further need to refer to the derived structure on the objects involved, so we pass to classical truncations.

\subsubsection{Closed stratum}\label{ssec: hecke aut closed} We analyze the rightmost term of \eqref{eq: Hecke aut exact triangle}. 

The closed stratum $Z \inj E \times_{\Bun_1} E$ may be identified with the (classical) stack $\Bun_{\Mir_2}^{\Omega}$, whose $R$-points form the groupoid of 
\begin{itemize}
\item $\cL$, a line bundle on $C \times R$.
\item An extension 
\begin{equation}\label{eq: Bun_Mir point}
\Omega_C \boxtimes R \rightarrow \cF \rightarrow \cL .
\end{equation}
\end{itemize} The stack $\Bun_{\Mir_2}^{\Omega}$ is a twisted version of the stack of torsors for the mirabolic group $\Mir_2 \inj \GL_2$.

Since $\ev|_Z$ is identified with projection to $0 \times \Bun_1 \inj \A^1 \times \Bun_1$, the $*$-restriction of $\ev^*\AS$ to $Z$ is identified with the constant sheaf $\ul{k}_{\Bun_{\Mir_2}^{\Omega}}$. Let $p_{\Mir}  = f \circ i \co \Bun_{\Mir_2}^\Omega \rightarrow \Bun_2$ be the map sending \eqref{eq: Bun_Mir point} to $\cF \otimes \Omega_C^{-1/2}$. So the rightmost term in \eqref{eq: Hecke aut exact triangle} is isomorphic to $p_{\Mir!}  (\ul{k}_{\Bun_{\Mir_2}^{\Omega}}) \tw{2 r_E}$. We rewrite this in terms of the Eisenstein functor. Let $B$ be the standard parabolic of $\GL_2$ containing $\Mir_2$. Then we have a commutative diagram 
\[
\begin{tikzcd}
& \Bun_{\Mir_2}^{\Omega} \ar[d, "g'"] \ar[dl, "q_{\Mir}"'] \ar[ddr, "p_{\Mir}"]  \\
\Bun_1  \ar[d, "g"] &  \Bun_{B} \ar[dl, "q"] \ar[dr, "p"'] \\
\Bun_1 \times \Bun_1 & & \Bun_2 
\end{tikzcd}
\]
in which the left square is Cartesian, and $g$ is the map $\cL \mapsto (\cL \otimes \Omega_C^{-1/2},  \Omega_C^{1/2})$. By proper base change, we have 
\begin{align*}
p_{\Mir !} (\ul{k}_{\Bun_{\Mir_2}^{\Omega}}) & \cong p_{\Mir !} q_{\Mir}^*(\ul{k}_{\Bun_1})   \cong p_! q^* g_! (\ul{k}_{\Bun_1})  =  \Eis_B (\ul{k}_{\Bun_1} \boxtimes \delta_{\mrm{spin}}) \in \Shv(\Bun_G)
\end{align*}
where $\delta_{\mrm{spin}}$ is the delta sheaf at $\Omega_C^{1/2} \in \Bun_1$, meaning $\mrm{spin}_{!} (\ul{k}_{\pt})$ for the map $\mrm{spin} \co \pt \rightarrow \Bun_1$ corresponding to the chosen spin structure $\Omega_C^{1/2}$. In conclusion, we have identified the term $f_! i_* i^* (\ev^* \AS)$ from \eqref{eq: Hecke aut exact triangle} with 
\[
p_{\Mir !} (\ul{k}_{\Bun_{\Mir_2}^{\Omega}}) \cong \Eis_B  (\ul{k}_{\Bun_1} \boxtimes \delta_{\mrm{spin}}) \in \Shv(\Bun_G). 
\]

\subsubsection{Drinfeld's compactifications} \label{sssec: drinfeld compactification}
In preparation for analyzing the open stratum, we take a digression to recall Drinfeld's compactifications. Let $G$ be a reductive group. Choose a Borel subgroup $B \subset G$ and let $N \subset B$ be the unipotent radical. Let $T = B/N$ and $\cF_T \in \Bun_T(\F)$.  

The map $u \co \Bun_N^{\cF_T} \rightarrow \Bun_G$ is not proper, even after quotienting out by the $T$-action on the source. Drinfeld constructed a stack $\ol{\Bun}_N^{\cF_T} $ fitting into a commutative diagram  
\[
\begin{tikzcd}
\Bun_N^{\cF_T} \ar[r, hook] \ar[dr, "u"'] & \ol{\Bun}_N^{\cF_T} \ar[d, "\ol{u}"] \\
& \Bun_G 
\end{tikzcd}
\]
such that the induced map $\ol{\Bun}_N^{\cF_T}/T \rightarrow \Bun_G$ is proper; this is explained in \cite[\S 2.2]{FGV01}. To explain its functor of points, we recall that the groupoid $\Bun_N^{\cF_T}(R)$ admits a description as classifying pairs $(\cF_G, \kappa)$ where $\cF_G$ is a $G$-torsor on $C_R$ and $\kappa$ is a system of (saturated) vector bundle embeddings
\[
\kappa^{\chlambda} \co \cF_T^{\chlambda} \inj \cV_{\cF_G}^{\chlambda}, \text{ for all } \chlambda \in X^*(T),
\]
satisfying the Pl\"{u}cker relations, where $\cF_T^{\chlambda} $ is the line bundle on $C_R$ obtained from $\cF_T$ via the homomorphism $\chlambda  \co T\rightarrow \G_m$; and $\cV_{\cF_G}^{\chlambda}$ is the vector bundle on $C_R$ obtained by twisting irreducible representation of $G$ with highest weight $\chlambda$ by $\cF_G$. Then $\ol{\Bun}_N^{\cF_T}(R)$ parametrizes tuples $(\cF_G, \kappa)$ defined similarly, except that the maps $\kappa^{\chlambda}$ are only required to be (not necessarily saturated) embeddings of coherent sheaves fiberwise over $\Spec R$. The maps $u$ and $\ol{u}$ are both described by $(\cF_G, \kappa) \mapsto \cF_G$. 

\begin{example}\label{ex: drinfeld compactification GL_2}
We explain the picture for $G = \GL_2$, which is the only case we will use in this paper. In this case, $\cF_T$ corresponds to a pair of line bundles $\cL_1, \cL_2$ on $C$. The stack $\Bun_N^{\cF_T}$ has $R$-points the groupoid of extensions
\[
\begin{tikzcd}
\cL_1 \boxtimes R \ar[r, "\kappa_1"] & \cF \ar[r, "\kappa_2"] &  \cL_2\boxtimes R 
\end{tikzcd}
\]
of vector bundles on $C_R$. The stack $\ol{\Bun}_N^{\cF_T}$ has $R$-points the groupoid of triples $(\cF, \kappa_1, \kappa_2)$ where 
\begin{itemize}
\item $\cF$ is a rank $2$ vector bundle on $C_R$,
\item $\kappa_1 \co \cL_1 \boxtimes R \inj \cF$ is an injection of coherent 	sheaves fiberwise on $\Spec R$, 
\item $\kappa_2 \co \cF \rightarrow \cL_2 \boxtimes R$ is a map of coherent sheaves which is non-zero fiberwise on $\Spec R$, and
\item the composition $\kappa_2 \circ \kappa_1$ is the zero map. 
\end{itemize}
\end{example}

Let $\Lambda := X_*(T)_{\geq 0}$. There is a stratification of $\ol{\Bun}_N^{\cF_T}$ indexed by $\Lambda$-colored divisors $D = \sum_i \nu_i \cdot x_i \in \Div^{\Lambda} C$; recall this means that $\nu_i \geq 0$ for each $i$.\footnote{We have reversed the conventions on positivity/negative from \cite{FGV01} for consistency with our discussion of graded configuration space.} The stratum $\ol{\Bun}_N^{\cF_T}(D)$ indexed by $D$ is characterized by the property that an $\F$-point $(\cF_G, \kappa) \in \ol{\Bun}_N^{\cF_T}$ lies in $\ol{\Bun}_N^{\cF_T}(D)$ precisely when 
\[
\text{the saturation of $\kappa^{\chlambda} \co \cF_T^{\chlambda}  \rightarrow \cV^{\chlambda}_{\cF_G}$ is }\cF_T^{\chlambda}\left( \sum_i \langle  \nu_i, \chlambda \rangle \cdot x_i \right) \text{ for each }\chlambda \in X^*(T).
\]
In particular, the stratum indexed by $D = 0$ is the open substack $\Bun_N^{\cF_T} \inj \ol{\Bun}_N^{\cF_T}$. 

\begin{example}
We explain the picture for $G = \GL_2$ in terms of Example \ref{ex: drinfeld compactification GL_2}. In this case, the stratification is indexed by a pair of effective divisors $D_1,D_2$ on $C$. The stratum corresponding to $(D_1, D_2)$ has $\F$-points the diagrams 
\[
\cL_1 \xrightarrow{\kappa_1} \cF \xrightarrow{\kappa_2} \cL_2
\]
such that the saturation of $\kappa_1$ is $\cL_1(D_1) \subset \cF$, and the image of $\kappa_2$ is $\cL_2(-D_2) \subset \cL_2$. 
\end{example}

 \subsubsection{Open stratum}\label{ssec: hecke aut open} We will analyze the leftmost term in \eqref{eq: Hecke aut exact triangle}.
 
 For $d \geq 0$, we abbreviate $\Hk_{2, C^{(d)}}^{\Std}$ for the Hecke stack $\Hk_{\GL_2, C^{(d)}}^{\Std}$ from \eqref{eq: Hecke Std}, and 
\[
\Hk_{2, \Sym^{\bu} C}^{\Std}  := \bigcup_{d \geq 0} \Hk_{2, C^{(d)}}^{\Std}
\]
for the stack $\Hk_{\GL_2, \Sym^{\bu} C}^{\Std}$ from \eqref{eq: Hecke Std}. We have a correspondence diagram
\[
\begin{tikzcd}
& \Hk_{2, \Sym^{\bu} C }^{\Std} \ar[dl, "h_l"'] \ar[dr, "h_r \times \pr"'] \ar[drr, "h_r"] \\
\Bun_2 & & \Bun_2 \times  \Sym^{\bu} C \ar[r] & \Bun_2
\end{tikzcd}
\]
where for a modification $(\cF \inj \cF') \in \Hk_{2, \Sym^{\bu} C}^{\Std}(R)$, we define 
\[
h_l(\cF \inj \cF') := \cF \in \Bun_2(R)
\]
\[
h_r(\cF \inj \cF') := \cF' \in \Bun_2(R)
\]
and $\pr(\cF \inj \cF')$ is the support divisor of $(\cF'/\cF)$. We write 
\begin{equation}\label{eq: Hecke hk operator}
\Hk_{2,  \Sym^{\bu}C}^{\Std} \star (-)  := \sum_{d \geq 0} \Hk_{2, C^{(d)}}^{\Std}  \star (-)
\end{equation}
for the endofunctor of $\Shv(\Bun_G)$ induced by the convolving with the kernel sheaf 
\begin{equation}\label{eq: Hecke Hk kernel}
\prod_{d \geq 0} \ulk\tw{d} \in \Shv(\bigcup_{d \geq 0} \Hk_{2, C^{(d)}}^{\Std})
\end{equation}
from \eqref{eq: Hk Std kernel sheaf}.

\begin{remark}[Satake sheaves]\label{rem: satake of sym std} In principle one should use the IC complexes coming from the Geometric Satake equivalence as kernel sheaves for Hecke convolution. Recall from Example \ref{ex: Hecke compatibility} that the definition \eqref{eq: Hecke hk operator} is justified by the fact that the Satake sheaf associated to $\Sym^d(\Std) \in \Rep(\GL_2)$ is the (shifted and twisted) \emph{constant sheaf} $\ulk\tw{d}$ on the corresponding Schubert variety. 
\end{remark}



Recall that $\sW \in \Shv(\Bun_2)$ is the Whittaker sheaf (cf. \S \ref{sssec: Whittaker compatibility}). We write 
$\tw{\deg} \co \Shv(\Bun_2) \rightarrow \Shv(\Bun_2)$ for the endofunctor given by $\tw{d}$ on 
$\Shv(\Bun_2^d)$, where $\Bun_2^d$ is the connected component parametrizing bundles with degree $d$. 

\begin{prop}\label{prop: Hecke-Whittaker}
With notation as in \S \ref{sssec: unfolding}, we have an isomorphism 
\[
\Hk_{2, \Sym^{\bu}(C)}^{\Std} \star \sW \cong f_! \pr_{1!} j_! j^*  (\ev^* \AS ) \tw{\deg}  \in \Shv(\Bun_G).
\]
\end{prop}

\begin{proof}
Let $\cF_T := \Omega_C^{1/2} \oplus \Omega_C^{-1/2}$. Then, with notation as in \S \ref{sssec: drinfeld compactification}, $\Bun_N^{\cF_T}$ classifies families of extensions
\[
0 \rightarrow \Omega_C^{1/2} \boxtimes R  \rightarrow \cE \rightarrow \Omega_C^{-1/2} \boxtimes R \rightarrow 0.
\]
Recall that the Whittaker sheaf $\sW$ is obtained by pulling back the Artin-Schreier sheaf $\AS$ via the map $\Bun_N^{\cF_T} \rightarrow \A^1$, obtaining $\Psi \in \Shv(\Bun_N^{\cF_T})$, and then $!$-pushing forward to $\Bun_2$. There is a commutative diagram 
\begin{equation}\label{eq: hecke unfold diagram}
\begin{tikzcd}
& \Bun_N^{\cF_T} \times_{\Bun_2} \Hk_{2, \Sym^{\bu}(C)}^{\Std} \ar[dl, "h_l"'] \ar[d, "u"]  \ar[dr, dashed, "h"]  \\ 
\Bun_N^{\cF_T} \ar[d, "u"]  & \Hk_{2, \Sym^{\bu}(C)}^{\Std}  \ar[dl, "h_l"'] \ar[dr, "h_r"] & \ol{\Bun}_N^{\cF_T} \ar[d, "\ol{u}"] \\
\Bun_2   & & \Bun_2
\end{tikzcd}
\end{equation}
where the left square is Cartesian. To describe the dashed arrow $h$, note that 
\[
\Bun_N^{\cF_T} \times_{\Bun_2} \Hk_{2, \Sym^{\bu}(C)}^{\Std} = \bigcup_{d \geq 0} 
\Bun_N^{\cF_T}  \times_{\Bun_2} \Hk_{2, C^{(d)}}^{\Std}
\]
has $R$-points the groupoid of diagrams of the form 
\begin{equation}\label{eq: hecke period HkW1}
\begin{tikzcd}
\Omega_C^{1/2} \boxtimes R \ar[d] \\
\cE \ar[r, hook] \ar[d] & \cE' \\
\Omega_C^{-1/2} \boxtimes R
\end{tikzcd}
\end{equation}
where the cokernel of $\cE \inj \cE'$ is a line bundle over some effective divisor $D \in  C^{(d)} (R)$. The dashed arrow $h$ in \eqref{eq: hecke unfold diagram} sends \eqref{eq: hecke period HkW1} to the inclusion of coherent sheaves $(\Omega_C^{1/2} \boxtimes R \inj \cE') \in  \ol{\Bun}_N^{\cF_T}(R)$ in terms of the description in Example \ref{ex: drinfeld compactification GL_2}.

There is an open embedding $\ol{j} \co U \rightarrow \ol{\Bun}_N^{\cF_T}$ sending \eqref{eq: hecke period diagram} to the inclusion of coherent sheaves $\Omega_C^{1/2} \boxtimes R \inj \cF \otimes \Omega_C^{-1/2}$. We claim that 
\begin{equation}\label{eq: hecke star psi}
h_! h_l^* (\Psi) \cong  \ol{j}_! \ev^* (\AS) \in \Shv(\ol{\Bun}_N^{\cF_T}).
\end{equation}
We will deduce this from the Geometric Casselman-Shalika formula of Frenkel-Gaitsgory-Vilonen \cite[Theorem 4]{FGV01}, which is a similar statement but for Hecke modification at a ``fixed leg''; thus our claim \eqref{eq: hecke star psi} can be viewed as a form of Geometric Casselman-Shalika ``with moving legs''. We first produce a map 
\begin{equation}\label{eq: hecke star psi 2}
\ev^* \AS \rightarrow \ol{j}^* h_! h_l^* (\Psi) \in \Shv(U). 
\end{equation}
Indeed, the fibered product $U'$ in the Cartesian diagram 
\begin{equation}
\begin{tikzcd}
U' \ar[r, "h'"] \ar[d, "\ol{j}'"] & U \ar[d, "\ol{j}"]  \\
 \Bun_N^{\cF_T} \times_{\Bun_2} \Hk_{2, \Sym^{\bu}(C)}^{\Std}  \ar[r, "h"] &  \ol{\Bun}_N^{\cF_T}  \\
\end{tikzcd}
\end{equation}
parametrizes divisors $D \in \Sym^{\bu}(C)(R)$ along with commutative diagrams 
\begin{equation}
\begin{tikzcd}
\Omega_C^{1/2}  \boxtimes R\ar[r, equals] \ar[d] &  \Omega_C^{1/2}  \boxtimes R\ar[d]  \\
\cE \ar[r, hook] \ar[d]  & \cE' \ar[d] \\
\Omega_C^{-1/2} \boxtimes R \ar[r, hook ] & \Omega_C^{-1/2}  \otimes_{\cO_C} \Cal{L}_2
\end{tikzcd}
\end{equation}
where $\coker(\Omega_C^{-1/2}  \boxtimes R \rightarrow \Omega_C^{-1/2}  \otimes \Cal{L}_2) \cong \coker(\cE \inj \cE')$ is a line bundle over $D$. The map $h' \co U' \rightarrow U$ projects to the data of $D$ along with the right column. But this data is enough to reconstruct the entire diagram, as the left column must then be the pullback of the right column along the composition 
\[
\Omega_C^{-1/2} \boxtimes R   \cong  \Omega_C^{-1/2} \otimes_{\cO_C} \cL_2(-D) \inj \Omega_C^{-1/2}  \otimes_{\cO_C} \Cal{L}_2.
\]
Hence the map $h' \co U' \xrightarrow{\sim} U$ is an isomorphism. This induces isomorphisms 
\begin{equation}\label{eq: Hk open eq 1}
\ev^* \AS \rightarrow  h'_! (\ol{j}')^* h_l^* ( \Psi ) \cong \ol{j}^* h_! h_l^* (\Psi ) \in \Shv(U),
\end{equation}
giving the claimed map \eqref{eq: hecke star psi 2}. 

As $\ol{j}$ is an open embedding, adjunction from \eqref{eq: Hk open eq 1} gives a map 
\begin{equation}
\ol{j}_! \ev^* \AS \rightarrow  h_! h_l^* (\Psi)  \in \Shv(\ol{\Bun}_N^{\cF_T})
\end{equation}
which we want to show is an isomorphism. This can be checked at stalks, so we may in particular restrict to the fiber over $D \in C^{(d)} \subset \Sym^{\bu} C$ with respect to the map $\ol{\Bun}_N^{\cF_T} \rightarrow \Sym^{\bu} C$. Now this is a special case of \cite[Theorem 4]{FGV01} (applied iteratively to the different points in the support of $D$), which says that the Hecke operator attached to $V \in \Rep(G)$ acting on $\Psi$ is the \emph{clean} extension of the stratum indexed by the highest weight. Indeed, the highest weight of $\Sym^d (\Std)$ is $(d,0)$, and the stratum indexed by $(d,0)$ in $\ol{\Bun}_N^{\cF_T}$ is the one where $\Omega_C^{1/2} \boxtimes R\rightarrow \Cal{E}'$ has no zeros, or in other words is already saturated, so that $\coker(\Omega_C^{1/2} \boxtimes R \inj \cE')$ is locally free of rank $d$ over $D$. This is precisely the fiber of $U$ over $D \in C^{(d)}$, which then establishes \eqref{eq: hecke star psi} by invoking the Geometric Casselman-Shalika formula, noting by Remark \ref{rem: satake of sym std} that the perverse sheaf associated to $V$ by the Geometric Satake equivalence is the shift and twist of the constant sheaf that we are using. 

Applying $\ol{u} \co \ol{\Bun}_N^{\cF_T} \rightarrow \Bun_2$ to both sides of \eqref{eq: hecke star psi} gives an isomorphism 
\begin{equation}\label{eq: Hk open eq 2}
\ol{u}_! h_! h_l^* (\Psi) \cong  \ol{u}_! \ol{j}_! \ol{j}_! (\ev^* \AS).
\end{equation}
From the commutativity of the diagram \eqref{eq: hecke unfold diagram}, the left side of \eqref{eq: Hk open eq 2} identifies with $\Hk_{2, \Sym^{\bu}(C)}^{\Std} \tw{-\deg} \star \sW \in \Shv(\Bun_G)$ since the kernel sheaf \eqref{eq: Hecke Hk kernel} is the constant sheaf shifted by $\tw{d}$ on the component over $C^{(d)}$. On the other hand, we see by inspection $\ol{u} \circ  \ol{j} = f \circ \pr_1 \circ j$, so the right side of \eqref{eq: Hk open eq 2} identifies with $f_! j_! j^* (\ev^* \AS) \in \Shv(\Bun_G)$. This completes the proof. 
\end{proof}

\subsubsection{Summary} Putting together \S \ref{ssec: hecke aut closed} and \S \ref{ssec: hecke aut open} we have produced an exact triangle in $\Shv(\Bun_G)$, 
\begin{equation}\label{eq: hecke triangle 2}
(\Hk_{2, \Sym^{\bu} C}^{\Std}) \star \sW \tw{2r_E -\deg}   \rightarrow \pi_{\aut !} (\ul{k}_{\BunX} ) \rightarrow \Eis_B  (\ul{k}_{\Bun_1} \boxtimes \delta_{\mrm{spin}}) \tw{2r_E} 
\end{equation}

\subsection{Comparison}

We compare the extensions \eqref{eq: hecke triangle 1} and \eqref{eq: hecke triangle 2} for the $L$-sheaf and period sheaf. 

\subsubsection{Rightmost terms} First we examine the terms on the right. By \eqref{eq: Tate whittaker normalization} and \eqref{eq: Tate perverse normalization}, we have 
\[
\Dmod(\Bun_1 \times \Bun_1) \ni   (\ul{k}_{\Bun_1} \boxtimes \delta_{\triv})  \tw{2g-2}\xmapsto{\LL_{\GL_1 \times \GL_1}} \cO_{\triv}  \boxtimes \omega_{\Loc_1}  \in \QCoh(\Loc_1 \times \Loc_1).
\]
Note that we are using $\QCoh(\Loc_1 \times \Loc_1) = \IndCoh_{\Nilp}(\Loc_1 \times \Loc_1)$ and the notational convention of \S \ref{ssec: singular support}. In this case Eisenstein compatibility (\S \ref{ssec: Eisenstein compatibility}) says that the diagram 
\[
\begin{tikzcd}[column sep = huge]
\Dmod(\Bun_1 \times \Bun_1) \ar[d, "\Eis_{\chB}"] \ar[r, "\LL_{\GL_1 \times \GL_1}", "\sim"'] & \QCoh(\Loc_{1} \times \Loc_1) \ar[d, "\Eis_{\chB}^{\spec}"] &   \\
 \Dmod(\Bun_2)  \ar[r, "\LL_G"] & \IndCoh_{\Nilp}(\Loc_{2}) 
\end{tikzcd}
\]
commutes after translating by $(\Omega_C^{1/2}, \Omega_C^{-1/2})$ on the source and twisting by $\tw{g-1 + \deg}$, so we deduce that
\begin{equation}\label{eq: Hecke right comparison}
\Dmod(\Bun_2) \ni \Eis_B (\ul{k}_{\Bun_1} \boxtimes \delta_{\mrm{spin}})\tw{3g-3+\deg}  \xmapsto{\LL_G} \Eis^{\spec}_{\chB} (\cO_{\triv}  \boxtimes \omega_{\Loc_1})  \in \IndCoh_{\Nilp}(\Loc_{2}) .
 \end{equation}
 
\subsubsection{Leftmost terms} Next we examine the terms on the left.

\begin{lemma}\label{lem: open match}
With notation as in Proposition \ref{prop: Hecke Loc} and Proposition \ref{prop: Hecke-Whittaker}, we have 
\[
\Dmod(\Bun_G) \ni \Hk_{2, \Sym^{\bu} C}^{\Std} \star \sW  \tw{3(g-1)} \xmapsto{\LL_G}   \Loc^{\spec} (\Fact(\Sym^{\bu} \Std)) \otimes \omega_{\LocG} \in \QCoh(\LocG).
\]
\end{lemma}

\begin{proof} The Whittaker normalization \eqref{eq: whittaker normalization} in this case reads 
\begin{equation}\label{eq: Hecke whit normalization}
\Dmod(\Bun_G) \ni \sW \tw{-(g-1) } \xmapsto{\LL_G} \omega_{\LocG} \tw{-(g-1) 4} \in \IndCoh_{\Nilp}(\LocG).
\end{equation}
By the special case of Hecke compatibility of $\LL_G$ in Example \ref{ex: Hecke compatibility}, $\LL_G$ intertwines the tensoring action of 
\[
 \Loc^{\spec}(\Fact(\Sym^{\bu} \Std)) \in \QCoh(\LocG) \text{ on } \IndCoh_{\Nilp}(\LocG)
 \]
with the Hecke convolution action of 
\[
\prod_{n \geq 0} \ulk_{\Hk_{G, C^{(n)}}^{\Std}} \tw{n} \in \Shv(\Hk_{G, \Sym^{\bu} C}^{\Std}) \text{ on }\Dmod(\Bun_G).
\]
Acting by these operators on either side of \eqref{eq: Hecke whit normalization} and twisting by $\tw{4(g-1)}$ yields the result.
\end{proof}

\subsubsection{Conclusion} Twisting \eqref{eq: hecke triangle 2} by $\tw{5(g-1)}$, we get an exact triangle in $\Shv(\Bun_G)$:
\[
(\Hk_{2, \Sym^{\bu} C}^{\Std}) \star \sW \tw{3(g-1)+\deg}   \rightarrow \pi_{\aut !} (\ul{k}_{\BunX} ) \tw{5(g-1)} \rightarrow \Eis_B  (\ul{k}_{\Bun_1} \boxtimes \delta_{\mrm{spin}}) \tw{2(g-1) + 2 \deg}  
\]
where we used that $2r_E +5(g-1) = 2(1-g + \deg) + 5(g-1)= 3(g-1) + 2 \deg$. By inspection, the twist $\tw{\deg}$ on the flanking terms correspond to the shearing operation $(-)^{\shear}$ on $\Loc^{\spec} (\Fact(\Sym^{\bu} \Std)) $ and $\Eis^{\spec}_{\check{B}}( \cO_{\triv} \boxtimes \omega_{\Loc_{1}}) \in \IndCoh_{\Nilp}(\LocG)$ for their respective $\G_m$-actions (coming from Proposition \ref{prop: Hecke Loc} in the former case), under Lemma \ref{lem: open match} and \eqref{eq: Hecke right comparison}. Therefore, we have proved the following theorem. 

\begin{thm}\label{thm: Hecke} The Geometric Langlands equivalence $\LL_G \co \Dmod(\Bun_2) \rightarrow \IndCoh_{\Nilp}(\Loc_2) $ intertwines the indicated objects in the diagram below, where each row is exact
\begin{equation}
\begin{tikzcd}[column sep = tiny
]
\Hk_{2, \Sym^{\bu} C}^{\Std} \star \sW  \tw{3(g-1) + \deg}  \ar[r] \ar[d, "\LL_G", "\text{Lemma \ref{lem: open match}}"', mapsto] &  \cP_{X} \tw{5(g-1)} \ar[d, "\LL_G?", dashed, mapsto] \ar[r] &  \Eis_B (\delta_{\mrm{spin}} \boxtimes \ul{k}_{\Bun_1}) \tw{3(g-1) + 2 \deg}\ar[d, "\LL_G", "\eqref{eq: Hecke right comparison}"', mapsto] \\
\Loc^{\spec} (\Fact(\Sym^{\bu} \Std))^{\shear} \ar[r] &   \cL_{\chX}  \ar[r] &   
\Eis^{\spec}_{\check{B}}( \cO_{\triv} \boxtimes \omega_{\Loc_{1}})^{\shear} 
\end{tikzcd}
\end{equation}
\end{thm}

\subsection{Microlocal interpretation}\label{ssec: microlocal} Here we give an explanation for the unfolding decomposition \eqref{eq: Hecke aut exact triangle} in terms of the geometry of $T^*X$. The fact that the stratification ultimately was seen in the Fourier dual space suggests that in general one should look for decompositions in the \emph{microlocal} geometry of $X$. Indeed, we shall see a decomposition of its hyperspherical variety $T^*X$ analogous to \eqref{eq: Lagrangian decomposition} that reflects the decomposition of periods. We thank David Ben-Zvi and Akshay Venkatesh for explaining this to us. 

Let $\mf{a} := \Lie A$ and $\mf{n} := \Lie N$; note that $\Lie \Mir_2 = \mf{a} \oplus \mf{n}$. With respect to the moment maps $\mu_A \co T^*G \rightarrow \mf{a}^*$ and $\mu_{\Mir} \co T^*G \rightarrow \mf{a}^* \oplus \mf{n}^*$, we have the following decomposition of the hyperspherical variety $T^*X$: 
\begin{equation}\label{eq: Hecke microlocal}
T^* X \cong  \frac{\mu_A^{-1}(0 )}{A} = \frac{\mu_{\Mir}^{-1}(0 \times \mf{n}^*)}{A} = \left(  \frac{\mu_{\Mir}^{-1}(0 \times 0)}{A} \right) \sqcup \left(\frac{\mu_{\Mir}^{-1}(0 \times (\mf{n}^* \setminus 0))}{A}  \right).
\end{equation}
(Note that $\mu^{-1}$ should be formed in the derived sense a priori, although it does not make a difference in this example.) Now observe that 
\[
\frac{\mu_{\Mir}^{-1}(0 \times (\mf{n}^* \setminus 0))}{A}  \cong T^*(G/N, \psi)
\]
is isomorphic to the hyperspherical variety of the Whittaker period (using that $A$ acts freely on $\mf{n}^* \setminus 0$)\footnote{This is a reflection of the equivalence between the Whittaker model and the Kirillov model.}, while $\mu_{\Mir}^{-1}(0 \times 0)/A $ can be thought of as the embedding in $T^*X$ of the Lagrangian correspondence 
\[
\begin{tikzcd}
T^*X  & \ar[l, hook'] X \times_{X/N} T^*(X/N) \ar[r, twoheadrightarrow] & T^*(X/N)  =T^*(G/\Mir_2)
\end{tikzcd}
\]
associated to the map $X \surj X/N$. Thus, \eqref{eq: Hecke microlocal} realizes a decomposition of $T^*X$ into $T^*(G/N,\psi)$ and $T^*(G/\Mir_2)$ in the sense of $G$-Hamiltonian spaces. These terms correspond to the decomposition of $\piaut_! (\ulk_{\BunX})$ in \eqref{eq: hecke triangle 2}.

\section{Singular Rankin-Selberg period}
\subsection{Setup} In this example we will study a duality which is closely related to the \emph{Rankin-Selberg} period. In the Rankin-Selberg period for $\GL_2$, $X^{\sf{RS}}= \GL_2 \times \A^2$ as a spherical variety for $G = \GL_2 \times \GL_2$ and $\chX^{\sf{RS}} = \Std \otimes \Std$. We shall pass to the ``rank one'' locus inside in $\chX^{\sf{RS}}$, which corresponds to passing to a certain open subset inside $X^{\sf{RS}}$. 

The Rankin-Selberg period is a fundamental construction, which is also the starting point for other approaches to automorphic periods such as the Relative Trace Formula. Our particular interest in passing from $\chX^{\sf{RS}}$ to $\chX$ owes to the latter being \emph{singular}. Singular spherical varieties lie beyond the scope of \cite{BZSV} but encompass many useful examples, so we would like a theory that covers them. 
Our analysis here suggests that \cite[Conjecture 12.1.1]{BZSV} extends well to the singular case, with modifications to the normalizations.

\subsubsection{Automorphic side} Let $X = (\GL_2 \times \GL_2) \times^{\GL_2} (\A^2 \setminus 0)$, where the quotient by $\GL_2$ is with respect to its diagonal embedding in $\GL_2 \times \GL_2$, and $G = \GL_2 \times \GL_2$ acting by left translation on the left factor. (In other words, $X$ is the induction of $\A^2 \setminus 0$ along the diagonal subgroup $\GL_2 \inj \GL_2 \times \GL_2$.)

\subsubsection{Spectral side} Let $\chX = (\Std \otimes \Std)^{\rank \leq 1}$ (the locus of rank $\leq 1$ tensors) and $\chG  = \GL_2 \times \GL_2$, acting factorwise through the standard representation.

\subsubsection{$\Ggr$-action} We will take action of $\Ggr$ on $X$ induced by the inverse standard scaling action of $\Ggr$ on $\A^2 \setminus 0$, and the trivial action of $\Ggr$ on $\chX$, which we propose to think of as the un-normalized $\Ggr$-actions. 

This is justified as follows: $X$ appears as an open subvariety of the Rankin-Selberg period $X^{\sf{RS}} := (\GL_2 \times \GL_2) \times^{\GL_2} \A^2$, and $\chX$ appears as a closed subvariety of the dual $\chX^{\sf{RS}} := \Std \otimes \Std$. As $X^{\sf{RS}}$ and $\chX^{\sf{RS}}$ are both smooth and affine, they fall under the scope of \cite{BZSV}, wherein the neutral action of $\Ggr$ on $X^{\sf{RS}}$ is induced by the scaling action on $\A^2$, while the neutral action of $\Ggr$ on $\chX^{\sf{RS}}$ is also the scaling action. Then the un-normalized action of $\Ggr$ on $X^{\sf{RS}}$ is induced by the inverse standard scaling on $\A^2$, while the un-normalized action of $\Ggr$ on $\chX^{\sf{RS}}$ is trivial. Passing to the subvarieties $X \inj X^{\sf{RS}}$ and $\chX \inj \chX^{\sf{RS}}$ then results in the actions of the preceding paragraph. 

\begin{remark} Let us compare this to the outcome of naively using the recipe of \cite{BZSV}. On $X$, the naive neutral action of $\Ggr$ should be induced by the scaling action on $\A^2 \setminus 0$, and the naive neutral action of $\Ggr$ on $\chX$ should be the standard scaling action (we use the adjective ``naive'' because strictly speaking $X$ and $\chX$ fall outside the scope of \cite{BZSV}, so we are performing some naive extrapolation from their formulas). Writing $X \cong \GL_2 \times (\A^2 \setminus 0)$ and taking Haar measure on $\GL_2$ and Lebesgue measure on $\A^2 \setminus 0$, we have $\eta = (1, \det) \co \GL_2 \times \GL_2 \rightarrow \G_m$. Writing $(\chX)^{\mrm{smooth}} = (\A^2 \setminus 0 ) \times^{\G_m} (\A^2 \setminus 0)$ and descending Lebesgue measure, we get $\check{\eta} = (\det, \det) \co \GL_2 \times \GL_2 \rightarrow \G_m$. Hence the ``naive un-normalized'' action of $\Ggr$ on $X$ is trivial, and the ``naive un-normalized'' action of $\Ggr$ on $\chX$ is also trivial. 
\end{remark}
 
\subsubsection{Automorphic period} We define the mirabolic subgroup 
\[
\Mir_2  := \left\{ \begin{pmatrix} 1 & * \\ & * \end{pmatrix}  \right\} \subset \GL_2.
\]
Then $\A^2 \setminus 0 \cong \GL_2/\Mir_2$ as a $\GL_2$-space, so $X = (\GL_2 \times \GL_2)  \times^{\GL_2} (\GL_2/\Mir_2) $ as a $G = \GL_2 \times \GL_2$-space. Thus $X/G \cong \rB\Mir_2$. Examining the recipe of \cite[\S 10.2]{BZSV}, we see that $\BunX$ can be identified with the space $\Bun_{\Mir_2}^\Omega$ with $R$-points the groupoid of extensions 
\begin{equation}\label{eq: RS Bun_Mir point}
\Omega_C^{1/2} \boxtimes R \rightarrow \cF \rightarrow \cL \boxtimes \Omega_C^{1/2}
\end{equation}
where $\cL \in \Bun_1(R)$ and $\Omega_C^{1/2}$ is the fixed spin structure (cf. \S \ref{ssec: spin}), in such a way that the map $\piaut\co \BunX \rightarrow \Bun_G$ is identified with the map sending \eqref{eq: RS Bun_Mir point} to $(\cF, \cF) \in \Bun_2 \times \Bun_2(R)$. The automorphic period sheaf is
\[
\PX := \piaut_! (\ulk_{\BunX}) \in \Shv(\Bun_G).
\]
In the language of \cite{BZSV}, this is the un-normalized period sheaf for the un-normalized action of $\Ggr$.

\subsubsection{Spectral period} By definition, $\LocX$ is the derived mapping stack
\[
\LocX = \Map(C_{\dR}, \chX/\chG)
\]
with $R$-points the groupoid of triples $(F_1, F_2,s)$, where: 
\begin{itemize}
\item $F_1, F_2 \in \Loc_2(R)$ are rank 2 vector bundles on $C_R$ with flat connections along $C$.
\item $s \in \hSect(C_{\dR}, (F_1 \otimes F_2)^{\rank \leq 1})(R)$, or informally speaking a flat section of the flat bundle $F_1 \otimes F_2$ landing pointwise in the locus of tensors of rank $\leq 1$. 
\end{itemize}
The (co-localized) spectral period sheaf is 
\[
\LX \cong (\pispec_* \omega_{\LocX})^{\shear} \in \IndCoh_{\Nilp}(\LocG).
\]
(Note that $\pispec_* \omega_{\LocX} = (  \pispec_* \omega_{\LocX})^{\shear}$ in this case, since the action of $\Ggr$ is trivial.) In the language of \cite{BZSV}, this is the un-normalized $\cL$-sheaf for the un-normalized action of $\Ggr$. 

\subsubsection{Duality} Since $\chX$ is singular, \cite[Conjecture 12.1.1]{BZSV} does not formally apply to this case. However, we propose that its statement should be extended as follows.  

\begin{conj} Under the Geometric Langlands equivalence for $G =\GL_2 \times \GL_2$, we have 
\[
\Dmod(\Bun_G) \ni \PX \tw{8(g-1)} \xmapsto{\LL_G} \LX \in \IndCoh_{\Nilp}(\LocG).
\]
\end{conj}

\subsection{Spectral side}

We let $Z \inj \LocX$ be the zero section and $U$ be the open complement of $Z$. We consider the exact triangle coming from filtering by sections with set-theoretic support in $Z$: 
\begin{equation}\label{eq: RS spectral triangle}
(\pispec)_{Z *}^{\wedge} (\omega_{\LocX}) \rightarrow \pispec_* (\omega_{\LocX}) \rightarrow (\pispec|_U)_* \omega_{U}.
\end{equation} 

\subsubsection{Open stratum}\label{sssec: RS spectral open}
We analyze the rightmost term in \eqref{eq: RS spectral triangle}, which corresponds to the open substack $U \subset \LocX$ where the flat section $s$ has rank \emph{exactly} equal to $1$. The locus of rank 1 vectors in the representation $\Std \otimes \Std$ of $\GL_2 \times \GL_2$ is homogeneous. Taking $\begin{pmatrix} 1 \\ 0 \end{pmatrix} \otimes \begin{pmatrix} 1 \\ 0 \end{pmatrix}$ for the basepoint, the stabilizer is generated by the subgroups $\Mir_2 \times \Mir_2$ and $A$ in $\chG$, where 
\[
\Mir_2 = \left\{\begin{pmatrix} 1 & * \\ 0 & * \end{pmatrix} \right\}  \subset \GL_2, \quad \text{and} \quad 
A  =  \left\{ \begin{pmatrix} a & \\ & 1 \end{pmatrix} , \begin{pmatrix} a^{-1} & \\ & 1 \end{pmatrix} \right\} \subset \GL_2 \times \GL_2.
\]
This choice of basepoint induces a presentation 
\[
(\A^2 \otimes \A^2)^{\rank = 1} \cong (\GL_2/\Mir_2)  \times^{A} (\GL_2/\Mir_2)
\]
identifying $U$ with the derived mapping stack $\Map(C_{\dR}, \rB\check{H}) = \LocH$ for the subgroup
\begin{equation}\label{eq: RS H}
\check{H} := \left\{ \begin{pmatrix} a & b \\ 0  & d \end{pmatrix} , \begin{pmatrix} a^{-1} & b' \\ 0 & d' \end{pmatrix}  \right\} \subset \GL_2 \times \GL_2 = \check{G}
\end{equation}
and the map $\LocX \rightarrow \LocG$ is the one induced by the homomorphism $\check{H} \inj \chG$. 

We will rewrite the contribution from this orbit in terms of Eisenstein series. Let $\chB$ be the standard parabolic in $\GL_2$ and $\chT$ its Levi quotient. Then we have a commutative diagram 
\begin{equation}\label{diag: sRS spectral open}
\begin{tikzcd}
& U \ar[dr] \ar[dl, "q'"'] \ar[dr, "g'"]  \\ 
 \Loc_A  \times \Loc_1 \times \Loc_1 \ar[dr, "g"'] & & \Loc_{\check{B}} \times \Loc_{\check{B}}  \ar[dl, "q"]  \ar[dr, "p"] \\
& \Loc_{\check{T}} \times \Loc_{\check{T}}    & & \LocG
\end{tikzcd}
\end{equation}
where the left square is derived Cartesian. Here: 
\begin{itemize}
\item The map $q$ is induced by the Levi quotient homomorphism $\chB \times \chB \rightarrow \chT \times \chT$; in coordinates, it sends
\[
\begin{pmatrix} a & b \\ 0 & d \end{pmatrix} \times \begin{pmatrix} a' & b' \\ 0 & d' \end{pmatrix} \mapsto  (a,d) \times (a',d').
\]
\item The map $q'$ is induced by the group homomorphism $\chH \rightarrow A \times \G_m \times \G_m$ sending \eqref{eq: RS H} to $a \times (d,d') \in A \times \G_m \times \G_m$. 
\item The map $g$ is induced by the group homomorphism $A \times \G_m \times \G_m \rightarrow \chT \times \chT$ which sends
\[
A \times \G_m \times \G_m \ni (a,d,d') \mapsto (a,d) \times ( a^{-1} ,d') \in \chT \times \chT.
\]
\item The map $g'$ is induced by the inclusion of $\chH$ in $\chB \times \chB$. 
\end{itemize}

By base change applied to the commutative diagram \eqref{diag: sRS spectral open}, we have isomorphisms
\begin{equation}\label{diag: sRS spectral open 0}
g'_* (\omega_U) \cong g'_* (q')^! ( \omega_{\Loc_A \times \Loc_1 \times \Loc_1
})  \cong q^! g_* (\omega_{\Loc_A \times \Loc_1 \times \Loc_1} ) \in \IndCoh(\LocB \times \LocB).
\end{equation}
Then applying $p_*$ to \eqref{diag: sRS spectral open 0} and co-localizing, we obtain an isomorphism
\begin{equation}\label{eq: RS spectral open}
(\pispec|_U)_* \omega_{U} \cong  p_* q^! g_* (\omega_{\Loc_A \times \Loc_1 \times \Loc_1}) = \Eis_{\chB \times \chB}^{\spec}(g_* \omega_{\Loc_A \times \Loc_1 \times \Loc_1} )\in \IndCoh_{\Nilp}(\LocG).
\end{equation}

\subsubsection{Closed stratum}  We begin with a general construction that amalgamates two graded (co)commutative (co)algebras. 

\begin{constr}\label{constr: diagonal graded algebra} Let $\Lambda$ be a finite rank free abelian monoid.

(1) Suppose $A= \bigoplus_{\lambda \in \Lambda } A_\lambda$ and $B = \bigoplus_{\lambda \in \Lambda} B_\lambda$ are $\Lambda$-graded commutative algebras in a symmetric monoidal category $(\msf{C}, \otimes)$. Then 
\[
A \diagotimes B := \bigoplus_{\lambda \in \Lambda} A_\lambda \otimes B_\lambda
\]
has a $\Lambda$-graded commutative algebra structure, with multiplication being the composition 
\begin{align*}
(A \diagotimes B) \otimes (A \diagotimes B)& = \left( \bigoplus_{i \in \Lambda} A_i \otimes B_i \right) \otimes
\left( \bigoplus_{j \in \Lambda} A_j \otimes B_j \right)  \\
 & \cong  \bigoplus_{i,j  \in \Lambda }  \left( (A_i \otimes B_i) \otimes  (A_j \otimes B_j)  \right) \\
& \inj  \bigoplus_{\lambda \in \Lambda} \left( \left( \bigoplus_{i+j=\lambda} A_i \otimes A_j \right) \otimes \left( \bigoplus_{i'+j'=\lambda} B_{i'} \otimes B_{j'} \right) \right) 
  \\
& \rightarrow \bigoplus_{\lambda \in  \Lambda} A_\lambda \otimes B_\lambda  = A \diagotimes B
\end{align*}
where the third line includes into the summands indexed by $i=i', j=j'$ and the fourth line uses the multiplications on $A,B$. 

(2) Similarly, if $A$ and $B$ are $\Lambda$-graded cocommutative coalgebras in $(\msf{C}, \otimes)$, then 
\[
A \diagotimes B := \bigoplus_{\lambda  \in \Lambda} A_\lambda \otimes B_\lambda
\]
has a $\Lambda$-graded cocommutative coalgebra structure, by the dual formulas. If $A_\lambda, B_\lambda$ are dualizable in $\msf{C}$ for all $\lambda$, then the graded duals $A^\vee := \bigoplus_{\lambda \in \Lambda}  A_\lambda^\vee$ and  $B^\vee := \bigoplus_{\lambda \in \Lambda} B_\lambda^\vee$ are $\Lambda$-graded commutative algebras, and $(A \diagotimes B)^\vee \cong (A^\vee \diagotimes B^\vee)$. 

\end{constr}

Next we analyze the leftmost term of \eqref{eq: RS spectral triangle}. We may view $Z \inj \LocX$ as the zero section of $\pispec \co \LocX \rightarrow \LocG$. 

We compute the coalgebra of distributions on $X$ supported at $0$. Note that unlike in the {\sf{Tate}} and {\sf{Hecke}} periods considered in the previous sections, $\chX$ is not a vector space, so Lemma \ref{lem: formal thickening of dvb} does not apply. 

\begin{lemma}\label{lem: RS distribution}
As a cocommutative coalgebra in $(\Rep(\GL_2) \otimes \Rep(\GL_2))^{\Z_{\geq 0}}$, we have 
\[
\Gamma_0(\chX; \omega_{\chX})  \cong  \Sym^{\bu} (\Std) \diagotimes \Sym^{\bu} (\Std)
\]
where $\Sym^{\bu} (\Std) \in \Rep(\GL_2)^{\Z_{\geq 0}}$ is the free $\Z_{\geq 0}$-graded cocommutative coalgebra on the standard representation placed in graded degree $1$.  
\end{lemma}

\begin{proof}
The variety $\chX$ is a quadric cone in $\A^4$. Let $r \co \chY \rightarrow \chX$ be the blowup of $\chX$ at the origin. We will use the following properties of $r$: 
\begin{itemize}
\item It is a smooth resolution, with exceptional fiber $W = r^{-1}(0)$ being a smooth quadric in $\PP^3$. 
\item It is a \emph{rational} resolution (as quadric cone singularities are rational in general) so the natural map $r_* \omega_{\chY}  \rightarrow \omega_{\chX}$ is an isomorphism. Dually, the natural map $\cO_{\chX} \rightarrow r_* \cO_{\chY}$ is an isomorphism. 
\item It restricts to an isomorphism on the complement of $0 \in \chX$. 
\end{itemize}

By adjunction, we have a commutative diagram of exact triangles in $\IndCoh(\chX)$,
\[
\begin{tikzcd}
 r_{W*}^\wedge(\omega_{\chY}) \ar[d, "\sim"] \ar[r] &  r_* \omega_{\chY} \ar[r] \ar[d, "\sim"] &   r|_{(\chY \setminus W)*} (\omega_{\chY \setminus W}) \ar[d, equals] \\
\omega_{\cV} \ar[r] & \omega_{\chX} \ar[r]  & \omega_{\chX \setminus 0} 
\end{tikzcd}
\]
where $\cV$ is the derived formal completion of $\chX$ at $0$. Hence we obtain that the map $\Gamma_W(\chY; \omega_{\chY} ) \rightarrow \Gamma_0(\chX; \omega_{\chX})$ is an isomorphism of cocommutative coalgebras. It therefore suffices to compute $\Gamma_W(\chY; \omega_{\chY} )$ as a $\chG$-equivariant graded cocommutative coalgebra. To do this we use \eqref{eq: derived thickenings}; note that $W \inj \chY$ is a regular embedding (unlike $0 \inj \chX$), so that its system of derived infinitesimal thickenings coincides with the classical such notion (whereas the analogous picture for $0 \inj \chX$ is in principle complicated, by Remark \ref{rem: non-LCI infinitesimal thickening}).

First, we compute the (underived) global functions $\cO(\chY)$ as a commutative algebra in $\Rep(\chG)$. Since $\cO_{\chX} \xrightarrow{\sim} r_* \cO_{\chY}$, we have $\cO(\chX) \xrightarrow{\sim} \cO(\chY) \in \Rep(\chG)$. Since $\chX$ is normal (it is $S_2$ like any hypersurface in $\A^n$, and regular in codimension 1) and its cone point is of codimension $2$, the restriction map $\cO(\chX) \rightarrow \cO(\chX \setminus 0)$ is an isomorphism. Noting that $\chX \setminus 0$ is the subspace of $(\A^2 \otimes \A^2)$ having rank exactly one, we have a presentation 
\[
[(\A^2 \setminus 0)\times (\A^2 \setminus 0)] / \G_m \xrightarrow{\sim} \chX \setminus 0 
\]
\[
  (v_1, v_2)  \mapsto  v_1 \otimes v_2
\]
where $\G_m$ acts via the standard scaling on the first factor and its inverse on the second factor. Now, $\cO(\A^2 \setminus 0) \xleftarrow{\sim} \cO(\A^2) \cong \Sym^{\bu} (\Std^\vee)$. Therefore, 
\begin{align}\label{eq: RS coalg 2}
\cO(\chX \setminus 0) & \cong \left(\Sym^{\bu} (\Std^\vee)  \otimes \Sym^{\bu} (\Std^\vee) \right)^{\G_m} \nonumber \\
& \cong \bigoplus_{n \geq 0}  \Sym^n(\Std^\vee) \boxtimes \Sym^n(\Std^\vee)
\end{align}
where in the last isomorphism we note that $\Sym^n(\Std^\vee) \boxtimes 1$ has weight $n$ for the $\G_m$-action and $1 \boxtimes \Sym^n(\Std^\vee)$ has weight $-n$ for the $\G_m$-action. 

Now we turn to calculating the cocommutative coalgebra $\Gamma_W(\chY; \omega_{\chY}) \in \Rep(\GL_2) \otimes \Rep(\GL_2)$. Let 
\[
W = W^{(0)} \inj W^{(1)} \inj \ldots \inj \cW \inj \chY
\]
be the system of formal neighborhoods of $W$ in $\chY$ as in \eqref{eq: derived thickenings}, with colimit the formal completion $\cW := Y_W^{\wedge}$. According to Lemma \ref{lemma: omega closed support} we have 
\[
\Gamma_W(\chY; \omega_{\chY}) \cong \Gamma(\cW; \omega_{\cW}) \in \Rep(\GL_2) \otimes \Rep(\GL_2),
\]
which in turn can be presented as 
\[
\Gamma(\cW; \omega_{\cW}) \cong \colim_n \Gamma(W^{(n)}; \omega_{W^{(n)}}).
\]
By Grothendieck duality, $\Gamma(W^{(n)}; \omega_{W^{(n)}})$ is the dual cocommutative coalgebra to the commutative algebra $\Gamma(W^{(n)}; \cO_{W^{(n)}})$. Since $W^{(n)}$ is the $n$th order neighborhood of $W \cong \PP^1 \times \PP^1$ (being a smooth quadric surface in $\PP^3$), its structure sheaf has no higher cohomology. Since $W \inj \chY$ is a regular embedding, $\cO_{W^{(n)}}$ is the quotient of \eqref{eq: RS coalg 2} by the $n$th power of the ideal of $W$. Dualizing and taking the colimit over $n$ gives the dual cocommutative coalgebra to \eqref{eq: RS coalg 2}, completing the proof. 
\end{proof}

\begin{prop}\label{prop: sRS Loc} Viewing $\Sym^{\bu} (\Std) \diagotimes \Sym^{\bu} (\Std)$ as a cocommutative coalgebra in $\Rep(\chG)^{\Z_{\geq 0}}$ with $\Std$ in grading degree $1$, let $\Fact(\Sym^{\bu} \Std \diagotimes \Sym^{\bu} \Std) \in \QCoh((\rB \chG)_{\Sym^{\bu} C})$ be the associated factorization algebra, augmented by the unit analogously to \eqref{eq: augment spectral fact algebra}. Then we have an isomorphism
\[
(\pispec)^\wedge_{Z*} (\omega_{\LocX}) \cong \Loc^{\spec}(\Fact(\Sym^{\bu} \Std \diagotimes \Sym^{\bu} \Std)) \otimes \omega_{\LocG}  \in \QCoh(\LocG). 
\]
\end{prop}

\begin{proof}
Applying Proposition \ref{prop: distribution algebra} to the origin inside the $\chG$-scheme $\chX$, with grading coming from the standard scaling action of $\Ggr$, we obtain an isomorphism
\begin{equation}\label{eq: RS non-unital Loc 1} 
\ol{(\pispec)_{Z*}^{\wedge} (\omega_{\LocX})} \cong \pr_{2*} (\Fact(\ol{\Gamma}_0(\chX; \omega))^{\univ} ) \otimes \omega_{\LocG}\in \QCoh(\LocG).
\end{equation}
From Lemma \ref{lem: RS distribution} we have 
\[	
\ol{\Gamma}_0(\chX; \omega) \cong  \ol{\Sym}^{\bu} (\Std) \diagotimes \ol{\Sym}^{\bu} (\Std) \in \Rep(\chG)^{\Z_{>0}}.
\]
Inserting this into \eqref{eq: RS non-unital Loc 1} yields an isomorphism
\begin{equation}\label{eq: RS non-unital Loc 2} 
\ol{(\pispec)_{Z*}^{\wedge} (\omega_{\LocX})} \cong \pr_{2*} (\Fact(\ol{\Sym}^{\bu} (\Std) \diagotimes \ol{\Sym}^{\bu} (\Std))^{\univ}) \otimes \omega_{\LocG} \in \QCoh(\LocG)
\end{equation}
where $\Fact(\ol{\Sym}^{\bu} (\Std) \diagotimes \ol{\Sym}^{\bu} (\Std) )$ is the universally twisted relative factorization algebra over $\LocG$,  analogous to \eqref{eq: augment spectral universal fact algebra}. The rest follows from adding a copy of the unit, as in the proof of Proposition \ref{prop: Tate Loc}. 

\end{proof}



\subsubsection{Summary}

We have produced an exact triangle in $\IndCoh_{\Nilp}(\LocG)$:
\begin{equation}\label{eq: RS triangle 1}
\Loc^{\spec}(\Fact(\Sym^{\bu} \Std \diagotimes \Sym^{\bu} \Std)) \otimes \omega_{\LocG}   \rightarrow \pispec_* (\omega_{\LocX})  \rightarrow  \Eis_{\chB}^{\spec}(g_* \omega_{\Loc_A \times \Loc_1 \times \Loc_1}).
\end{equation}

\subsection{Automorphic side} We seek a filtration of $\cP_X$ corresponding to the filtration of $\cL_{\chX}$ in \eqref{eq: RS triangle 1}, which came from the stratification of $\chX$ by $\chG$-orbits. But, as in the case of the Hecke period, $X$ is a homogeneous $G$-variety, so the stratification of $X$ by $G$-orbits is trivial. Therefore, we will again look for a stratification in the Fourier dual space. 

\subsubsection{Unfolding}\label{sssec: sRS unfolding}
We may identify $\Bun_{\Mir_2}^{\Omega}$ with the derived vector bundle over $\Bun_1$, associated to the perfect complex $\ul{\RHom(\cL^{\univ}, \Omega_C[1])}$ whose pullback to $\cL \in \Bun_1(R)$ is $\RHom_{C_R}(\cL, \Omega_C \boxtimes R[1])$ as an animated $R$-module. We abbreviate $E := \Bun_{\Mir_2}^\Omega$ viewed in this way. Then there is a factorization of $\piaut$ as 
\begin{equation}
\begin{tikzcd}
E \ar[r, "\Delta"] \ar[rr, bend left, "\piaut"] & E \times_{\Bun_1} E \ar[r, "f"] & \Bun_2 \times \Bun_2
\end{tikzcd}
\end{equation}
where $f$ sends 
\begin{equation}
 E \times_{\Bun_1} E(R) \ni \left\{ \begin{tikzcd}[row sep = tiny] \Omega_C \boxtimes R \rightarrow \cF_1 \rightarrow \cL \\  \Omega_C \boxtimes R \rightarrow \cF_2 \rightarrow \cL \end{tikzcd}  \right\} \mapsto \left\{ \begin{tikzcd}[row sep = tiny] \cF_1 \otimes_{\cO_C} \Omega_C^{-1/2} \\\cF_2 \otimes_{\cO_C} \Omega_C^{-1/2}\end{tikzcd}  \right\} \in \Bun_2 \times \Bun_2(R).
 \end{equation}
This induces an isomorphism 
\begin{equation}
\piaut_! (\ul{k}_E) \cong f_! (\Delta_! \ul{k}_E) \in \Shv(\Bun_G).
\end{equation}
From \S \ref{sssec: FT functoriality} we have a natural isomorphism 
\begin{equation}\label{eq: sRS unfold 1}
\Delta_! \FT_{E^\vee}  \cong \FT_{E^\vee \times_{\Bun_1} E^\vee} (\Delta^\vee)^* [r_E](r_E)
\end{equation}
where $r_E := \rank(E) = 1-g + \deg \cL^{\univ}$ is the virtual rank of $E$; and $\Delta^\vee \co  E^\vee \times_{\Bun_1} E^\vee \rightarrow E^\vee $ is dual map of $\Delta$, which can be identified with the addition on $E^\vee$. In particular, applying \eqref{eq: sRS unfold 1} to the $\delta$-sheaf on $E^\vee$, and using that $\FT_{E^\vee} (\delta_{E^\vee}) \cong \ulk_E[r_E]$ (cf. Example \ref{ex: FT delta}), we obtain an isomorphism
\begin{equation}\label{eq: RS strat 1}
\piaut_! (\ul{k}_E) \cong f_! \FT_{E^\vee \times_{\Bun_1} E^\vee}((\Delta^\vee)^* \delta_{E^\vee})  (r_E) \in \Shv(\Bun_G).
\end{equation}
We have a derived Cartesian square 
\begin{equation}\label{eq: sRS unfold 2}
\begin{tikzcd}
E^\vee \ar[r, "z'"] \ar[d] & E^\vee \times_{\Bun_1} E^\vee \ar[d, "\Delta^\vee"] \\
\Bun_1 \ar[r, "z"] & E^\vee
\end{tikzcd}
\end{equation}
where the map $z$ is the zero section, and $z'$ is the anti-diagonal embedding $e \mapsto (e,-e)$. Applying proper base change, we obtain an isomorphism
\begin{equation}\label{eq: sRS unfold 3}
(\Delta^\vee)^* (\delta_{E^\vee}) = (\Delta^\vee)^* z'_! (\ulk_{\Bun_1}) \cong   z'_! ( \ul{k}_{E^\vee})  \in \Shv(E^\vee \times_{\Bun_1} E^\vee). 
\end{equation}
Putting this into \eqref{eq: RS strat 1} gives an isomorphism
\begin{equation}\label{eq: RS strat 2}
\piaut_! (\ul{k}_E) \cong f_! \FT_{E^\vee \times_{\Bun_1} E^\vee}(z'_! \ul{k}_{E^\vee})(r_E) \in \Shv(\Bun_G).
\end{equation}
Now we unravel the meaning of the right side of \eqref{eq: RS strat 2}. Consider the commutative diagram 
\begin{equation}
\begin{tikzcd}
& & E \times_{\Bun_1} E \times_{\Bun_1} E^\vee \ar[dr, "\pr_r'"]  \ar[d, "z''"] \\
& E \ar[d, "\Delta"]  & E \times_{\Bun_1} E \times_{\Bun_1} E^\vee \times_{\Bun_1} E^\vee \ar[dl, "\pr_l"']  \ar[dr, "\pr_r"]  & E^\vee \ar[d, "z'"]  \\
& E \times_{\Bun_1} E \ar[dl, "f"] \ar[dr] & & E^\vee \times_{\Bun_1} E^\vee  \ar[dl] \ar[d, "\Delta^\vee"]  \\
\Bun_2 \times \Bun_2 & & \Bun_1  & E^\vee 
\end{tikzcd}
\end{equation}
where the parallelograms are derived Cartesian. By base change and the projection formula, we have
\begin{align}\label{eq: RS strat 3}
\FT_{E^\vee \times_{\Bun_1} E^\vee}(z'_! \ul{k}_{E^\vee}) & = \pr_{l!} (\ev^* \AS \otimes \pr_r^* z'_!  \ul{k}_{E^\vee})[2r_E] \nonumber \\
& \cong \pr_{l!} (\ev^* \AS \otimes  z''_! \pr_r'^* (\ul{k}_{E^\vee}) )[2r_E] \nonumber \\
& \cong \pr_{l!} (z''_! (z''^* \ev^* \AS))[2r_E].
\end{align} 

Explicitly, $E \times_{\Bun_1} E \times_{\Bun_1} E^\vee$ is a derived Artin stack whose classical truncation has $R$-points the groupoid of diagrams
\begin{equation}\label{eq: RS big diag}
\left\{ \begin{tikzcd}
& \Omega_C \boxtimes R \ar[d] \\
& \cF_1 \ar[d] \\
\cO_C \boxtimes R \ar[r, "s"] & \cL 
\end{tikzcd} \quad 
\begin{tikzcd}
& \Omega_C \boxtimes R \ar[d] \\
& \cF_2 \ar[d] \\
\cO_C \boxtimes R  \ar[r, "-s"] & \cL 
\end{tikzcd} \right\}
\end{equation}
where 
\begin{itemize}
\item $\cF_1$, $\cF_2$ are rank $2$ vector bundles on $C_R$, 
\item $\cL$ is a line bundle on $C_R$,
\item The columns are extensions of $\cL$ by $\Omega_C \boxtimes R$, and 
\item $s \co \cO_{C} \boxtimes R \rightarrow \cL$ is any global section of $\cL$ (possibly zero). 
\end{itemize}
The map 
\begin{equation}
\Pi := f \circ \pr_l \circ z'' \co E \times_{\Bun_1} E \times_{\Bun_1} E^\vee \rightarrow \Bun_G
\end{equation}
sends \eqref{eq: RS big diag} to the underlying bundles $(\cF_1 \otimes_{\cO_C} \Omega_C^{-1/2} , \cF_2 \otimes_{\cO_C} \Omega_C^{-1/2} ) \in \Bun_2 \times \Bun_2(R)$. The map $\ev \circ z''$ sends \eqref{eq: RS big diag} to the sum of the extension classes in $\Ext^1(\cO_C \boxtimes R, \Omega_C \boxtimes R) \cong R$ obtained by pulling back the extensions in the columns via the horizontal maps $s$ and $-s$ respectively. We abuse notation by abbreviating 
\begin{equation}
\ev^* \AS := z''^* \ev^* \AS \in \Shv(E \times_{\Bun_1} E \times_{\Bun_1} E^\vee).
\end{equation}

Putting together \eqref{eq: RS strat 2} and \eqref{eq: RS strat 3}, we have established:

\begin{prop}\label{prop: RS period FT} 
We have an isomorphism 
\[
\piaut_! (\ul{k}_{\BunX}) \cong \Pi_! \ev^* \AS \tw{2r_E} \in \Shv(\Bun_G).
\]
\end{prop}

Now we contemplate the decomposition coming from the stratification by the vanishing of $s$. Let $i \co Z \inj E \times_{\Bun_1} E \times_{\Bun_1} E^\vee$ be the closed substack where $s=0$, which may be viewed as the pullback of the zero section of $E^\vee$ under $\pr_r$. Let $j \co U \inj E \times_{\Bun_1} E \times_{\Bun_1} E^\vee$ be its open complement. Then, using Proposition \ref{prop: RS period FT}, we have the following exact triangle in $\Shv(\Bun_G)$:
\begin{equation}\label{eq: RS period decomp}
\Pi_! j_! j^* (\ev^* \AS)\tw{2r_E} \rightarrow \piaut_! (\ul{k}_{\BunX})  \rightarrow  \Pi_! i_* i^* (\ev^* \AS)\tw{2r_E}.
\end{equation}
Henceforth we will have no further need of the derived structure on the objects involved, so we pass to classical truncations.

\subsubsection{Closed stratum}\label{sssec: RS aut closed} We analyze the rightmost term of \eqref{eq: RS period decomp}. Note that $Z \cong E \times_{\Bun_1} E$ and $i^* \ev^* \AS \cong \ul{k}_Z$. Let $B$ be the standard Borel subgroup of $\GL_2$ containing $\Mir_2$, and $T$ its Levi quotient. We have an equality of maps 
\begin{equation}\label{eq: sRS aut closed 1}
\Pi \circ i = p \circ g'   \co Z \rightarrow \Bun_G
\end{equation}
referring to the diagram 
\begin{equation}
\begin{tikzcd}
E \times_{\Bun_1} E \ar[rr, "\Pi", bend left] \ar[r, "g'"] &  \Bun_B \times \Bun_B \ar[r, "p"] &  \Bun_2 \times \Bun_2
\end{tikzcd}
\end{equation}
where $g'$ sends \eqref{eq: RS big diag} to 
\begin{equation}
(\Omega_C^{1/2} \boxtimes R \inj \cF_1 \otimes \Omega_C^{-1/2} ) \times (\Omega_C^{1/2} \boxtimes R \inj \cF_2 \otimes \Omega_C^{-1/2} ) \in \Bun_B(R) \times \Bun_B(R).
\end{equation}
The map $g'$ fits into a derived Cartesian square 
\begin{equation}
\begin{tikzcd}
E \times_{\Bun_1} E \ar[r, "g'"] \ar[d, "q'"] & \Bun_B \times \Bun_B  \ar[d, "q"]\\
\Bun_1 \ar[r, "g"] & \Bun_T \times \Bun_T
\end{tikzcd}
\end{equation}
where 
\begin{itemize}
\item The map $q$ is induced by the Levi quotient $B \rightarrow T$. 
\item The map $q'$ sends \eqref{eq: RS big diag} to $\cL$. 
\item The bottom map $g$ sends 
\begin{equation}
\cL \in \Bun_1(R) \mapsto (\Omega_C^{1/2} \boxtimes R , \cL \otimes_{\cO_C} \Omega_C^{-1/2}) \times (\Omega_C^{1/2} \boxtimes R, \cL\otimes_{\cO_C} \Omega_C^{-1/2}) \in \Bun_T \times \Bun_T  (R).
\end{equation}
\end{itemize}
Hence by proper base change, we have natural isomorphisms 
\begin{equation}\label{eq: sRS aut closed 2}
g'_! (\ul{k}_Z) \cong g'_! (q')^* (\ul{k}_{\Bun_1}) \cong q^* g_! (\ul{k}_{\Bun_1} )\in \Shv(\Bun_B \times \Bun_B).
\end{equation}
Combining \eqref{eq: sRS aut closed 2} with the identity of maps \eqref{eq: sRS aut closed 1} and the observation $i^* \ev^* \sA \cong \ulk_Z$, we have an isomorphism
\begin{equation}\label{eq: RS aut closed 1}
\Pi_! i_* i^* (\ev^* \AS) \cong p_! g'_! (\ul{k}_Z ) \cong p_! q^* (g_! \ul{k}_{\Bun_1}) \cong \Eis_{B \times B}(g_! \ul{k}_{\Bun_1}).
\end{equation}

\subsubsection{Open stratum}\label{sssec: RS aut open} We analyze the leftmost term in \eqref{eq: RS big diag}.

Twisting by $\Omega_C^{1/2}$, we may rewrite the $R$-points of $U$ (recall that we replaced the original $U$ by its classical truncation) as the groupoid of diagrams
\begin{equation}\label{RS eq: open}
\left\{ \begin{tikzcd}
& \Omega_C^{1/2} \boxtimes R \ar[d] \\
& \cE_1 \ar[d] \\
\Omega_C^{-1/2} \boxtimes R  \ar[r, "-s", hook] & \cL \otimes \Omega_C^{-1/2} 
\end{tikzcd} \hspace{1cm} 
\begin{tikzcd}
& \Omega_C^{1/2}  \boxtimes R \ar[d] \\
& \cE_2 \ar[d] \\
\Omega_C^{-1/2}  \boxtimes R  \ar[r, "s", hook] & \cL \otimes \Omega_C^{-1/2} 
\end{tikzcd} \right\}
\end{equation}
where the columns are extensions of vector bundles and where $s$ is required to be injective as a map of coherent sheaves fiberwise over $\Spec R$.  

Let $\Hk_{2, C^{(d)}}^{\Std}$ and $\Hk_{2, \Sym^{\bu} C}^{\Std}$ be as in \S \ref{ssec: hecke aut open}. Write 
\[
\Hk_{G, C^{(d)}}^{\Std \diagotimes \Std} := \Hk_{2, C^{(d)}}^{\Std} \times_{C^{(d)}} \Hk_{2, C^{(d)}}^{\Std}
\]
which parametrizes pairs of modifications $(\cF_1 \inj \cF_1' \in \Hk_{2, C^{(d)}}^{\Std}(R), \cF_2 \inj \cF_2' \in \Hk_{2, C^{(d)}}^{\Std}(R))$ such that 
\[
\supp(\coker(\cF_1 \inj \cF_1') ) = \supp(\coker(\cF_2 \inj \cF_2') )  \in C^{(d)}(R).
\]
Write
\[
\Hk_{G, \Sym^{\bu} C}^{\Std \diagotimes \Std} := \bigcup_{d \geq 0} 
\Hk_{G, C^{(d)}}^{\Std \diagotimes \Std} \cong \Hk_{2, \Sym^{\bu} C}^{\Std} \times_{\Sym^{\bu} C} \Hk_{2, \Sym^{\bu} C}^{\Std}.
\]
We have a correspondence diagram
\begin{equation}\label{eq: sRS Hk correspondence}
\begin{tikzcd}
& \Hk_{G, \Sym^{\bu} C}^{\Std \diagotimes \Std} \ar[dl, "h_l"'] \ar[dr, "h_r \times \pr"'] \ar[drr, "h_r"] \\
\Bun_G & & \Bun_G \times \Sym^{\bu} C \ar[r] & \Bun_G
\end{tikzcd}
\end{equation}
where 
\begin{align*}
h_l(\cF_1 \inj \cF_1', \cF_2 \inj \cF_2') &:= (\cF_1, \cF_2) \in \Bun_2 \times \Bun_2(R)\\
h_r(\cF_1 \inj \cF_1', \cF_2 \inj \cF_2') &:= (\cF_1', \cF_2') \in \Bun_2 \times \Bun_2(R) \\
\pr(\cF_1 \inj \cF_1', \cF_2 \inj \cF_2') &:= \supp(\coker(\cF_i \inj \cF_i')) \in \Sym^{\bu} C(R).
\end{align*}
Consider the complex 
\begin{equation}\label{eq: sRS kernel sheaf}
\prod_{d \geq 0} \ulk\tw{2d} \in \Shv\left(\bigcup_{d \geq 0} 
\Hk_{G, C^{(d)}}^{\Std \diagotimes \Std} \right) =  \Shv(\Hk_{G, \Sym^{\bu} C}^{\Std \diagotimes \Std})
\end{equation}
as a kernel sheaf on $\Hk_{G, \Sym^{\bu} C}^{\Std \diagotimes \Std}$. We write $\Hk_{2, \Sym^{\bu}(C)}^{\Std \diagotimes \Std} \star (-)$ for the endofunctor $\pr_{r!} (\eqref{eq: sRS kernel sheaf} \otimes \pr_l^*(-))$ of $\Shv(\Bun_G)$ induced by \eqref{eq: sRS kernel sheaf}.

Let $N$ be the unipotent radical of $B$, so the Whittaker space for $G$ is $G/(N \times N, \psi \times \psi)$. Let $\cF_T := \Omega_C^{1/2} \oplus \Omega_C^{-1/2} \in \Bun_T(\F)$. Then, with notation as in \S \ref{sssec: drinfeld compactification}, $\Bun_N^{\cF_T}(R)$ classifies extensions
\[
0 \rightarrow \Omega_C^{1/2} \boxtimes R \rightarrow \cE \rightarrow \Omega_C^{-1/2} \boxtimes R \rightarrow 0.
\]
The Whittaker sheaf $\sW \in \Shv(\Bun_G)$ is obtained by pulling back the Artin-Schreier sheaf $\AS$ via the map $\Bun_N^{\cF_T} \times \Bun_N^{\cF_T} \rightarrow \A^1$, to get $\Psi \in \Shv(\Bun_N^{\cF_T} \times \Bun_N^{\cF_T})$, and then $!$-pushing forward to $\Bun_G$. 

We write $\tw{\deg} \co \Shv(\Bun_G) \rightarrow \Shv(\Bun_G)$ for the endofunctor given by $\tw{d_1+d_2}$ on $\Shv(\Bun_G^{(d_1,d_2)})$, where $\Bun_G^d$ is the connected component parametrizing pairs of bundles $(\cF_1, \cF_2)$ with $\deg \cF_i = d_i$. 

\begin{prop}\label{prop: RS-Whittaker}
With notation as in \S \ref{sssec: sRS unfolding}, we have a natural isomorphism 
\[
\Hk_{G, \Sym^{\bu}(C)}^{\Std \diagotimes \Std} \star \Psi \cong \Pi_! j_! j^*  (\ev^* \AS) \tw{\deg} \in \Shv(\Bun_G).
\] 
\end{prop}

\begin{proof}Consider the commutative diagram 
\begin{equation}\label{eq: RS unfold diagram}
\begin{tikzcd}
& (\Bun_N^{\cF_T}  \times \Bun_N^{\cF_T})\times_{\Bun_G} \Hk_{G, \Sym^{\bu}(C)}^{\Std \diagotimes \Std} \ar[dl, "h_l"'] \ar[d, "u"]  \ar[dr, dashed, "h"]  \\ 
(\Bun_N^{\cF_T} \times \Bun_N^{\cF_T})  \ar[d, "u"]  & \Hk_{G, \Sym^{\bu}(C)}^{\Std \diagotimes \Std}  \ar[dl, "h_l"'] \ar[dr, "h_r"] & \ol{\Bun}_N^{\cF_T} \times \ol{\Bun}_N^{\cF_T}  \ar[d, "\ol{u}"] \\
\Bun_G   & & \Bun_G
\end{tikzcd}
\end{equation}
where: 
\begin{itemize}
\item The left square is Cartesian. 
\item The dashed arrow $h$ is as follows. Observe that $(\Bun_N^{\cF_T}  \times \Bun_N^{\cF_T}) \times_{\Bun_G} \Hk_{G, \Sym^{\bu}(C)}^{\Std}$ parametrizes diagrams of the form 
\begin{equation}\label{eq: RS period HkW1}
\begin{tikzcd}
\Omega_C^{1/2} \boxtimes R \ar[d] \\
\cE_1 \ar[r, hook] \ar[d] & \cE'_1 \\
\Omega_C^{-1/2} \boxtimes R
\end{tikzcd} \hspace{1cm} 
\begin{tikzcd}
\Omega_C^{1/2} \boxtimes R \ar[d] \\
\cE_2 \ar[r, hook] \ar[d] & \cE'_2 \\
\Omega_C^{-1/2} \boxtimes R
\end{tikzcd}
\end{equation}
where $\coker(\cE_1 \inj \cE_1')$ and $\coker(\cE_2 \inj \cE_2')$ have the same divisor of supports, which is some effective divisor $D \in C^{(d)}(R)$. The dashed arrow $h$ sends \eqref{eq: RS period HkW1} to the pair of injections of coherent sheaves 
\[
(\Omega_C^{1/2} \boxtimes R \inj \cE_1', \Omega_C^{1/2} \boxtimes R \inj \cE_2') \in  \ol{\Bun}_N^{\cF_T}\times \ol{\Bun}_N^{\cF_T}(R).
\]
\end{itemize}

There is an open embedding $\ol{j} \co U \rightarrow \ol{\Bun}_N^{\cF_T} \times \ol{\Bun}_N^{\cF_T}$ sending \eqref{RS eq: open} to the inclusions of coherent sheaves $(\Omega_C^{1/2} \boxtimes R \stackrel{s}\inj \cE_1', \Omega_C^{1/2} \boxtimes R \stackrel{s}\inj \cE_2')$. We claim that 
\begin{equation}\label{eq: RS star psi}
h_! h_l^* (\Psi) \cong  \ol{j}_! (\ev^* \AS) \in \Shv(\ol{\Bun}_N^{\cF_T}\times \ol{\Bun}_N^{\cF_T}).
\end{equation}

Consider the commutative diagram 
\begin{equation}\label{eq: RS aut open 1}
\begin{tikzcd}[column sep = tiny]
& (\Bun_N^{\cF_T} \times \Bun_N^{\cF_T}) \times_{\Bun_G} \Hk_{2, \Sym^{\bu}(C)}^{\Std \diagotimes \Std} \ar[ddl, "h_l"', bend right] \ar[d, "\iota"]  \ar[ddr, dashed, "h", bend left]  \\ 
   & (\Bun_N^{\cF_T} \times \Bun_N^{\cF_T}) \times_{\Bun_G}  (\Hk_{2, \Sym^{\bu} C}^{\Std} \times \Hk_{2, \Sym^{\bu} C}^{\Std} )\ar[dl, "\wt{h}_l"] \ar[dr, "\wt{h}"'] &  \\
(\Bun_N^{\cF_T} \times \Bun_N^{\cF_T})   & & (\ol{\Bun}_N^{\cF_T} \times \ol{\Bun}_N^{\cF_T})
\end{tikzcd}
\end{equation}
where the map $\iota$ is the base change along the diagonal map $\Sym C \xrightarrow{\Delta} \Sym C \times \Sym C$, which is the diagonal embedding $C^{(d)} \rightarrow C^{(d)} \times C^{(d)}$ on the $d$th component. 

We will use the Geometric Casselman-Shalika formula \cite[Theorem 4]{FGV01} to analyze the action of the kernel sheaf 
\[
\ul{k}\tw{\deg} := \prod_{d \geq 0} \ulk\tw{2d} \in \Shv \left(\Hk_{2, \Sym^{\bu} C}^{\Std} \times \Hk_{2, \Sym^{\bu} C}^{\Std}\right)
\]
through the correspondence \eqref{eq: RS aut open 1}. The proof of Proposition \ref{prop: Hecke-Whittaker} shows that $\wt{h}_! \wt{h}_l^* (\Psi \otimes \ul{k})$ is the clean extension from the stratum of $ (\ol{\Bun}_N^{\cF_T} \times \ol{\Bun}_N^{\cF_T} )$ parametrizing diagrams 
\begin{equation}\label{eq: HkRS 1}
\begin{tikzcd}
\Omega_C^{1/2} \boxtimes R \ar[r, equals] \ar[d] &  \Omega_C^{1/2} \boxtimes R \ar[d]  \\
\cF_1 \ar[r, hook] \ar[d]  & \cF_1' \ar[d] \\
\Omega_C^{-1/2} \boxtimes R \ar[r, "s_1", "\neq 0"'] & \Cal{L}_1 \otimes \Omega_C^{-1/2}
\end{tikzcd}\hspace{1cm} 
\begin{tikzcd}
\Omega_C^{1/2} \boxtimes R \ar[r, equals] \ar[d] &   \Omega_C^{1/2} \ar[d]  \\
\cF_2 \ar[r, hook] \ar[d]  & \cF_2' \ar[d] \\
\Omega_C^{-1/2} \boxtimes R \ar[r, "s_2", "\neq 0"'] & \Cal{L}_2 \otimes \Omega_C^{-1/2}
\end{tikzcd}
\end{equation}
Therefore, base changing along the diagonal map $\iota \co C^{(d)} \rightarrow C^{(d)} \times C^{(d)}$, we find that $h_! h_l^* (\Psi) $ is the clean extension along the locus where $\coker(\cF_1 \inj \cF_1')$ has the same divisor of supports as $\coker(\cF_2 \inj \cF_2')$. This is precisely $U$, which completes the proof of the claim \eqref{eq: RS star psi}.

Applying $\ol{u} \co \ol{\Bun}_N^{\cF_T}  \times \ol{\Bun}_N^{\cF_T} \rightarrow \Bun_G$ to both sides of \eqref{eq: RS star psi} gives an isomorphism 
\begin{equation}\label{eq: RS open eq 2}
\ol{u}_! h_! h_l^* (\Psi) \cong  \ol{u}_! \ol{j}_! \ol{j}_! (\ev^* \AS) \in \Shv(\Bun_G). 
\end{equation}
From the commutativity of the diagram \eqref{eq: RS unfold diagram}, the left side of \eqref{eq: RS open eq 2} identifies with $\Hk_{G, \Sym^{\bu}(C)}^{\Std \diagotimes \Std} \tw{-\deg} \star \sW \in \Shv(\Bun_G)$ since the kernel sheaf \eqref{eq: sRS kernel sheaf} is the constant sheaf shifted by $\tw{2d}$ on the component over $C^{(d)}$. On the other hand, by inspection $\ol{u} \circ  \ol{j} = \Pi \circ j$, so the right side of \eqref{eq: RS open eq 2} identifies with $\Pi_! j_! j^* (\ev^* \AS) \in \Shv(\Bun_G)$. This completes the proof. 
\end{proof}

\subsubsection{Summary}In summary, combining \S \ref{sssec: RS aut closed} and \S \ref{sssec: RS aut open} we have produced an exact triangle in $\Shv(\Bun_G)$, 
\begin{equation}\label{eq: RS triangle 2}
\Hk_{G, \Sym^{\bu}(C)}^{\Std \diagotimes \Std}  \star \sW \tw{2r_E - \deg} \rightarrow \piaut_! (\ulk_{\BunX}) \rightarrow \Eis_{B \times B}(g_! \ul{k}_{\Bun_1})\tw{2r_E}.
\end{equation}

\subsection{Comparison}
We compare \eqref{eq: RS triangle 1} and \eqref{eq: RS triangle 2} for the $L$-sheaf and period sheaf.

\subsubsection{Rightmost terms} First we examine the terms on the right. On the spectral side, let us rename the object in question by $\Eis^{\spec}_{\chB \times \chB}(g_*^{\spec} \omega_{\Loc_A \times \Loc_1 \times \Loc_1})$, as the same name $g$ is regrettably used for a different map on the automorphic side. To spell out the meaning of $g_*^{\spec} (\omega_{\Loc_A \times \Loc_1 \times \Loc_1})$, let us label the two factors of the maximal torus $\chT \times \chT \subset \chG$ as $\chT_1 \times \chT_2$, and then write $\chT_1 \cong \G_{m}^{(1)} \times \G_{m}^{(2)}$ and $\chT_2 \cong \G_m^{(3)} \times \G_m^{(4)}$.  In these terms, the map 
\[
g^{\spec} \co \G_m \times \G_m \times \G_m  \rightarrow \G_{m}^{(1)}  \times \G_{m}^{(2)}  \times \G_{m}^{(3)}  \times \G_{m}^{(4)} 
\]
sends
\[
(a,d,d') \mapsto (a,d, a^{-1}, d').
\]
As the anti-diagonal embedding $\check{\G}_m \inj \check{\G}_m \times \check{\G}_m$ is dual to the pullback along the map $\G_m \times \G_m \rightarrow \G_m$ given by $(x,y) \mapsto x/y$, we have using \eqref{eq: Tate perverse normalization} that 
\begin{equation}\label{eq: RS torus match}
\Dmod(\Bun_{T \times T}) \ni \Delta^{(13)}_! (\ul{k}_{\Bun_1}) \boxtimes \delta_{\triv}^{(2)} \boxtimes \delta_{\triv}^{(4)}\tw{4(g-1)} \xrightarrow{\LL_{\chT \times \chT}}   g_*^{\spec} (\omega_{\Loc_A \times \Loc_1 \times \Loc_1}) \in \QCoh(\Loc_{\chT \times \chT})
\end{equation}
where the superscripts on the LHS indicate the factors in which they lie, and $\Delta^{(13)}$ is the diagonal embedding of $\Bun_1$ into the $\Bun_{\G_m^{(1)}} \times \Bun_{\G_m^{(3)}}$.

In this case Eisenstein compatibility (\S \ref{ssec: Eisenstein compatibility}) says that the diagram 
\[
\begin{tikzcd}[column sep = huge]
\Dmod(\Bun_T \times \Bun_T) \ar[d, "\Eis_{B \times B}"] \ar[r, "\LL_{T \times T}", "\sim"'] & \QCoh(\Loc_{\chT} \times \Loc_{\chT}) \ar[d, "\Eis_{\chB \times \chB}^{\spec}"] &   \\
 \Dmod(\Bun_G)  \ar[r, "\LL_G"] & \IndCoh_{\Nilp}(\Loc_{\chG}) 
\end{tikzcd}
\]
commutes after translating by $(\Omega_C^{1/2} \oplus \Omega_C^{-1/2}) \times (\Omega_C^{1/2} \oplus \Omega_C^{-1/2})  \in \Bun_{\chT} \times \Bun_{\chT}$ and twisting by $\tw{2(g-1) + \deg}$. Combining this with \eqref{eq: RS torus match}, we deduce that
\begin{equation}\label{eq: RS right comparison}
\Dmod(\Bun_G) \ni \Eis_{B \times B} (g_! \ul{k}_{\Bun_1})\tw{6(g-1)+\deg}  \xmapsto{\LL_G} \Eis^{\spec}_{\chB \times \chB} (g_*^{\spec} \omega_{\Loc_A \times \Loc_1 \times \Loc_1})  \in \IndCoh_{\Nilp}(\LocG)
 \end{equation}
 where we have used the notational convention of \S \ref{ssec: singular support}.

\subsubsection{Leftmost terms} Next we examine the terms on the left. In preparation we observe here a binary operation on (co)commutative $\Lambdapos$-graded factorization (co)algebras, which will be seen to be parallel to Construction \ref{constr: diagonal graded algebra}.

\begin{lemmaconstr}\label{constr: fact diagotimes} Let $\Lambda$ be a finite rank free abelian monoid with a chosen basis, as in \S \ref{ssec: graded factorization algebras}. 
 
(1) Let $\cA = \{\cA^{\lambda}\}_{\lambda \in \Lambdapos}$ and $\cB = \{\cB^{\lambda}\}_{\lambda \in \Lambdapos}$ be commutative $\Lambdapos$-graded factorization algebras in $\Shv(\Div^{\Lambdapos} C)$. Then 
\[
\cA \diagotimes \cB := \{ \cA^{\lambda} \sotimes \cB^\lambda\}_{\lambda \in \Lambdapos} \in \Shv(\Div^{\Lambdapos} C)
\]
has a natural structure of commutative $\Lambdapos$-graded factorization algebra.

(2) Let $\cA = \{\cA^{\lambda}\}_{\lambda \in \Lambdapos}$ and $\cB = \{\cB^{\lambda}\}_{\lambda \in \Lambdapos}$ be cocommutative $\Lambdapos$-graded factorization coalgebras in $\Shv(\Div^{\Lambdapos} C)$. Then 
\[
\cA \diagotimes \cB = \{ \cA^{\lambda} \otimes \cB^\lambda\}_{\lambda \in \Lambdapos} \in \Shv(\Div^{\Lambdapos} C)
\]
has a natural structure of cocommutative $\Lambdapos$-graded factorization coalgebra.
\end{lemmaconstr}

\begin{proof}
(1) Define the factorization algebra structure (cf. \eqref{eq: graded factorization algebra})
\begin{equation}
(\cA^{\lambda_1+\lambda_2} \sotimes \cB^{\lambda_1+\lambda_2})|_{(C^{(\lambda_1)} \times C^{(\lambda_2)})_{\mrm{disj}}} \cong (\cA^{\lambda_1} \sotimes \cB^{\lambda_1}) \boxtimes (\cA^{\lambda_2} \sotimes \cB^{\lambda_2})|_{(C^{(\lambda_1)} \times C^{(\lambda_2)})_{\mrm{disj}}}
\end{equation}
as the $\sotimes$ of the respective isomorphisms for $\cA$ and $\cB$. Define the commutative factorization algebra structure 
(cf. \eqref{eq: commutative graded fact algebra}) 
\begin{equation}
(\cA^{\lambda_1}  \sotimes \cB^{\lambda_1}) \boxtimes (\cA^{\lambda_2}\sotimes \cB^{\lambda_2})  \rightarrow \add^!_{\lambda_1, \lambda_2} (\cA^{\lambda_1+\lambda_2} \sotimes \cB^{\lambda_1 + \lambda_2} )
\end{equation}
as the $\sotimes$ of the respective maps for $\cA$ and $\cB$. It is straightforward to check that these define a commutative $\Lambdapos$-graded factorization algebra. 

(2) Similar, with dual formulas. 
\end{proof}

\begin{lemma}\label{lem: diag for factorization algebras}(1) Let $A, B \in \CAlg(\Shv(C), \sotimes)^{\Lambdapos}$. Let $\Fact(A), \Fact(B) \in \CAlg^{\star}(\Shv(\Div^{\Lambdapos} C))$ be the corresponding commutative factorization algebras (under Theorem \ref{thm: commutative graded fact alg}). Then there is a natural isomorphism 
\[
\Fact(A \diagotimes B) \cong \Fact(A) \diagotimes \Fact(B).
\]

(2) Let $A, B \in \CoCAlg(\Shv_{\hol}(C), \otimes)^{\Lambdapos}$. Let $\Fact(A), \Fact(B) \in \CoCAlg^{\star}(\Shv_{\hol}(\Div^{\Lambdapos} C))$ be the corresponding cocommutative factorization coalgebras (under Theorem \ref{thm: cocommutative graded fact coalg}). Then there is a natural isomorphism 
\[
\Fact(A \diagotimes B) \cong \Fact(A) \diagotimes \Fact(B).
\]
\end{lemma}

\begin{proof}
(1) Since $\Fact$ is the inverse of $\Delta^!$ from Theorem \ref{thm: commutative graded fact alg}, an equivalent statement is that for any commutative $\Lambdapos$-graded factorization algebras $\cA, \cB \in \Shv(\Div^{\Lambdapos} C)$, there is a natural isomorphism 
\[
(\Delta^! \cA) \diagotimes (\Delta^! \cB) \cong \Delta^! (\cA  \diagotimes \cB) \in \CAlg(\Shv(C)^{\Lambdapos}, \sotimes)
\]
where the ``$\diagotimes$'' on the left side refers to Construction \ref{constr: diagonal graded algebra} and the ``$\diagotimes$'' on the right side refers to Construction \ref{constr: fact diagotimes}. This follows directly from the combinatorics of the definitions. 

(2) Follows from dualizing (1). 
\end{proof}

\begin{lemma}\label{lem: RS open match}
With notation as in Proposition \ref{prop: sRS Loc} and Proposition \ref{prop: RS-Whittaker}, we have 
\[
\Dmod(\Bun_G) \ni \Hk_{G, \Sym^{\bu}(C)}^{\Std \diagotimes \Std} \star \sW  \tw{6(g-1)} \xmapsto{\LL_G}   \Loc^{\spec} (\Fact(\Sym^{\bu} \Std \diagotimes \Sym^{\bu} \Std)) \otimes \omega_{\LocG} \in \QCoh(\LocG).
\]
\end{lemma}

\begin{proof} The Whittaker normalization \eqref{eq: whittaker normalization} in this case reads 
\begin{equation}\label{eq: RS whit normalization}
\Dmod(\Bun_G) \ni \sW \tw{6(g-1) } \xmapsto{\LL_G} \omega_{\LocG} \in \IndCoh_{\Nilp}(\LocG).
\end{equation}
By Example \ref{ex: Hecke compatibility}, specialized to $r=2$, $\LL_{\GL_2}$ intertwines the tensoring action of 
\[
 \Loc^{\spec}(\Fact(\Sym^{\bu} \Std)) \in \QCoh(\Loc_2) \text{ on } \IndCoh_{\Nilp}(\Loc_2)
 \]
with the Hecke convolution action of 
\[
\prod_{n \geq 0} \ulk \tw{n} \in \Shv(\Hk_{2, \Sym^{\bu} C}^{\Std}) \text{ on }\Dmod(\Bun_2).
\]
Hence by Lemma \ref{lem: diag for factorization algebras}(2), we deduce that $\LL_G$ intertwines the tensoring action of 
\[
 \Loc^{\spec}(\Fact(\Sym^{\bu} \Std \diagotimes \Sym^{\bu} \Std)) \in \QCoh(\Loc_{\chG}) \text{ on } \IndCoh_{\Nilp}(\Loc_{\chG})
 \]
 with the Hecke convolution action of 
\[
\prod_{n \geq 0} \ulk \tw{2n} \in \Shv(\Hk_{G, \Sym^{\bu} C}^{\Std \diagotimes \Std}) \text{ on }\Dmod(\Bun_G).
\]
Acting by these operators on either side of \eqref{eq: RS whit normalization} gives the result.
\end{proof} 

\subsubsection{Conclusion} 
Twisting \eqref{eq: RS triangle 2} by $\tw{8(g-1)}$, we get an exact triangle in $\Shv(\Bun_G)$:
\[
(\Hk_{G, \Sym^{\bu} C}^{\Std \diagotimes \Std}) \star \sW \tw{6(g-1)+\deg}   \rightarrow \pi_{\aut !} (\ul{k}_{\BunX} ) \tw{8(g-1)} \rightarrow \Eis_{B \times B}(g_! \ul{k}_{\Bun_1}) \tw{8(g-1) + 2 \deg}  
\]
where we used that $2r_E +8(g-1) = 2(1-g)+\deg + 8(g-1)= 6(g-1) + \deg$. 

Therefore, we have proved the following theorem. 

\begin{thm}\label{thm: RS} The Geometric Langlands equivalence $\LL_G \co \Dmod(\Bun_G) \rightarrow \IndCoh_{\Nilp}(\LocG) $ intertwines the indicated objects in the diagram below, where each row is exact
\begin{equation}
\begin{tikzcd}[column sep = tiny
]
\Hk_{G, \Sym^{\bu} C}^{\Std \diagotimes \Std} \star \sW  \tw{6(g-1)}  \ar[r] \ar[d, "\LL_G", "\text{Lemma \ref{lem: RS open match}}"', mapsto] &  \cP_{X}  \tw{8(g-1)} \ar[d, "\LL_G?", dashed, mapsto] \ar[r] &  \Eis_{B \times B} (g_! \ul{k}_{\Bun_1}) \tw{6(g-1) +  \deg}\ar[d, "\LL_G", "\eqref{eq: RS right comparison}"', mapsto] \\
\Loc^{\spec} (\Fact(\Sym^{\bu} \Std \diagotimes \Sym^{\bu} \Std))^{\shear} \ar[r] &   \cL_{\chX} \ar[r] &   \Eis^{\spec}_{\check{B} \times \chB}(g_*^{\spec} \omega_{\Loc_A \times \Loc_1 \times \Loc_1})
\end{tikzcd}
\end{equation}
\end{thm}

\bibliographystyle{amsalpha}
\bibliography{Bibliography}

\end{document}